\documentclass[a4paper,10pt]{article} \usepackage{amsmath,amssymb,amsthm,dsfont}
\usepackage[english]{babel}
\newtheorem{theo}{Theorem}[section] 
\newtheorem{defi}[theo]{Definition}
\newtheorem{lemm}[theo]{Lemma} 
\newtheorem{prop}[theo]{Proposition}
\newtheorem{coro}[theo]{Corollary}

\newcommand{\Na}{\mathbb N}                   

\newcommand{\Ra}{\mathbb R}                   

\newcommand{\Ca}{\mathbb C}                   

\newcommand{\scal}[1]{\langle #1 \rangle}

\newcommand{\finpreuve}{\hfill $\Box$}

\newcommand{\name}{$\underline{\qquad \qquad}$}

\begin{document}

\title{\bf \sc Strichartz inequalities on surfaces with cusps}

\author{Jean-Marc Bouclet \\
Institut de Math\'ematiques de Toulouse \\ 
118 route de Narbonne 
\\ 
F-31062 Toulouse Cedex 9 \\
 jean-marc.bouclet@math.univ-toulouse.fr}

\maketitle

{\abstract We prove Strichartz inequalities for the wave and Schr\"odinger equations on noncompact surfaces with ends  of finite area, {\it i.e.} with ends isometric to $ \big(  (r_0,\infty) \times {\mathbb S}^1 , dr^2 + e^{- 2 \phi (r)}d \theta^2 \big) $ with $ e^{-\phi} $ integrable.
We prove first that all Strichartz estimates, with any derivative loss, fail to be true in such ends. We next show for the wave equation that, by projecting off the zero mode of $ {\mathbb S}^1 $, we recover the same inequalities as on $ \Ra^2 $. On the other hand, for the Schr\"odinger equation, we prove that even by projecting off the zero angular modes we have to consider additional losses of derivatives compared to the case of closed surfaces;  in particular, we show that the semiclassical estimates of Burq-G\'erard-Tzvetkov do not hold in such geometries. Moreover our semiclassical estimates with loss are sharp.}

\section{Introduction}
\setcounter{equation}{0}
Strichartz inequalities are well known a priori estimates on linear dispersive partial differential operators which are particularly interesting to solve nonlinear equations at low regularity. Let us recall their usual form for the wave and Schr\"odinger equations on $ \Ra^n $. For $n \geq 2$, if $ (p,q) $ is a  wave admissible pair, namely
\begin{eqnarray}
 p,q \geq 2, \qquad (p,q,n) \ne (2,\infty ,3 ), \qquad \frac{2}{p} + \frac{n-1}{q} \leq \frac{n-1}{2} \label{waveadmissibleiin}
\end{eqnarray}
then the Strichartz inequalities on the solutions  to the wave equation $ \partial_t^2 \Psi - \Delta \Psi = 0 $ are
\begin{eqnarray}
 || \Psi ||_{L^p ( [0,1],L^q (\Ra^n) )} \leq C || \Psi (0) ||_{H^{\sigma_{\rm w}}(\Ra^n)} + C || \partial_t \Psi (0) ||_{H^{\sigma_{\rm w}-1}(\Ra^n)} , \qquad \sigma_{\rm w} = \frac{n}{2} - \frac{n}{q} - \frac{1}{p} . \label{ondeslibres}
\end{eqnarray}
Note that $ \sigma_{\rm w} \geq  \frac{n+1}{2} \left( \frac{1}{2} - \frac{1}{q} \right)  $, with equality for sharp wave admissible pairs, {\it i.e.} when the last inequality in (\ref{waveadmissibleiin}) is an equality.  
 Schr\"odinger admissible pairs are defined by
$$  p, q \geq 2, \qquad (p,q,n) \ne (2,\infty,2) \qquad \frac{2}{p} + \frac{n}{q} = \frac{n}{2} , $$
in any dimension $n \geq 1$, and for such pairs the Strichartz inequalities on solutions to the Schr\"odinger equation $ i \partial_t \Psi + \Delta \Psi = 0 $ are
\begin{eqnarray}
  || \Psi ||_{L^p ( [0,1],L^q (\Ra^n) )} \leq C || \Psi (0) ||_{L^2(\Ra^n)} . \label{Schrodingerlibre}
\end{eqnarray}  
 We refer to \cite{KeelTao} for complete proofs of the above estimates and classical references. 
We recall that the interest of Strichartz inequalities is to guarantee that $ \Psi (t) \in L^q $ for a.e. $t$ (and more precisely in $ L^p $ mean) without using as many derivatives on the initial data as would require  the usual Sobolev estimates 
$$ || \psi ||_{L^q} \leq C || \psi ||_{H^{\frac{n}{2}- \frac{n}{q}}} \qquad (q \in [2,\infty)). $$
The extension of Strichartz inequalities to curved backgrounds has attracted a lot of activity  since many nonlinear dispersive equations are posed on manifolds or domains. In the setting of asymptotically flat or hyperbolic manifolds with non (or weakly \cite{BGH})  trapped  geodesic flow, several  papers have shown that the above estimates still hold (see \cite{Boucletref} for references), including globally in time \cite{MMTa,MetcTata,HaZh,Zhang}.  
Such situations are the most favorable ones since they correspond to {\it large} ends; heuristically, the waves escape to infinity where there is room enough for the dispersion to play in the optimal way. This holds for both the wave and Schr\"odinger equations. In other geometries, the results are as follows. For the wave equation, it is known that Strichartz inequalities are the same as (\ref{ondeslibres}) for smooth enough closed manifolds, or reasonable manifolds with non vanishing injectivity radius (see \cite{Kapi} in the smooth case and \cite{TaIII} for metrics with optimal regularity). In most other cases, one has in general to consider Strichartz inequalities  with losses, meaning that the initial data have to be smoother than what is required in the free cases (\ref{ondeslibres}) or (\ref{Schrodingerlibre}). For the wave equation, this is known for low regularity metrics  \cite{BaCh,TaII} and for manifolds with boundary \cite{ILP}. Furthermore the losses are unavoidable in the sense that there are counterexamples  \cite{SmTa,Ivanovici}. 
For the Schr\"odinger equation, the situation is similar but the losses are more dramatic in compact domains due to the infinite speed of propagation. The general result of \cite{BGT}  says that
\begin{eqnarray}
 \big| \big| \Psi \big| \big|_{L^p([0,1],L^q)} \lesssim || \Psi (0) ||_{H^{1/p}} := || (1-\Delta)^ {\frac{1}{2p}} \Psi_0 ||_{L^2} \label{upBGT} 
\end{eqnarray}
when $ \Delta $ is the Laplace-Beltrami operator on a compact manifold $ ({\mathcal M},G) $. The loss is unavoidable at least on $ {\mathbb S}^3$, though it can be strongly weaken on $ {\mathbb T}^2 $ \cite{Bourgain}. The upper bound (\ref{upBGT})  holds in fairly large generality provided that the injectivity radius of the manifold is positive \cite{BGT}. It also holds for  polygonal domains \cite{BFHM} or manifolds with strictly concave boundaries \cite{OIvan}. For general manifolds with boundary (or low regularity metrics) the losses are worse than $ 1/ p $ \cite{Anton,BSS1} (see also the recent improvement \cite{BSS2} for subadmissible pairs).

Schematically, the usual strategy to address such issues  (for time independent operators) is to prove semiclassical Strichartz inequalities of the form
\begin{eqnarray}
 \big| \big| S (h) e^{it (-\Delta)^{\nu}} \psi  \big| \big|_{L^p ( [0, T(h)],L^q )} \leq C h^{-\sigma} || \psi ||_{L^2} , \label{contientBGT}
\end{eqnarray}
for some spectral localization $ S (h) $ ({\it e.g.} $ S (h) = \varphi (-h^2 \Delta_G) $ with $ \varphi \in C_0^{\infty}(0,+\infty) $) and some suitable time scale $   T (h) $. Here $ \nu = 1 $ for the Schr\"odinger equation and $ \nu = 1/2 $ for the wave equation. In practice $ T(h) $ is dictated by the range of the times  over which one has a good parametrix for the evolution operator (by Fourier integral operators or wave packets); see {\it e.g.} \cite{BaCh,BGT,Anton,BSS1,BSS2} where similar or closely related estimates appear explicitly.
For smooth manifolds without boundary, if we let $ \varrho_{\rm inj} $ be the injectivity radius, one can basically take
$ \sigma = \frac{n+1}{2} $ and $ T (h) \approx \varrho_{\rm inj} $ if $ \nu = 1/2 $, or $ \sigma = 0 $ and $ T (h) \approx h \varrho_{\rm inj}  $ if $ \nu = 0 $.
This leads for instance in \cite{BGT} to the following estimates on closed manifolds
\begin{eqnarray}
  \big| \big| \varphi (-h^2 \Delta) e^{it \Delta}  \psi  \big| \big|_{L^p ( [0,h],L^q ({\mathcal M}) )} \leq C  || \varphi (-h^2 \Delta) \psi ||_{L^2({\mathcal M})} . \label{estimeequonvanier} 
\end{eqnarray}
 Cumulating $  O (1/h) $ such estimates to replace $ [0,h] $ by $ [0,1] $ leads to (\ref{upBGT}). 
If the metric is $ C^s $ with $ 0< s <2 $ or even Lipschitz (a case to which manifolds with boundary can be reduced), one has to 
consider smaller $ T (h) $ and thus to  cumulate more estimates which cause additional losses.


In view of this general picture, it is natural to seek which type of Strichartz inequalities can hold on manifolds with {\it small} ends, where the injectivity radius vanish. In this paper, we will consider the case of surfaces with cusps. They can be thought to as complete\footnote{we will also consider surfaces with boundary, thus non geodesically complete but this won't actually play any role in our result} noncompact surfaces with finite area. An example is $ {\mathcal S} = \Ra_r \times {\mathbb S}^1_{\theta} $ equipped with the metric $ dr^2 + d \theta^2 / \cosh^2 (r)  $.
Our results  are roughly the following ones. The first one is that, due to zero modes on the angular manifold $ {\mathbb S}^1 $, no Strichartz estimate can hold on such surfaces (weighted versions thereof could however hold). This is closely related to the well known fact that even standard Sobolev estimates fail in such geometries. The second result is that, by removing zero angular modes ({\it i.e.} essentially by considering functions with zero mean on $ {\mathbb S}^1 $), the Strichartz inequalities for the wave equation are the same as on $ \Ra^2 $. Thus, in this case, the vanishing of the injectivity radius does not destroy the usual estimates. In other words, the only obstruction to standard inequalities is due to zero angular modes. The situation is more subtle for the Schr\"odinger equation since our third result says that for the Schr\"odinger equation (and after the removal of zero angular modes) we have to consider new losses in the Strichartz inequalities, even at the semiclassical level where they are unavoidable. In this sense, the situation is different from the general one  considered in \cite{BGT}. 



Here are the  precise framework and  results.  Our model for the cusp end is $ \big( {\mathcal S}_0 , G_{0} \big) $ with
\begin{eqnarray}
 {\mathcal S}_0 = [ r_0,+\infty)_r \times {\mathcal A}, \qquad G_{0} = dr^2 + e^{-2 \phi (r)} g_{\mathcal A} , \label{modelwarpedproduct}
\end{eqnarray}
where $ r_0 $ is some real number, $ ( {\mathcal A} , g_{\mathcal A} ) $ is a compact Riemannian manifold of dimension $ 1 $, that is a disjoint union of circles, and $\phi $
is a real valued function such that, 
\begin{eqnarray}
 \int_{r_0}^{+\infty} e^{-\phi(r)} dr < \infty , \label{condvolfini} 
\end{eqnarray}
which means that $ {\mathcal S}_0 $ has finite area (see the Riemannian density in (\ref{seedensity})).  At a more technical level, we will also require that $ \phi $ extends to a smooth  function on $ \Ra $ such that, 
\begin{eqnarray}
|| \phi^{(j)} ||_{L^{\infty}(\Ra)} \leq C_j , \qquad j \geq 1 . \label{order}
\end{eqnarray}
For instance, the functions $ e^{-r} $ and $ \cosh (r)^{-1} $ are of the form $ e^{-\phi(r)} $ with $ \phi $ satisfying (\ref{condvolfini}) and (\ref{order}). We can take any $ r_0 \in \Ra$ in those cases. Other examples are $ \sinh r $ and $ r^{\sigma} $ with $ \sigma > 1 $ , on $ [r_0 , \infty ) $ with $ r_0 > 0 $.
More generally, we will state our main results on surfaces $ \big( {\mathcal S},G \big) $ of the form
\begin{eqnarray}
 {\mathcal S} = {\mathcal K} \sqcup \stackrel{\circ}{{\mathcal S}}_0 \label{surfaceplusgenerale}
\end{eqnarray}
where $ \stackrel{\circ}{{\mathcal S}}_0 = (r_0, \infty) \times {\mathcal A} $ is glued smoothly along $ \{r_0 \} \times{\mathcal A} $ to  a compact surface $ {\mathcal K} $, and where $G$ is a smooth metric on  $ {\mathcal S} $ such that
$$ G = G_{0} \ \ \mbox{on} \ {\mathcal S} \setminus {\mathcal K}. $$
In practice, we shall focus on the analysis on $ {\mathcal S}_0 $ 
but we shall state our main results on $ {\mathcal S} $, seeing $ {\mathcal S}_0 $ as a special case. 

We denote by $ \Delta $ the (non positive) Laplace-Beltrami operator on $ {\mathcal S} $ and by $ d {\rm vol} $ the associated volume density. The same objects on $ {\mathcal S}_0 $ will be denoted with a $ 0 $ index; the Laplacian on $ {\mathcal S}_0 $ is then
\begin{eqnarray}
\Delta_{0} & = & \frac{\partial^2}{\partial r^2} -  \phi^{\prime}(r) \frac{\partial}{\partial r} + e^{2 \phi (r)} \Delta_{\mathcal A} , 
\end{eqnarray}
where $ \Delta_{\mathcal A} $ is the Laplacian on $ {\mathcal A} $, and the volume density is
\begin{eqnarray}
d {\rm vol}_{0} & = &  e^{-\phi(r)} dr d{\mathcal A} ,  \label{seedensity} 
\end{eqnarray}
 where $ d {\mathcal A} $ is the line element on $ {\mathcal A} $. They coincide respectively with $ \Delta $ and $ d {\rm vol} $ on  $ {\mathcal S} \setminus {\mathcal K}  $.
 For $ q \in [ 1 , \infty ] $, we denote
\begin{eqnarray}
 L^q_{G_0} := L^q \big( {\mathcal S}_0 , d {\rm vol}_{0} \big), \qquad L^q_{G} := L^q \big( {\mathcal S} , d {\rm vol} \big) ,
\end{eqnarray}
and will use the shorter notation $ L^q $ for  $ L^q (\Ra) $, $ L^q ((r_0,\infty),dr) $ or  $ L^q \big( (r_0,\infty) \times {\mathcal A} , dr d {\mathcal A} \big) $ {\it i.e.} when the measure is equivalent to the standard Lebesgue measure.
  We will keep the notation $ \Delta $ (resp. $ \Delta_0 $) for the Friedrichs extension of the Laplacian, defined a priori on $ C_0^{\infty} ({\mathcal S} \setminus \partial {\mathcal S}) $ (resp. $C_0^{\infty} \big( (r_0,\infty) \times {\mathcal A} \big)$).   For $ \sigma \in \Ra $ and $ \psi \in \cap_{j \geq 0} \mbox{Dom} (\Delta^j) $, we denote
\begin{eqnarray}  
 || \psi ||_{H^{\sigma}_G} := || (1-\Delta)^{\sigma/2} \psi ||_{L^2_G} , \label{Sobolevnorm}
\end{eqnarray} 
and define the Sobolev space $ H^{\sigma}_G $ as the completion of  $ \cap_{j} \mbox{Dom} (\Delta^j) $ for this norm. Of course, $ H^{\sigma}_{G_0} $ and $ || \cdot ||_{H_{G_0}^{\sigma}} $ are defined analogously on $ {\mathcal S}_0 $.

We denote by $ (e_k)_{k \geq 0} $ an orthonormal basis of $ L^2 ({\mathcal A}, d {\mathcal A}) $ of eigenfunctions of $ \Delta_{\mathcal A} $, with
\begin{eqnarray}
 - \Delta_{\mathcal A} e_k = \mu_k^2 e_k, \qquad 0 \leq \mu_0 \leq \mu_1 \leq \cdots   \label{eigenbasis}
\end{eqnarray}
the  eigenvalues $ \mu_k^2 $ being repeated according to their multiplicities. In particular, we set
\begin{eqnarray}
 k_0 = \mbox{dim Ker} (\Delta_{\mathcal A}) ,  \label{defk0}
\end{eqnarray} 
so that $ \mu_0 = \cdots = \mu_{k_0 - 1} = 0 $ and $ \mu_k \geq \mu_{k_0} > 0 $ for $ k \geq k_0 $. Here $ k_0 $ may be larger than $ 1 $ since we do not assume that $ {\mathcal A} $ is connected. We also define
$$ \pi_0 = \mbox{projection on Ker} (\Delta_{\mathcal A}), \qquad \Pi = 1 \otimes \pi_0, \qquad \Pi^c = I - \Pi $$
where we see $ \Pi $ as an operator on $ L^2_{G_0} $ or on  $ L^2 \big( (r_0,\infty) , e^{-\phi(r)} dr \big) \otimes L^2 ({\mathcal A},d {\mathcal A}) $ which is isomorphic to $ L^2_{G_0} $. It is then  an orthogonal projection. Note that $ \Pi $ is also  an orthogonal projection on the (isomorphic) spaces $ L^2 ((r_0,\infty) \times {\mathcal A},dr d {\mathcal A}) $ and   $ L^2 \big( (r_0,\infty) , dr \big) \otimes L^2 ({\mathcal A},d {\mathcal A}) $.
For clarity, we record that
$$ (\Pi \psi ) (r,\alpha) = \sum_{k < k_0} \left( \int_{\mathcal A} \overline{e_k ( \cdot )} \psi (r, \cdot ) d {\mathcal A} \right) e_{k} (\alpha) , \qquad r > r_0 , \ \alpha \in {\mathcal A} , $$
for functions $ \psi $ on $ {\mathcal S}_0 $ (in $ L^q_{G_0} $ or $ L^q \big( (r_0,\infty) \times {\mathcal A} ,dr d {\mathcal A} \big) $ with $ q \in [1,\infty] $). Note that the dependence on $ \alpha $ is somewhat artificial. If $ {\mathcal A} $ is connected then $ k_0 = 1 $ and $ e_0 $ is a constant function, so that $ \Pi $ is the projection on radial functions and $ \Pi \psi $ is independent of $ \alpha  $. In the general case, $ {\mathcal S}_0 $ has $k_0$ connected components and $ \Pi $ is the sum of projections on radial functions on each component, so $ \alpha $ only labels the component one is looking at.
 
We can now state our main results. The first one says that, due to zero modes on the angular manifold $ {\mathcal A} $, no {\it global} (in space) Sobolev neither Strichartz estimates can hold on cusps.
Regarding the Sobolev estimates, this phenomenon is essentially well known. It is for instance proved in \cite{Hebey} that, on noncompact manifolds of finite volume, the usual Sobolev estimates ({\it i.e.} as on $  \Ra^n$) fail. In Theorem \ref{theorem1} below, we first remark that  we actually never have an embedding of the form $  H_{G}^{\sigma} \subset L^q_G $ for some $q > 2 $ and $ \sigma > 0 $, {\it i.e.} even when $ \sigma $ is large. More originally, we also show that no Strichartz estimates, with any loss, can hold.

\begin{theo}[Zero angular modes destroy Sobolev and Strichartz estimates] \label{theorem1} Let us fix real numbers $ p \geq 1 $, $ q >2 $ and $ \sigma \geq 0 $. 
\begin{enumerate} \item{Sobolev estimates:} there exists a sequence $ (\psi_n )_{n \geq 0} $ of non zero functions in $H^{\sigma}_{G_0} \cap \emph{Ran} ( \Pi ) $ such that
$$ \sup_{n \geq 0} \frac{|| \psi_n||_{L^q_{G_0}}}{|| \psi_n ||_{H^{\sigma}_{G_0}}} = + \infty . $$
\item{Strichartz estimates for the wave equation:} there exists a sequence $ (\psi_n )_{n \geq 0} $ of non zero functions in $H^{\sigma}_{G_0} \cap \emph{Ran} (\Pi) $ such that, if we set
$$ \Psi_n (t) = \cos ( t \sqrt{-\Delta_0} ) \psi_n  , $$
we have
$$ \sup_{n \geq 0} \frac{|| \Psi_n||_{L^p ( [0,1] ; L^q_{G_0} ) }}{|| \psi_n ||_{H^{\sigma}_{G_0}}} = + \infty . $$
\item{Strichartz estimates for the Schr\"odinger equation:} consider the case $ e^{\phi (r)} = e^r $ and $ r_0 = 0 $.  There exists a sequence $ (\psi_n )_{n \geq 0} $ of non zero functions in $H^{\sigma}_{G_0} \cap \emph{Ran} (\Pi) $ such that, if we set
$$ \Psi_n (t) =  e^{-i t \Delta_0} \psi_n  , $$
we have
$$ \sup_{n \geq 0} \frac{|| \Psi_n||_{L^p ( [0,1] ; L^q_{G_0} ) }}{|| \psi_n ||_{H^{\sigma}_{G_0}}} = + \infty . $$
\end{enumerate}
\end{theo} 
 
In contrast to this theorem, we emphasize that Sobolev estimates, and Strichartz estimates likewise, hold however {\it locally} in space as on any Riemannian manifold. One could actually show  {\it weighted} versions of such estimates (at least away from the boundary) with a weight going to zero at infinity.  We are here in the opposite situation to \cite{BanicaDuyckaerts} where working on manifolds with large ends allows to improve the standard Strichartz estimates for radial functions by a growing weight. 

Theorem \ref{theorem1} suggests to look at Strichartz estimates on the range of $ \Pi^c $. To guarantee that $ \Pi^c $ is well defined on a surface $ {\mathcal S} $, we use a spatial cutoff
$$  \mathds{1}_{[r_1,\infty )} (r) = \mbox{  restriction operator to } \ [r_1,\infty ) \times {\mathcal A} \subset {\mathcal S}_0 , \qquad r_1 > r_0 , $$
 to be supported in $ {\mathcal S}_0 $.   For the wave equation, the next theorem states that such a localization allows to recover the same Strichartz estimates as on $ \Ra^2 $ or on closed surfaces.

\begin{theo}[Wave-Strichartz estimates at infinity away from zero angular modes] \label{theoremeonde} Let $ p,q \geq 2 $ be real numbers such that $ (p,q) $ is sharp wave admissible, {\it i.e.}
\begin{eqnarray}
 \frac{2}{p} + \frac{1}{q} = \frac{1}{2} , 
\end{eqnarray}
and set
$$ \sigma_{\rm w} = \frac{3}{2} \left( \frac{1}{2} - \frac{1}{q} \right) . $$
 Then, for any $ r_1 > r_0 $ there exists $ C $ such that, if we set
$$ \Psi (t) = \cos ( t \sqrt{-\Delta} ) \psi_0 + \frac{\sin ( t \sqrt{-\Delta} )}{\sqrt{-\Delta}} \psi_1 , $$
we have
$$ \big| \big| \Pi^c \mathds{1}_{[r_1,\infty )} (r) \Psi  \big| \big|_{L^p ([0,1];L^q_{G_0})} \leq C || \psi_0 ||_{H_G^{\sigma_{\rm w}}} + C || \psi_1 ||_{H_G^{\sigma_{\rm w}-1}} $$ 
for all $ \psi_0, \psi_1 \in \cap_j \emph{Dom}( \Delta^j )$. 
\end{theo} 
 
Note that we consider the $ L^p ([0,1] ; L^q_{G_0} ) $ norm since, thanks to the spatial cutoff, we see $ \mathds{1}_{[r_1,\infty )} (r) \Psi (t) $ as a function on $ {\mathcal S}_0 $ to which we can apply $ \Pi^c $. We shall use this convention everywhere in this paper.
 
We next consider Strichartz estimates for the Schr\"odinger equation. Due to the infinite speed of propagation, we consider both semiclassical and non semiclassical estimates. We shall see here that, even by working away from zero angular modes, we don't recover the general Strichartz estimates of \cite{BGT}, even at the semiclassical level. There is an unavoidable additional derivative loss.

\begin{theo}[Semiclassical Schr\"odinger-Strichartz estimates at infinity away from zero angular modes] \label{semiclassicalSchrodinger} Let $ p,q \geq 2 $ be real numbers such that $ (p,q) $ is  Schr\"odinger admissible, {\it i.e.}
\begin{eqnarray}
 \frac{1}{p} + \frac{1}{q} = \frac{1}{2} , 
\end{eqnarray}
and set
$$ \sigma_{\rm S} = \frac{1}{2} \left( \frac{1}{2} - \frac{1}{q} \right) . $$
Fix $ r_1 > r_0 $ and $ \varphi \in C_0^{\infty} (\Ra) $. Then, there exists $ C $ such that, if we set
\begin{eqnarray}
 \Psi_h (t) = e^{it \Delta} \varphi (- h^2 \Delta) \psi  , \label{notationPsih}
\end{eqnarray}
we have
$$ \big| \big| \Pi^c \mathds{1}_{[r_1,\infty )} (r) \Psi_h  \big| \big|_{L^p ([0,h];L^q_{G_0})} \leq C || \psi ||_{H_G^{\sigma_{\rm S}}} $$ 
for all $ \psi \in L^2_G $ and all $ h \in (0,1] $. 
\end{theo} 

\begin{coro} \label{corollaireSchrodinger} Let $ r_1 > r_0 $ and $ (p,q) \in [2,\infty )^2 $ be  Schr\"odinger admissible. There exists $ C > 0 $ such that, if we set 
\begin{eqnarray}
 \Psi (t) = e^{it \Delta} \psi  , \label{defPsiSchrodinger}
\end{eqnarray}
we have
$$ \big| \big| \Pi^c \mathds{1}_{[r_1,\infty )} (r) \Psi  \big| \big|_{L^p ([0,1];L^q_{G_0})} \leq C || \psi ||_{H_G^{\frac{1}{p} + \sigma_{\rm S}}} $$ 
for all $ \psi \in \cap_j \emph{Dom}(\Delta^j) $. 
\end{coro}
Notice that the loss $ \sigma_{\rm S} + \frac{1}{p} = \frac{3}{2} \big(\frac{1}{2} - \frac{1}{q} \big) $ is larger than the general one obtained in \cite{BGT}, but it remains better than the Sobolev index $ 2 \big( \frac{1}{2} - \frac{1}{q} \big) $ in two dimensions. We do not know whether these non semiclassical Strichartz estimates are sharp. However, at the semiclassical level, Theorem \ref{theorem5} below shows that the estimates of Theorem \ref{semiclassicalSchrodinger} are sharp, which is already a difference with the general situation of \cite{BGT} and suggests that the estimates of Corollary \ref{corollaireSchrodinger} may be natural. 

\begin{theo} \label{theorem5} Let $ e^{\phi(r)} (r) = e^r $, $ r_0 = 0 $ and fix $ r_1 > 0  $.  We can find $ \varphi \in C_0^{\infty} (\Ra) $ and a family $ (\psi_{0}^h)_{h \in (0,h_0]} $ of non zero functions in $ C_0^{\infty} \big( (r_0,\infty) \times {\mathcal A} \big)  $ such that, if we set
$$ \Psi_h (t) = e^{it \Delta} \varphi (-h^2 \Delta) \psi_{0}^h  , $$
then, for any  sharp Schr\"odinger admissible  pair $ (p,q) $ (with $q > 2$) and any $ \sigma < \sigma_{\rm S} $, we have
$$ \lim_{h \rightarrow 0} \frac{ \big| \big| \Pi^c \mathds{1}_{[r_1,\infty )} (r) \Psi_h  \big| \big|_{L^p ([0,h];L^q_{G_0})}}{|| \psi_0^h ||_{H^{\sigma}_{G_0}}} = + \infty . $$
\end{theo}
 
The plan of the paper is as follows. In Section \ref{section2}, we explain how to separate variables and prove some useful elliptic estimates. In Section \ref{section3}, we provide a suitable pseudo-differential description of $ \varphi (-h^2 \Delta_G) $ which we use in particular in Section \ref{section4} to derive a  Littlewood-Paley decomposition. Theorems \ref{theoremeonde}, \ref{semiclassicalSchrodinger} and Corollary \ref{corollaireSchrodinger} are proved in Section \ref{proofs} while the counterexamples, {\it i.e.} Theorems \ref{theorem1} and \ref{theorem5}, are proved in Section \ref{modezero}.

\medskip

\noindent {\bf Acknowledgments.} It is a pleasure to thank Maciej Zworski for suggesting the study of Strichartz estimates on manifolds with cusps. We also thank Semyon Dyatlov for helpful discussions about cusps.
 
\section{Separation of variables and resolvent estimates} \label{section2} 
\setcounter{equation}{0}
\subsection{Separation of variables}
This paragraph is devoted to the model warped product (\ref{modelwarpedproduct}) for which we describe basic objects, mostly for further notational purposes. We explain in particular how to separate variables. 
We consider the unitary mapping 
\begin{eqnarray}
 {\mathcal U} : L^2_{G_0} \rightarrow L^2 ,
\qquad {\mathcal U} \psi = e^{- \frac{\phi(r)}{2}} \psi , \label{theunitarymapping}
\end{eqnarray}
where $ L^2 =L^2 \big( (r_0,\infty),dr \big) \otimes L^2 ({\mathcal A}, d {\mathcal A}) $ 
and we let
\begin{eqnarray}
 P := {\mathcal U} (-\Delta_0) {\mathcal U}^* & = & - \partial_r^2 - e^{2 \phi (r)} \Delta_{\mathcal A} + 
\frac{\phi^{\prime}(r)^2 - 2 \phi^{\prime \prime}(r)}{4}  \label{pourunitaritereference} \\
& = : &  - \partial_r^2 - e^{2 \phi (r)} \Delta_{\mathcal A} + w (r) . \label{referenceaw}
\end{eqnarray}
This defines in passing the potential $w$ which is bounded as well as its derivatives by (\ref{order}).

 For both $ \Delta_0 $ and $ P $, we consider the Dirichlet boundary condition at $ r_0 $, namely their Friedrichs extension from $ C_0^{\infty} \big( (r_0,\infty) \times {\mathcal A} \big) $ on $L^2_{G_0}$ and $ L^2 $ respectively. The domain of $ P $ can be described as follows. We introduce the space
$$ {\mathcal H}_0^1 = \mbox{closure of} \ C_0^{\infty}((r_0,\infty) \times {\mathcal A}) \ \mbox{for} \ \ \left( || \partial_r u ||_{L^2}^2 + \big| \big| e^{\phi(r)} |\Delta_{\mathcal A}|^{1/2} u  \big| \big|_{L^2}^2 +  (u , (w(r) + c_0)u)_{L^2} \right)^{1/2} $$
with a constant $ c_0 > 0 $ such that $ - || w||_{L^{\infty}} + c_0 \geq 1 $. 
 We also consider the sesquilinear form
$$ Q (u,v) = (\partial_r u , \partial v)_{L^2} + \big( e^{\phi(r)} |\Delta_{\mathcal A}|^{1/2} u , e^{\phi(r)} |\Delta_{\mathcal A}|^{1/2} v \big)_{L^2} + (u,w(r)v)_{L^2}, \qquad u , v \in {\mathcal H}_0^1 , $$
where $ \partial_r $ and $ e^{\phi(r)} |\Delta_{\mathcal A}|^{1/2} $ are  the continuous extensions of those operators from $ C_0^{\infty} $ to $ {\mathcal H}_0^1 $. Then, 
$$ \mbox{Dom}(P) = \{ u \in {\mathcal H}_0^1 \ | \ |Q(u,v)| \leq C_u || v ||_{L^2} \ \mbox{for all} \ v \in {\mathcal H}_0^1 \} , $$
and
$$   {\mathcal U} \big( \mbox{Dom} (\Delta_0) \big) = \mbox{Dom}(P) . $$
Notice that, since $ P $ is unitarily equivalent to $ - \Delta_0 $ which is non-negative, we have
\begin{eqnarray}
 P \geq 0 . \label{Pnonnegatif}
\end{eqnarray}
To describe the separation of variables,  we introduce for any integer $ k \geq 0 $ the sesquilinear form
$$ q_k ( f , g ) = \int_{r_0}^{+\infty} \overline{f^{\prime}} g^{\prime} + \big( \mu_k^2 e^{2 \phi(r)} + w (r) \big) \overline{f} g dr $$
first for $ f , g \in C_0^{\infty}  (r_0,\infty )  $ and then in $ {\mathfrak h}_{0,k}^1 $ with
$$ {\mathfrak h}_{0,k}^1 := \mbox{closure of} \ C_0^{\infty} (r_0,\infty)  \ \mbox{for the norm} \ \ \left( q_k (f,f) + c_0 || f||_{L^2(r_0,\infty)}^2 \right)^{1/2} . $$
The choice of $ c_0 $ implies that
\begin{eqnarray}
 || f ||_{H_0^1 (r_0,\infty)} \leq || f ||_{{\mathfrak h}_{0,k}^1} , \label{injectionSobolevexplicite}
\end{eqnarray}
and in particular that $ H_0^1 (r_0,\infty) \subset {\mathfrak h}_{0,k}^1 $.
We also consider the related one dimensional Schr\"odinger operator
\begin{eqnarray}
 {\mathfrak p}_k = - \partial_r^2 + \mu_k^2 e^{2 \phi (r)} + w (r) , \label{defpk}
\end{eqnarray} 
which is selfadjoint on $ L^2 (r_0,\infty) $ if we define it on the domain
$$ \mbox{Dom} ({\mathfrak p}_k) = \left\{ f \in \ {\mathfrak h}_{0,k}^1 | \ |q_k (f,g)| \leq C_f || g ||_{L^2(r_0,\infty)} \ \ \mbox{for all} \ g \in {\mathfrak h}_{0,k}^1 \right\} .$$
 Using the eigenbasis of $ \Delta_{\mathcal A} $ (see (\ref{eigenbasis})), any $u$ in $L^2$ can be decomposed as
\begin{eqnarray}
 u (r,\alpha) = \sum_{k \geq 0} u_k (r) e_k (\alpha) , \qquad u_k (r) = \int_{ {\mathcal A}} \overline{e_k(\alpha)} u (r,\alpha) d {\mathcal A}  , \label{Pythagore}
\end{eqnarray}
where the sum converges in $ L^2 $. All this leads to the following

\begin{prop}[Separation of variables] \label{separationdevariablesprop}  The map 
\begin{eqnarray}
  L^2 \big( (r_0,\infty) \times {\mathcal A} , dr d {\mathcal A}\big) \ni u \mapsto (u_k)_{k \geq 0}  \in \bigoplus_{k \geq 0} L^2 (r_0,\infty) \label{isometryl2}
\end{eqnarray}  
is an isometry. In particular
$$ || u ||_{L^2}^2 = \sum_{k \geq 0} || u_k ||_{L^2 (r_0,\infty)}^2 . $$
We also have the following characterizations
\begin{enumerate}
\item{$$ u \in {\mathcal H}_0^1 \qquad \Longleftrightarrow \qquad u_k \in {\mathfrak h}_{0,k}^1 \ \ \mbox{for all} \ k \ \ \mbox{and} \ \ \ \sum_{k \geq 0} q_k (u_k , u_k) + c_0 || u_k ||^2_{L^2(r_0,\infty)} < \infty $$
in which case the last sum equals $ || u ||_{{\mathcal H}_0^1}^2 $,}
\item{ $$ u \in \emph{Dom}(P) \qquad \Longleftrightarrow \qquad u_k \in \emph{Dom}({\mathfrak p}_k) \ \ \mbox{for all} \ k \geq 0 \ \ \mbox{and} \ \sum_{k \geq 0} || {\mathfrak p}_k u_k ||^2_{L^2 (r_0,\infty)} < \infty , $$
in which case we have
$$ || P u ||_{L^2}^2 = \sum_{k \geq 0 } || {\mathfrak p}_k u_k||^2_{L^2(r_0,\infty)} , \qquad (P u)_k = {\mathfrak p}_k u_k . $$ }
\end{enumerate}
\end{prop}

The isometry (\ref{isometryl2}) allows to constuct operators on $ L^2_{G_0} $ from a sequence of operators on $ L^ 2 (r_0,\infty) $: if $ (A_k)_{k \geq 0} $ is a bounded sequence of bounded operators on $ L^2 (r_0, \infty) $, then we can define $ A $ on $ L^2 \big( (r_0,\infty) \times {\mathcal A} , d r d {\mathcal A} \big)  $ by
\begin{eqnarray}
 A u := \sum_{k \geq 0} (A_k u_k ) \otimes e_k , \label{pourbornoppenheimer}
\end{eqnarray} 
and we have
\begin{eqnarray}
 ||A||_{L^2 \rightarrow L^2} = \sup_k || A_k ||_{L^2(r_0,\infty) \rightarrow L^2(r_0,\infty)} . \label{bornebornoppenheimer}
\end{eqnarray} 
The associated operator on $ L^2_{G_0} $ will be $ e^{\phi(r)/2} A e^{-\phi(r)/2} $. If $ \varphi $ is a bounded Borel function, we record that 
\begin{eqnarray}
\varphi(-\Delta_0) = e^{\frac{\phi(r)}{2}} \varphi (P) e^{-\frac{\phi(r)}{2}}, \qquad \varphi (P) u = \sum_{k \geq 0} \varphi ({\mathfrak p}_k)u_k \otimes e_k , \label{opdiaglapbis}
\end{eqnarray}
the expression of $ \varphi (P) $ following from the unitary equivalence of $ P$ with the sequence $ ({\mathfrak p}_k)_{k \geq 0} $ through (\ref{isometryl2}).

\subsection{Resolvent estimates}
The purpose of this paragraph is to prove $ L^2$ elliptic a priori estimates  away from the zero angular modes. 
Everywhere, we consider a smooth function $ \xi = \xi (r) $ such that
\begin{eqnarray}
 \mbox{supp}(\xi) \subset [ r_1 , \infty ) \subset (r_0,\infty ), \qquad \xi (r) \equiv 1 \ \ r \gg 1 . \label{introductiondexi}
\end{eqnarray} 
It will be used on $ {\mathcal S}_0 $ as a cutoff away from its boundary. More generally, if $ {\mathcal S} $ is a surface as (\ref{surfaceplusgenerale}), $ \xi $ will serve as a localisation in the interior of $ {\mathcal S}_0 $ so that the operator $ \Pi^c \xi : L^2_{G} \rightarrow L^2_{G_0}$ is well defined on $ {\mathcal S} $.
 
\begin{prop} \label{theoremesurresolvante} 
 For all integers $ N \geq 1 $ and  $ N_1 , N_2 $ such that
\begin{eqnarray}
 N_1 + 2 N_2 \leq 2 N , \label{conditionN} 
\end{eqnarray} 
there exists $ C > 0 $ such that
$$ \left| \left| \Delta_{\mathcal A}^{N_2} e^{2 N_2 \phi (r)} D_r^{N_1} \Pi^c \xi  (P+1)^{-N}  u \right| \right|_{L^2  } \leq C  ||  u ||_{ L^2  } , $$
for all  $ u \in C_0^{\infty} ((r_0,\infty) \times {\mathcal A})  $. Here $ L^2 = L^2 ( (r_0,\infty) \times {\mathcal A} , dr d {\mathcal A} ) $.
\end{prop}

Note that we can get rid of the factor $ \Delta_{\mathcal A}^{N_2} $ (which is the usual shorthand for $ 1 \otimes \Delta_{\mathcal A}^{N_2} $) since it is invertible on the range of $ \Pi^c $ and since $ e^{2 N_2 \phi(r)} D_r^{N_1}  $ commutes with $ \Pi^c $.

The estimates of Proposition \ref{theoremesurresolvante} are elliptic in the usual sense of smoothness, but also in the spatial sense since (\ref{condvolfini}) and (\ref{order}) imply that $ e^{\phi(r)} \rightarrow \infty $ as $ r \rightarrow \infty $ according to the following proposition.

\begin{prop} \label{propositionlemme} We have $  \lim_{r \rightarrow + \infty} e^{\phi(r)} = + \infty $ and $ \sum_{L \geq L_0 } e^{-\phi (L)} < \infty $  for any integer $ L_0 > r_0  $.
\end{prop} 

\noindent {\it Proof.} We observe that (\ref{order}) for $j=1$ implies that there exists $ C > 0 $ such that
\begin{eqnarray}
C^{-1} e^{\phi(L)} \leq  e^{\phi (r)} \leq C e^{\phi(L)} , \qquad |r-L| \leq 1 . \label{encadrement}
\end{eqnarray}
 We have in particular
$$ \int_{L_0}^{\infty} e^{-\phi(r)}dr \geq C^{-1 }\sum_{L \geq L_0} e^{-\phi(L)} . $$
By (\ref{condvolfini}), the sum is finite. This  implies that $ e^{-\phi(L)} \rightarrow 0 $ hence that $ e^{\phi(r)} \rightarrow + \infty $ by (\ref{encadrement}). \finpreuve

\bigskip

Before proving Proposition \ref{theoremesurresolvante}, we proceed to a few reductions. We first observe that the estimate  is only non obvious where $ r $ is large, otherwise it is a simple consequence of standard local elliptic regularity. We may thus assume in the proof that 
\begin{eqnarray}
  \mbox{supp}(\xi) \subset [ r_1 ,\infty ) \ \ \mbox{with} \ \  r_1 \geq L_0 \ \ \mbox{large enough (to be chosen)} \ .  \label{fixcutoff}
\end{eqnarray}
By separation of variables, it then suffices to show that for all $ k \geq k_0 $ (see (\ref{defk0}))
\begin{eqnarray}
\left| \left| \mu_k^{2 N_2} e^{2 N_2 \phi (r)} D_r^{N_1} \xi ( {\mathfrak p}_k + 1)^{-N}  f \right| \right|_{L^2 (L_0,\infty) } \leq C  ||  f ||_{ L^2 (r_0,\infty) }, \qquad f \in C_0^{\infty} (r_0,\infty),
\label{resolvante1Dfacile}
\end{eqnarray}
with a constant $ C $ independent of $ k  $. 

To deal easily with the possible exponential growth of $ e^{\phi(r)} $ by mean of standard microlocal methods, we shall reduce this problem to a family of problems localized in spatial shells where $ r \sim L \geq L_0 $. For this purpose, we will use the following two simple estimates
\begin{eqnarray}
 \sum_{L \geq L_0} || f ||^2_{L^2 (L-2,L+2)} \leq 4 || f ||_{L^2 (L_0-2,\infty)}^2 , \label{donneorthogonale}
\end{eqnarray}
and  
\begin{eqnarray}
 || v ||_{L^2(L_0,+\infty)}^2  \leq \sum_{L \geq L_0} || v ||^2_{L^2 (L - 1,L+1)} . \label{solutionorthogonale}
\end{eqnarray}
Let us proceed to the detailed analysis. We consider the semiclassical parameter 
\begin{eqnarray}
 \epsilon = e^{-\phi(L)} \mu_k^{-1}   , \label{parametresemiclassique}
\end{eqnarray} 
and write
\begin{eqnarray}
\mathfrak{p}_k = \epsilon^{-2} \big( \epsilon^2 D_r^2 + V_{k,L} (r) \big), \qquad V_{k,L}(r) := e^{2 (\phi(r) - \phi(L))} + \epsilon^2 w(r) .  
\label{definitPLH}
\end{eqnarray}

\begin{lemm} \label{pourprolongerpotentiel} We can choose $ L_0 > r_0 $ large enough and a real number $ m > 0 $ such that, for all $ k \geq k_0 $ and all $ L \geq L_0 $,
\begin{eqnarray}
V_{k,L} (r) \geq m, \qquad r \in [L-3,L+3] . \label{bornefinfV}
\end{eqnarray}
Moreover, for all $ \alpha \geq 0 $,
\begin{eqnarray}
| \partial_r^{\alpha} V_{k,L} (r) | \leq C_{\alpha}, \qquad r \in [L-3,L+3] , \label{bornesupV}
\end{eqnarray}
with $ C_{\alpha} $ independent of $ L $ and $ k $.
\end{lemm}

\noindent {\it Proof.} By (\ref{order}), $ \phi (r)-\phi (L) $ is bounded if $ |r-L| $ is bounded, so there exists $ m > 0 $ such that
$$ e^{2 (\phi(r) - \phi(L))} \geq 2 m, \qquad L > r_0, \ k \geq k_0, \ r \in [L-3,L+3] .$$
 Using that $w$ is bounded (by (\ref{order}) too), we have $ \epsilon^2 || w ||_{L^{\infty}} \leq m $ if $ \epsilon $ is small enough. Since $ e^{-\phi(L)} $ goes to zero as $ L $ goes to infinity, we can choose $ L_0 $ large enough (independent of $k$) to guarantee that $ \epsilon $ is small enough. This shows (\ref{bornefinfV}). The estimates (\ref{bornesupV}) follow easily from (\ref{order}). \finpreuve

\bigskip

We fix $ \chi , \tilde{\chi} , \tilde{\tilde{\chi}} \in C_0^{\infty} (-2,2) $ such that
$$ \chi \equiv 1 \ \ \mbox{near} \ [-1,1], \qquad \tilde{\chi} \equiv 1 \ \ \mbox{near supp}(\chi), \qquad \tilde{\tilde{\chi}} \equiv 1 \ \ \mbox{near supp}(\tilde{\chi}), $$
and define
$$ \chi_L (r) = \chi (r-L), \qquad \tilde{\chi}_L (r) = \tilde{\chi}(r-L), \qquad \tilde{\tilde{\chi}}_L(r) = \tilde{\tilde{\chi}}(r-L) . $$
In the next proposition, we use the usual semiclassical quantization
\begin{eqnarray}
O \! p_{\epsilon} (a) f = (2\pi)^{-1} \int e^{i r \rho} a (x, \epsilon \rho) \hat{f}(\rho)d \rho , \label{quantificationusuelle}
\end{eqnarray}
where $ \hat{f}(\rho) = \int e^{-ir\rho}f(r)dr $.
\begin{prop}[Spatially localized parametrix] Fix $ M \geq 0 $. There are two families of symbols  
$$ a_{k,L} \in S^{-2N}(\Ra^2) , \qquad r_{k,L} \in S^{-M} (\Ra^2) , $$
bounded in $ S^{-2N} $ and $ S^{-M} $ respectively as $ k \geq k_0 $ and $ L \geq L_0 $ vary, such that
\begin{eqnarray}
 \chi_{L} ( \mathfrak{p}_{k} + 1 )^{-N} = \epsilon^{2N} \chi_{L} O \! p_{\epsilon}\big( a_{k,L} \big)  \tilde{\chi}_{L} + \epsilon^M \chi_{L} O \! p_{\epsilon} \big( r_{k,L} \big) \tilde{\tilde{\chi}}_L ( \mathfrak{p}_{k} + 1 )^{-N} . \label{resolvente1} 
\end{eqnarray} 
\end{prop}

\noindent {\it Proof.}
Choose $  \chi_0 \in C_0^{\infty}(-3,3) $ with values in $ [0,1] $, such that $ \chi_0 \equiv 1 $ near $ [-2,2]  $. Define the potential
$$ \widetilde{V}_{k,L} (r) = \chi_0 (r-L) V_{k,L} (r) + m (1 - \chi_0 ) (r-L) $$
and the semiclassical operator on $ \Ra $
$$ \widetilde{\mathfrak{p}}_{k,L} = \epsilon^2 D_r^2 + \widetilde{V}_{k,L} (r) . $$
 By construction,  $ \widetilde{\mathfrak{p}}_{k,L} $ coincides with $ \epsilon^2 \mathfrak{p}_k $ on $ [L-2,L+2] $. By Lemma \ref{pourprolongerpotentiel}, its symbols satisfies
$$ \rho^2 + \widetilde{V}_{k,L} (r) \geq \rho^2 + m , \qquad (r,\rho) \in \Ra^2, \ k \geq k_0, \ L \geq L_0 , $$
and belongs to a bounded family in $ S^2 (\Ra^2) $ as $ k $ and $ L  $ vary. By standard elliptic parametrix construction, the above lower bound ensures that, if we fix a compact subset $ K $ of $ \Ca $ such that $ \mbox{dist}(K,[m,+\infty)) > 0 $, we can find  symbols $ \tilde{q}_{k,L,z} \in S^{-2N} $ and $ \tilde{r}_{k,L,z} \in S^{-M} $ such that, for all $z \in K$,
\begin{eqnarray}
 (\widetilde{\mathfrak{p}}_{\epsilon,L} - z )^N O \! p_{\epsilon} ( \tilde{q}_{k,L,z} ) = 1 + \epsilon^M  O \! p_{\epsilon} ( \tilde{r}_{k,L,z} ) , 
 \label{parametrixeboitenoire}
\end{eqnarray} 
 with  $ \tilde{q}_{k,L,z} $ (resp. $ \tilde{r}_{k,L,z} $) in a bounded subset of $ S^{-2N} $ (resp. $ S^{-M} $)  when $ k \geq k_0 $, $ L \geq L_0 $ and $ z \in K $ vary. By possibly increasing $ L_0 $ to make $ \epsilon $ small enough, we may take
 $$ z = - \epsilon^2  , $$
 which is in a neighborhood of $0$ hence at positive distance from $ [m,\infty) $. 
Then, using (\ref{parametrixeboitenoire}) and that $ \widetilde{\mathfrak{p}}_{k,L} = \epsilon^2 \mathfrak{p}_{k} $ on $ \mbox{supp} (\tilde{\chi}_L) $, we get
$$ \epsilon^{2N} ( \mathfrak{p}_{k} +1 )^N \tilde{\chi}_{L} O \! p_{\epsilon} ( q_{k,L,-\epsilon^2} )  \chi_{L} = \chi_{L} + \epsilon^M   O \! p_{\epsilon} ( \tilde{r}_{k,L,-\epsilon^2} ) \chi_{L} + \epsilon^{2N} \big[ ( \mathfrak{p}_{k} +1 )^N , \tilde{\chi}_{L} \big]  O \! p_{\epsilon} ( q_{k,L,-\epsilon^2} )  \chi_{L} . $$
The last term is smoothing and $ O (\epsilon^{\infty}) $ since the commutator and $ \chi_L $ have disjoint supports. Then, by multiplying to the left by $ \tilde{\tilde{\chi}}_L $, taking the adjoint and applying $ ( \mathfrak{p}_{k} +1 )^{-N} $ to the right of the resulting identity, we get (\ref{resolvente1}). \finpreuve
 
\bigskip 

\noindent {\it Proof of Proposition \ref{theoremesurresolvante}.} It suffices to show that
\begin{eqnarray}
 \left| \left| \mu_k^{2 N_2} e^{2 N_2 \phi (r)} D_r^{N_1} ({\mathfrak p}_k+1)^{-N} f \right| \right|_{L^2 ( L-1,L+1 ) } \leq C \left( || f ||_{L^2 ( L-2,L+2 ) } + e^{-\phi(L)} || f ||_{L^2 (r_0,\infty)} \right) ,  \label{reductionshell}
\end{eqnarray}
with a constant $ C $ independent of $  k \geq k_0 $, $ L \geq L_0 $ and $ f \in C_0^{\infty} (r_0,\infty) $. Indeed, with (\ref{reductionshell}) at hand, (\ref{resolvante1Dfacile}) follows easily from (\ref{donneorthogonale}), (\ref{solutionorthogonale}) and the summability of $ e^{-\phi(L)} $ given by Proposition \ref{propositionlemme}. Since $ \chi_{L}=1 $ on $ (L-1,L+1) $, it suffices to estimate
$$ \left| \left| \mu_k^{2 N_2} e^{2 N_2 \phi (r)}  D_r^{N_1} \chi_{L} ({\mathfrak p}_k+1)^{-N} f \right| \right|_{L^2 ( \Ra ) } , $$
by the right hand side of (\ref{reductionshell}).  Using (\ref{parametresemiclassique}),
\begin{eqnarray}
 \mu_k^{2 N_2} e^{2 N_2  \phi (r)} D_r^{N_1} = \epsilon^{-2N_2 - N_1} e^{2N_2 (\phi (r) - \phi (L))} (\epsilon D_r)^{N_1}, \label{semiclassiquereecrit}
\end{eqnarray} 
where  $ \phi (r) - \phi (L) $ is bounded on the support of $ \chi_{L} $ (uniformly in $ L \geq L_0 $). We then write the resolvent using (\ref{resolvente1}). 
By (\ref{conditionN}), (\ref{semiclassiquereecrit}) and the standard $ L^2 $ boundedness of  pseudodifferential operators, we have
\begin{eqnarray}
 \left| \left| \mu_k^{2 N_2} e^{2 N_2  \phi (r)} D_r^{N_1}  \epsilon^{2N} \chi_{L} O \! p_{\epsilon} ( a_{k,L}) \tilde{\chi}_{L} f \right| \right|_{L^2 (\Ra) } & \leq & C || \tilde{\chi}_{L} f ||_{L^2 (\Ra)}  \nonumber \\
& \leq & C || f ||_{L^2 (L- 2, L + 2)}  \label{estimeesurdiagexp}
\end{eqnarray}
 with a constant $ C $ independent of $ k $ and $L$. Similarly, by choosing $ M $ such that
 $$ M \geq 2 N + 1 , $$
 we obtain 
 \begin{eqnarray}
  \left| \left|  \mu_k^{2 N_2} e^{2 N_2  \phi (r)} D_r^{N_1} \chi_{L} \epsilon^M O \! p_{\epsilon}( r_{k,L}) \tilde{\tilde{\chi}}_L ({\mathfrak p}_k+1)^{-N} f \right| \right|_{L^2 ( \Ra ) } & \leq & C  \epsilon ||  f ||_{L^2 (r_0,\infty)}  \nonumber \\
 & \leq &  C e^{- \phi (L)} ||f ||_{L^2 (r_0,\infty)} . \label{estimeehorsdiagexp}
\end{eqnarray}
Here we have also used the bound
$ || ( {\mathfrak p}_k + 1 )^{-N} ||_{L^2 (r_0,\infty) \rightarrow L^2 (r_0,\infty)} \leq 1 $.
Using (\ref{estimeesurdiagexp}) and (\ref{estimeehorsdiagexp}), we obtain (\ref{reductionshell}).  \finpreuve

\bigskip

In the next proposition, we convert the result of Proposition \ref{theoremesurresolvante} for $ P $ into estimates for the Laplacian $ \Delta $ on a surface $ {\mathcal S} $ as in (\ref{surfaceplusgenerale}) (hence in particular for $ {\mathcal S}_0 $ itself).

\begin{prop} \label{aprioriexp}
For all integers $ N , N_1 , N_2 \geq 0 $ such that 
$$ 2 N_2 + N_1  \leq 2 N , $$
there exists $ C > 0 $ such that
$$ \left| \left| \Delta_{\mathcal A}^{N_2} e^{2 N_2 \phi (r)} D_r^{N_1} \Pi^c \xi (1-\Delta)^{-N}  \psi \right| \right|_{L^2_{G_0} } \leq C  ||  \psi ||_{ L^2_G } , $$
for all $ \psi \in C_0^{\infty}({\mathcal S}) $. 
\end{prop}

We use the norm $ L^2_{G_0} $ in the left hand side to emphasize that we consider a function supported in $ {\mathcal S}_0 $, although $ (1-\Delta)^{-N} \psi$ belongs to $ L^2_G $.

\bigskip

\noindent {\it Proof.} Define $ u = {\mathcal U} \Pi^c \xi (1-\Delta)^{-N} \psi $. Then, using that $ \Delta = \Delta_{0} $ near the support of $ \xi $, that $ P {\mathcal U} = {\mathcal U} (-\Delta_{0}) $, and that $ \Pi^c  $ commutes with $ {\mathcal U} $ and $ \Delta_{0} $,
$$ (P +1)^N u = {\mathcal U} \Pi^c \xi \psi + {\mathcal U} \Pi_c \big[ (1-\Delta)^N, \xi \big] (1-\Delta)^{-N} \psi .  $$
The commutator is a differential operator of order $ 2N-1 $ with compactly supported coefficients so, by standard (local) elliptic regularity, $ \big[ (1-\Delta)^N, \xi \big] (1-\Delta)^{-N} $ is bounded on $ L^2_{G}  $. 
Then, using that $ \Pi^c $ is a projection which commutes with $ P $ (hence with $ (P+1)^{-N} $),
$$ u = \Pi^c (P+1)^{-N} {\mathcal U} \left( \xi \psi + \big[ (\Delta_G +1)^N, \xi \big] (\Delta_G+1)^{-N} \psi \right) .  $$
Since $ u $ is supported on $ \mbox{supp}(\xi) $, we can multiply both sides of the above equality by $ \tilde{\xi}(r) $ for some smooth function $ \tilde{\xi} $ supported in $  (r_0,\infty ) $ and equal to $1$ near $ \mbox{supp} (\xi) $. We then conclude by using Proposition \ref{theoremesurresolvante} (with $ \tilde{\xi} $) together with the fact that the operator $ \Delta_{\mathcal A}^{N_2} e^{2 N_2 \phi(r)} D_r^{N_1} {\mathcal U}^* $ (which we want to apply to $u$) is a linear combination of operators of the form
$$ {\mathcal U}^* \phi^{(m_1)}(r) \cdots  \phi^{(m_M)}(r) \Delta_{\mathcal A}^{N_2} e^{2 N_2 \phi (r)}  D_r^n , \qquad n + m_1 + \cdots + m_M = N , $$ 
and where $ m_1, \ldots, m_N \geq 1 $ so that each  $ \phi^{(m_j)} $ is bounded by (\ref{order}). The result follows. \finpreuve

\bigskip

We end up this subsection with the following rough Sobolev estimates. 
\begin{prop} \label{aprioriL1} There exists $ N_0 \in \Na $ such that, for all $ q \in [2, \infty] $ and all $ N \geq N_0 $, 
$$ \left| \left| e^{2 N \phi (r)}\Pi^c (1- \Delta_0)^{-2 N} \right| \right|_{ L^2_{G_0} \rightarrow L^q_{G_0} } < \infty . $$ In particular,
$ \big| \big|  (1-\Delta_0)^{-2 N_0}  \Pi^c\big| \big|_{L^1_{G_0} \rightarrow L^2_{G_0}} < \infty $. On $ {\mathcal S} $, we have
\begin{eqnarray}
 \left| \left| e^{2 N \phi (r)}\Pi^c \xi (1- \Delta)^{-2 N} \right| \right|_{ L^2_{G} \rightarrow L^q_{G_0} } < \infty . \label{onS}
\end{eqnarray}
\end{prop}

\noindent {\it Proof.} We only prove the result for $ \Delta_0 $ since it implies (\ref{onS}) by the same trick as in the proof of Proposition \ref{aprioriexp}. It suffices to prove the result for $q=\infty$. The other $q$ are treated by interpolation and the $ L^1_{G_0} \rightarrow L^2_{G_0} $ boundedness follows by taking the adjoint.  By (\ref{pourunitaritereference}), the problem is equivalent to prove the  $ L^2 \rightarrow L^{\infty} $ boundedness of $ e^{(2N + 1/2) \phi (r)}\Pi^c (P+1)^{-2N} $. Let $ f \in L^2 $ and set
$$ u = e^{(2N+2)\phi(r)} \Pi^c (P+1)^{-N}  f $$
with $ N \geq N_0 $ large enough to be chosen. For convenience we have replaced $ 2N + 1/2 $ by $ 2N +2 $ which will be sufficient. We study first the contribution away from the boundary. It follows from Proposition \ref{theoremesurresolvante} that
$$ \big| \big|  \big( 1+ D_r^2 - \Delta_{\mathcal A} \big)^{n_0} \xi u \big| \big|_{L^2} \leq C || f ||_{L^2} , $$
provided that $ 2 n_0 + 2 N + 2 \leq 4 N  $.
Here we also use that commutations between $ D_r^j $ and powers of $ e^{\phi(r)} $ are harmless since all derivatives of $ \phi $ are bounded. Then, by standard Sobolev estimates in the cylinder $ \Ra \times {\mathcal A} $, we see that if $ n_0  $ is large enough\footnote{here $n_0 \geq 2$ since $ {\mathcal A} $ has dimension 1, but the same analysis would work for higher dimensional manifolds and larger $n_0$}, $ \xi u $ belongs $ L^{\infty} $ with norm controlled by $ ||f||_{L^2} $. 
We now consider $ (1 - \xi) u $.  We can drop the weight $ e^{(2N+2) \phi(r)} $ which is bounded on $ \mbox{supp}(1-\xi) $ and thus consider, according to (\ref{opdiaglapbis}),  
$$ (1-\xi) (r) \Pi^c (P_0+1)^{-2 N} f (r,\alpha) = \sum_{k \geq k_0} \mu_k^{2(2 N-1)}(1-\xi)(r)({\mathfrak p}_k + 1)^{-2 N} f_k (r) \mu_k^{-2(2N-1)}e_k (\alpha) . $$
Using that, for some $ c > 0  $ 
 $$ \mu_{k}^{2} e^{2 \phi(r)} + w(r) = \mu_k^{2} \big(e^{2 \phi(r)} + \mu_k^{-2} w(r) \big) \geq c \mu_k^2, \qquad k \gg 1 , \ r \geq r_0 ,  $$
 the operator $ {\mathfrak p}_k $ has spectrum contained in $ [ c \mu_k^2 , \infty ) $ for $k$ large enough. Thus
$$ \big| \big| (1-\xi) \Pi^c (P_0+1)^{-2N} f \big| \big|_{L^{\infty}} \lesssim \left( \sup_{k \geq k_0} ||  ({\mathfrak p}_k + 1)^{-1} ||_{L^2 \rightarrow L^{\infty}} \right) \left( \sum_{k \geq k_0} \mu_k^{-2(2 N-1)} || e_k ||_{L^{\infty}({\mathcal A})} \right) ||f||_{L^2} . $$
For $ N $ large enough, the second factor in the right hand side is bounded, using the boundedness of $ \mu_k^{1-2N} || e_k ||_{L^{\infty}} $ and the summability of $ \mu_k^{1-2N} $ for  $ N $ large enough which follow from standard rough estimates on eigenvalues and eigenfunctions on compact manifolds\footnote{as previously for $n_0$, this analysis holds whatever the dimension of $ {\mathcal A} $ is, provided $N$ is large enough.}. The finiteness of the first factor follows from the one dimensional Sobolev embedding $ H_0^1 (r_0,\infty) \subset L^{\infty} $ and from (\ref{injectionSobolevexplicite}) using that
\begin{eqnarray*}
 || ({\mathfrak p}_k + 1)^{-1} ||_{L^2 \rightarrow L^{\infty}} & \leq & C || ({\mathfrak p}_k + 1)^{-1} ||_{L^2 \rightarrow H_0^1 (r_0,\infty)}  \\
 & \leq & C || ({\mathfrak p}_k + 1)^{-1} ||_{L^2 \rightarrow {\mathfrak h}_{0,k}^1} 
\end{eqnarray*}
which is  is bounded uniformly in $k$ by construction of $ {\mathfrak h}_{0,k}^1 $ and $ {\mathfrak p}_k $. This completes the proof. 
\finpreuve

\bigskip

\noindent {\bf Remark.} The projection $ \Pi^c $ in Proposition \ref{aprioriL1} is needed to prove {\it global} $ L^2_{G_0} \rightarrow L^q_{G_0} $ estimates. If one is only interested in $ L^2_{G_0} \rightarrow L^q_{G_0, {\rm loc}} $ bounds, then we can drop $ \Pi^c $.

\bigskip

Let finally point out that (\ref{onS}) implies that, for any  $ r_1 > r_0 $,
$$ \psi \in \cap_j \mbox{Dom}(\Delta^j) \qquad \Longrightarrow \qquad  {\mathds 1}_{[r_1 , \infty )}(r) \Pi^c \psi \in L^{\infty}({\mathcal S}_0) . $$
This justifies the interest of considering  $ \cap_j \mbox{Dom}(\Delta^j) $ in Theorem \ref{theoremeonde} and Corollary \ref{corollaireSchrodinger} since in both cases this implies that $  {\mathds 1}_{[r_1 , \infty )}(r) \Pi^c \Psi $ belongs to $ C ( [0,1] , L^{q}_{G_0} ) $.

\section{Functional calculus} \label{section3}
\setcounter{equation}{0}
In this section, we provide asymptotic expansions of $ \varphi (-h^2 \Delta) $ and $ \varphi (-h^2 \Delta_0) $ in term of the semiclassical parameter $h \in (0,1]$, when $ \varphi \in C_0^{\infty}(\Ra) $. We shall use them in particular to justify the Littlewood-Paley decomposition.

We start by fixing some definitions and notation about properly supported operators and operator valued symbols. When $a$ is a scalar symbol, say in $ S^{-2} (\Ra^2) $,  we shall  replace the usual quantization (\ref{quantificationusuelle})
by a properly supported one, which has the advantage to preserve exponential decay or growth. For $ \kappa \in C_0^{\infty} (\Ra) $, we thus define the quantization $ O \! p_h^{\kappa} $ by
\begin{eqnarray}
 O \!p_h^{\kappa} (a) v (r) = (2 \pi)^{-1} \int \! \! \int e^{ i(r - s) \rho} a (r,h \rho) \kappa (s-r) v (s) d \rho ds  , \label{derivativesfalling}
\end{eqnarray} 
which defines a properly supported operator since its Schwartz kernel vanishes for $ s-r $ outside $ \mbox{supp}(\kappa) $.

If $a = a (r,\rho, \mu^2)$ depends on  $ \mu^2 \geq 0 $,  such that its seminorms in $ S^{-2}(\Ra^2_{r,\rho}) $ are uniform  with respect to $ \mu $,
one can define the operator valued symbol $ a (r,\rho,-h^2 \Delta_{\mathcal A}) $ by the spectral theorem for $ \Delta_{\mathcal A} $ and then the associated operator $  {\mathbf O}{\mathbf p}_h^{\kappa}(a)  $  by
\begin{eqnarray} 
 {\mathbf O}{\mathbf p}_h^{\kappa}(a) u (r,\alpha) &  = &  (2 \pi)^{-1} \int \! \! \int e^{ i(r - s) \rho} a (r,h \rho, -h^2 \Delta_{\mathcal A}) \kappa (s-r) u (s,\alpha) d \rho ds \label{operatorvaluedsymbol} \\ 
 & = & \sum_{k \geq 0} O \! p_h^{\kappa} (a(.,.,h^2 \mu^2_k)) u_k (r) e_{j}(\alpha) . \nonumber
\end{eqnarray}
This corresponds to (\ref{pourbornoppenheimer}) when $(A_k)_k = \big(O \! p_h^{\kappa} (a(.,.,h^2 \mu^ 2_k)) \big)_k $ is a sequence of pseudodifferential operators on the real line (with scalar symbols). To avoid to deal with the (possible) boundary and to be able to project away from the zero modes by $ \Pi^c $, we localize our operators inside $ {\mathcal S}_0 $.  To this end, we consider the cutoff $ \xi $ introduced in (\ref{introductiondexi}) and let
$ \kappa \in C_0^{\infty} (\Ra) $ be such that 
\begin{eqnarray}
\kappa \equiv 1 \ \ \mbox{near} \  0   , \qquad\mbox{supp}(\xi) + \mbox{supp} (\kappa) \subset (r_0,+\infty), 
\label{conditionsupport}
\end{eqnarray}
which implies that the Schwartz kernel of any operator of the form $ \xi O \! p_h^{\kappa}(a) $ is supported in $ (r_0,\infty)^2 $.

\begin{prop} \label{calculfonctionneltheoreme} 
 Let $ \varphi  \in C_0^{\infty} (\Ra) $. Then, we can find  symbols 
\begin{eqnarray*}
 a_{\varphi,0} \big(r,\rho,\mu^2\big) & = & \varphi \big(\rho^2 + \mu^2 e^{2 \phi(r)} \big) , \\
 a_{\varphi,j} \big(r,\rho,\mu^2\big) & = & \sum_{l=1}^{2j} p_{lj}(r,\rho,\mu^2 e^{2\phi(r)}) \varphi^{(l)} \big(\rho^2 + \mu^2 e^{2\phi(r)} \big) , \qquad j \geq 1 ,
\end{eqnarray*}
where each  $ p_{lj} (r,\rho,\eta) \equiv 0  $ if $ 2 l - j < 0 $ and otherwise is a universal ({\it i.e.} $ h,\mu $ and $ \varphi $ independent) linear combination  of 
\begin{eqnarray*}
( \mbox{product of derivatives of order $ \geq 1 $ of} \ \phi(r)) \times \rho^{m_1} \eta^{m_2} \ , & & \qquad 0 \leq m_1 + m_2 \leq 2 l -j , 
\end{eqnarray*}
such that, for all integers $ N,M $ such that $ N  - 4 M \geq 0 $
\begin{eqnarray}
  \Pi^c \xi \varphi (h^2 P) =  \Pi^c \xi \sum_{j=0}^{N-1} h^j {\mathbf O} {\mathbf p}_h^{\kappa} ( a_{\varphi,j}) + h^{N-4M}  \Pi^c \xi  R_{N,M}(h) , \label{pourlerestef}
\end{eqnarray}
 with a  remainder of the form
 \begin{eqnarray*}
 R_{N,M}(h)  =  (P+1)^{-M} B_{N,M}(h) (P+1)^{-M} 
\end{eqnarray*}
where $ B_{N,M}(h)$ is a bounded operator commuting with $ \Pi^c $ and such that  
\begin{eqnarray} 
  || B_{N,M}(h)   ||_{L^2 \rightarrow L^2} & \leq & C_{N,M} , \qquad  h \in (0,1] . \label{propertytoget}
\end{eqnarray}
\end{prop}

To prove Proposition 
 \ref{calculfonctionneltheoreme}, the point is to construct pseudo-differential parametrices for $ \varphi (h^2 {\mathfrak p}_k) $. This is very standard and somewhat elementary since $ h^2 {\mathfrak p}_k $ is a one dimensional Schr\"odinger operator. The only subtleties are that the potential $ h^2 \mu_k^2 e^{2 \phi (r)} $ depends unboundedly on $k$ and that $ e^{\phi(r)} $ may grow exponentially. Our proof below mainly focusses on these two issues.

Before turning to the proof, we state the corresponding result for $ \Delta $ on $ {\mathcal S} $ (keeping in mind that we see $ {\mathcal S}_0 $ as a special case). We recall that $ \xi $ localizes in the interior of $ {\mathcal S}_0  $ which is considered as a subset of $ {\mathcal S} $.

\begin{prop} \label{calculfonctionnelreecrit}
 For all $ \varphi \in C_0^{\infty} (\Ra) $ and all $ N \geq 0 $
$$ \Pi^c \xi \varphi (-h^2 \Delta) = \Pi^c \xi e^{\frac{\phi(r)}{2}} \left( \sum_{j=0}^{N-1} h^j {\mathbf O} {\mathbf p}_h^{\kappa} ( a_{\varphi,j})   \right) e^{- \frac{\phi(r)}{2}} + h^N \Pi^c \xi R_N (h) $$
where, for any $ M \geq 0 $,
$$  \big| \big| (1-\Delta)^M R_N (h) \big| \big|_{L^2_G \rightarrow L^2_G} \leq C h^{-2M} , \qquad h \in (0,1] . $$
\end{prop}

\bigskip

\noindent {\it Proof of Proposition \ref{calculfonctionneltheoreme}.} 
For each $ k \geq k_0 $, we consider the one dimensional Schr\"odinger operator 
$$ h^2 {\mathfrak p}_k = h^2 D_r^2 + h^2 \mu_k^2 e^{2 \phi (r)} + h^2 w (r) $$
and first build a parametrix for  $ (h^2 {\mathfrak p}_k - z)^{-1} $. To this end we recall the procedure for a one dimensional Schr\"odinger operator
$$ H (h) := h^2 D_r^2 + V (r) + h^2 b (r) = \sum_{j=0}^2 h^j O \! p_{h} (p_j) , $$
with $ V,b $  smooth, $ V > 0 $, $b$ bounded and  
$$ p_0 = \rho^2 + V (r), \qquad p_1 \equiv 0 , \qquad p_2 = b (r) . $$
For $ z \in \Ca \setminus  [ 0 , + \infty ) $ we can  construct iteratively
\begin{eqnarray*}
q_0 & = & \frac{1}{p_0-z} \\
q_{j} & = & - \frac{1}{p_0 - z} \sum_{n + m + l = j \atop m \leq j-1 } (p_n \# q_m )_{l} , \qquad j \geq 1
\end{eqnarray*}
where $ (a \# b)_l $ is the $ l$-th term in the expansion of the symbol of $ O \! p_h (a) O \! p_h (b) $. Explicitly, we have
\begin{eqnarray*}
q_j & = & - \frac{1}{p_0-z} \left( \frac{2}{i} \rho \partial_r q_{j-1}   - \partial_r^2 q_{j-2} + b q_{j-2} \right), \qquad j \geq 1 ,
\end{eqnarray*}
provided we set $ q_{-1} \equiv 0 $ to handle the contribution of $q_{j-2}$ when $j=1$. We then have
\begin{eqnarray}
 \big( H (h) - z \big) \sum_{j=0}^{N-1} h^j O \! p_{h} (q_j) = 1 + h^N O \! p_{h} (r_N) + h^{N+1} O \! p_{h} (\tilde{r}_N)  , \label{replacingquantization}
\end{eqnarray} 
where
$$ r_N = 2 \rho D_r q_{N-1} + b q_{N-2} , \qquad \tilde{r}_N = D_r^2 q_{N-1} + b q_{N-1} . $$
By induction, we see that for $ j \geq 1 $ we have
\begin{eqnarray}
 q_j = \sum_{ l = 1 }^{ 2 j} \frac{d_{l j}}{(p_0-z)^{1+l}}  \label{sumthosel}
\end{eqnarray} 
 where
 $$ d_{lj} = \ \mbox{universal linear combination of} \ \ \ \rho^{2 (l-l_1)-j} \prod_{s=1}^{l_1} \partial_r^{\nu_s} V \prod_{s=1}^{l_2} \partial_r^{\delta_s} b  $$
 with
 $$ \sum_{s=1}^{l_1} \nu_s + \sum_{s=1}^{l_2} \delta_s = j - 2 l_2 , \qquad 2 l - j \geq l_1 , \qquad l_1 + l_2 \leq l . $$
 Note in particular that we must have $ 2 l - j \geq 0 $ hence we may actually restrict the sum in (\ref{sumthosel}) to those $l$ such that $ \frac{j}{2} \leq l \leq 2 j $.
  Assuming that $ V > 0 $, we have $ V / (p_0+1) \leq 1 $ and $ |\rho| / (p_0+1)^{1/2} \leq 1 $. It is then not hard to check that
for all $ \gamma , \beta $, $ \big|\partial_r^{\gamma} \partial_{\rho}^{\beta} q_j \big| $ can be estimated by a constant (independent of $V$, $b$ and $z$) times
 $$  \scal{p_0}^{-1 - \frac{j}{2} - \frac{\beta}{2}} \left(1+ \max_{\delta \leq j + \gamma} \big| \big| b^{(\delta)}  \big| \big|_{L^{\infty}} \right)^{2j} 
\left(1 + \max_{\nu \leq j + \gamma} || V^{(\nu)} / V ||_{L^{\infty}} \right)^{2j+\gamma}  \left( \frac{\scal{z}}{|{\rm Im}(z)|} \right)^{2j + 1 + \gamma + \beta} . $$
Specializing this construction to $ V = h^2 \mu_k^2 e^{2 \phi(r)} $ and $b = w $, the above estimate and  (\ref{order}) show that 
\begin{eqnarray}
 \big|\partial_r^{\gamma} \partial_{\rho}^{\beta} q_j (r,\rho) \big| \leq C \frac{1}{(\rho^2 + \mu^2_k e^{2 \phi (r)} + 1 )^{1+(j+\beta)/2}} \left( \frac{\scal{z}}{|{\rm Im}(z)|} \right)^{2j+1+\gamma + \beta} \label{thankstoepsilon}
\end{eqnarray}
with a constant independent of $ h,k,z$, which is the main point.  Replacing the quantization $ O \! p_h $ by $ O\!p_h^{\kappa} $ in (\ref{replacingquantization}), we obtain
$$ \big( h^2 {\mathfrak p}_k - z \big) \sum_{j=0}^{N-1} h^j O \! p^{\kappa}_{h} (q_j) = 1 + h^N O \! p^{\kappa}_{h} (r_N) + h^{N+1} O \! p^{\kappa}_{h} (\tilde{r}_N) + \tilde{\tilde{R}} , $$
where $ \tilde{\tilde{R}} $ is an additional remainder term which is the contribution of the derivatives from $ h^2 D_r^2 $ falling on $ \kappa (s-r) $ (see (\ref{derivativesfalling})). By off diagonal decay, {\it i.e.} by integrating by part with $ h \partial_{\rho}/|r-s| $, this term is $ O (h^{\infty}) $. Note that we keep a uniform control of the symbol with respect to $k$ after such integrations by part thanks to (\ref{thankstoepsilon}). The interest of the properly supported quantization is that applying $ h^2 \mu_k^2 e^{2 \phi} $ to the right of  $ O \! p_h^{\kappa} (r_N) $ (or $ O \! p_h^{\kappa} (\tilde{r}_N) $ or $ \tilde{\tilde{R}} $) corresponds to the multiplication by 
$$ h^2 \mu_k^2 e^{2 \phi (s)} = h^2 \mu_k^2 e^{2 \phi(r)}  e^{2 ( \phi (r) - \phi(s) ) } $$
where $ h^2 \mu^2_k e^{2 \phi (r)} $ is controlled by $ (\rho^2 + h^2 \mu^2_k e^{2 \phi (r)} - z)^{-1} $, while $ e^{2 (\phi(r)-\phi(s))} $ is bounded on the support of $ \kappa (s-r) $. Compositions to the left do not cause any trouble. The interest of this remark is that operators of the form $ (h^2{\mathfrak p})_k^{M} O \! p^{\kappa}(r_N) (h^2{\mathfrak p})_k^{M}  $ (when $4M \leq  N  $) are bounded on $ L^2 $, with norm of polynomial growth in $ \scal{z}/ |{\rm Im}(z)| $. We need this property to get (\ref{propertytoget}). The rest of the proof is standard by using the Helffer-Sj\"ostrand formula (see {\it e.g.} \cite{DiSj1}) to pass from the resolvent of $ h^2 {\mathfrak p}_k $  to $ \varphi (h^2 {\mathfrak p}_k) $. The estimates on the remainder for $ \varphi (h^2 P) $ follow from the ones for the remainders of $ \varphi (h^2 {\mathfrak p}_k) $ by using (\ref{pourbornoppenheimer})-(\ref{bornebornoppenheimer}). \finpreuve

\bigskip

\noindent {\it Proof of Proposition \ref{calculfonctionnelreecrit}.} In the case when $ {\mathcal S} = {\mathcal S}_0 $, the result is a direct consequence of Proposition \ref{calculfonctionneltheoreme} and (\ref{pourunitaritereference}).  For a more general manifold $ {\mathcal S} $, it suffices to observe that the same parametrix as the one used on $ {\mathcal S}_0 $ will work since $ \xi $ localizes inside $ {\mathcal S}_0 $.  This is again fairly standard. We recall the main idea for the convenience of the reader.
We note first that the construction described in the proof of Proposition \ref{calculfonctionneltheoreme}  provides  a parametrix for  $ \Pi^c \xi (-h^2 \Delta_0 -z)^{-1} $.  Choosing $ \tilde{\xi} $ supported in the interior of $ {\mathcal S}_0 $ and equal to $1$ near $ \mbox{supp}(\xi) $, we compute
$$ \Pi^c \xi (-h^2 \Delta_0 - z)^{-1} \tilde{\xi} (-h^2 \Delta -z) = \Pi^c \xi \left(  1 +  (-h^2 \Delta_0 - z)^{-1} [\tilde{\xi},h^2 \Delta]  \right) . $$
Using the parametrix for $ \Pi^c \xi (-h^2 \Delta_0 - z)^{-1} $ and the fact  that $ \xi $ and $ [\tilde{\xi},h^2 \Delta] $  have disjoint supports, we get that for all $N  $ there exists $ C $ and $ M $ such that 
$$ \big| \big| (1-\Delta_0)^N \Pi^c \xi (-h^2\Delta_0 - z)^{-1} [\tilde{\xi},h^2 \Delta] \big| \big|_{L^2_G \rightarrow L^2_{G_0}} \leq C h^N \frac{\scal{z}^M}{|{\rm Im}(z)|^M}, \qquad z \notin [ 0 , \infty ), \ h \in (0,1] . $$ 
We thus obtain that, for all given $ N $,
$$  \Pi^c \xi (-h^2 \Delta_0 - z)^{-1} \tilde{\xi}  =  \Pi^c \xi (-h^2 \Delta -z)^{-1} +  O_{L^2_G \rightarrow H_{G_0}^{2N}} \left( h^N \frac{\scal{z}^{M}}{|{\rm Im}(z)|^{M}} \right) . $$
Using the Helffer-Sj\"ostrand formula and the parametrix for the left hand side obtained in the proof of Proposition \ref{calculfonctionneltheoreme},  we get the result.  \finpreuve

\bigskip

We end up this section with a  result on (microlocal) finite propagation speed. This will be useful to localize spatially our Strichartz estimates.
Let us fix $ r_1 > r_0 $ and $ \delta > 0 $ such that $ r_1 - 2 \delta > r_0 $.  Let us define
$$ {\mathds 1}_L (r) := {\mathds 1}_{[r_1,r_1+1]} (r-L), \qquad \widetilde{{\mathds 1}}_L (r) := {\mathds 1}_{[r_1-\delta,r_1 + 1 + \delta]} (r-L)  $$
for all $ L \geq 0 $. In particular, the multiplication by $ \widetilde{\mathds 1}_L $ maps $ L^2_G $ into $ L^2_{G_0} $. 

\begin{prop} \label{propspeed2eq} Let $ \nu \in \{ 1 , 1/2 \} $ and $ \varphi \in C_0^{\infty}(0,+\infty) $. Let $ r_1 > r_0 $ and $ \delta > 0 $ be as above. There exists $ t_0 > 0 $ such that, for all $ N \geq 0 $ and all $ q \geq 2 $, there exists $ C > 0 $ such that
$$ \left| \left| \Pi^c {\mathds 1}_L (r)  \varphi (-h^2 \Delta) e^{i \frac{t}{h} (-h^2\Delta)^{\nu}} - \Pi^c {\mathds 1}_L (r) \varphi (-h^2 \Delta_0) e^{i \frac{t}{h} (-h^2\Delta_0)^{\nu}} \widetilde{\mathds 1}_L (r) \right| \right|_{L^2_G \rightarrow  L^q_{G_0}} \leq C h^N e^{-N \phi(L)}  , $$
for all $ t \in [-t_0, t_0] $, $ h \in (0,1] $ and $ L \geq 0 $.
\end{prop}

The meaning of this proposition is that we have an upper bound for the propagation speed  in the radial direction which is uniform with respect to $ L $, both for the wave and Schr\"odinger equations (localized in frequency). It does not give  any information on the propagation speed in the cross section $ {\mathcal A} $ 
but it will be sufficient for our purpose.

\bigskip

We first reduce the problem to a question involving only $ \Delta_0 $. Let us consider the following property (P) 
\begin{enumerate}
\item[{\bf (P)}]{for all $ \chi_0 , \breve{\chi}_0 \in C_0^{\infty}(-r_1 - 2 \delta , r_1 + 1 + 2 \delta) $  such that
$$ \chi_0 \equiv 1 \ \ \mbox{near} \ [r_1,r_1 + 1], \qquad \breve{\chi}_0 \equiv 0 \ \ \mbox{near} \ \mbox{supp}(\chi_0) $$
there exists $ t_0 > 0 $ such that, for each $N, N_1 , N_2 \geq 0 $ there exists $ C $ such that, if we set $ \chi_L (r) = \chi_0 (r-L) $ and $ \breve{\chi}_L (r) = \breve{\chi}_0 (r-L) $,  
\begin{eqnarray}
\left| \left| \partial_r^{N_1} \breve{\chi}_L \varphi (-h^2 \Delta_0) e^{i\frac{t}{h}(-h^2 \Delta_0)^{\nu}} \Pi^c \chi_L (1-\Delta_0)^{N_2}  \right| \right|_{L^2_{G_0} \rightarrow L^2_{G_0}} \leq C h^N
\end{eqnarray} 
for all $ L \geq 0 $, all $ h \in (0,1] $ and all $ |t| \leq t_0 $. }
\end{enumerate}

\bigskip

\begin{lemm} \label{reductionEgorov} The property {\bf(P)} implies  Proposition \ref{propspeed2eq}.
\end{lemm}

\noindent {\it Proof.} We assume {\bf(P)} and prove Proposition \ref{propspeed2eq}. We consider first the case $ \nu = 1/2 $. We let $ \chi_0 , \tilde{\chi}_0 \in C_0^{\infty}(r_1 - \delta, r_1 + 1 + \delta) $ be such that $ \chi_0 \equiv 1 $ near $ [r_1 , r_1 + 1] $ and $ \tilde{\chi}_0  \equiv 1$ near $ \mbox{supp}(\chi_0) $. Let us set $ \tilde{\chi}_L (r) = \tilde{\chi}_0 (r-L) $ and define
\begin{eqnarray*}
 \widetilde{W}_0(t,h) & = & \tilde{\chi}_L\varphi (-h^2 \Delta_0) e^{i\frac{t}{h}(-h^2 \Delta_0)^{1/2}} \Pi^c \chi_L (1-\Delta_0)^{N_2}  , \\
  W (t,h) & = & \varphi (-h^2 \Delta) e^{i\frac{t}{h}(-h^2 \Delta)^{1/2}} \Pi^c \chi_L (1-\Delta_0)^{N_2} ,
\end{eqnarray*} 
for some $ N_2 > 0 $ to be chosen below.
Using that $ [\Delta,\tilde{\chi}_L] = 2 \tilde{\chi}_L^{\prime} \partial_r + \big(\tilde{\chi}^{\prime \prime}_L - \phi^{\prime} \tilde{\chi}_L^{\prime} \big) $ has bounded coefficients with supports disjoint from $ \mbox{supp}(\chi_L) $, it follows from {\bf (P)} that
\begin{eqnarray}
(\partial_t^2 - \Delta) \widetilde{W}_0 (t,h) & = & -
   [\Delta,\tilde{\chi}_L]  \varphi (-h^2 \Delta_0) e^{i\frac{t}{h}(-h^2 \Delta_0)^{1/2}} \Pi^c \chi_L (1-\Delta_0)^{N_2} \nonumber \\
   & = & O_{L^2_{G_0} \rightarrow L^2_G}(h^{\infty}) \label{Duhamelnew}
\end{eqnarray} 
for $ |t| \leq t_0 $ small enough independent of $ L $.
At $ t = 0 $, it follows from Proposition \ref{calculfonctionnelreecrit} that
\begin{eqnarray}
\widetilde{W}_0 (0,h)   & =& \tilde{\chi}_L \varphi (-h^2 \Delta_0)  \Pi^c \chi_L (1-\Delta_0)^N
\nonumber \\
& = & W (0,h) + O_{L^2_{G_0} \rightarrow L^2_G} (h^{\infty}) , \label{Duhamel0new}
\end{eqnarray}
since $ \chi_L \Pi^c \varphi (-h^2 \Delta_0) $ and $ \chi_L \Pi^c \varphi (-h^2 \Delta) $ have the same pseudo-differential parametrix. Here again, the remainder in (\ref{Duhamel0new}) is also uniform in $ L $.  Similarly, for the first derivative
\begin{eqnarray}
 \partial_t \widetilde{W}_0(0,h) 
   & = & \partial_t W (0,h) + O_{L^2_G \rightarrow L^2_G}(h^{\infty}) . \label{Duhamel1new}
\end{eqnarray}
By (\ref{Duhamelnew}), (\ref{Duhamel0new}), (\ref{Duhamel1new}) and the Duhamel formula, we obtain for all $ N $ the existence $ C > 0 $ such that 
$$ \left| \left|  \widetilde{W}_0 (t,h) - W (t,h) \right| \right|_{L^2_{G_0} \rightarrow L^2_G} \leq C_N h^N ,  $$
for all  $|t| \leq t_0$, $ h \in (0,1]$ and $ L \geq 0 $, that is, by taking the adjoint
\begin{eqnarray}
  \left| \left| (1-\Delta_0 )^{N_2}  \Pi^c \left( \chi_L \varphi (-h^2 \Delta_0)  e^{i\frac{t}{h}(-h^2 \Delta_0)^{1/2}} \tilde{\chi}_L - \chi_L  \varphi (-h^2 \Delta)  e^{i\frac{t}{h}(-h^2 \Delta)^{1/2}} \right) \right| \right|_{L^2_{G} \rightarrow L^2_{G_0}} \lesssim h^N . 
  \label{localisationsharptilde}
\end{eqnarray}  
Choose next  $ \tilde{\tilde{\chi}}_0 \in C_0^{\infty}(r_1 - 2 \delta, r_1 + 1 + 2 \delta) $ such that  $ \tilde{\tilde{\chi}}_0 \equiv 1 $ near $ [r_1 - \delta, r_1 + 1 + \delta] $ so that, if we set $ \tilde{\tilde{\chi}}_L (r) = \tilde{\tilde{\chi}}(r-L) $, we have
\begin{eqnarray}
 \tilde{\chi}_L - \widetilde{\mathds 1}_L  = \big( \tilde{\chi}_L - \tilde{\tilde{\chi}}_L \big) \widetilde{\mathds 1}_L 
 , \qquad \chi_L {\mathds 1}_L = {\mathds 1}_L . \label{firstandsecondconditions}
\end{eqnarray} 
The first condition in (\ref{firstandsecondconditions}) and {\bf (P)} imply, upon possibly decreasing $ t_0 $, that
$$ \left| \left| (1-\Delta_0 )^{N_2} \chi_L \Pi^c  \varphi (-h^2 \Delta_0)  e^{i\frac{t}{h}(-h^2 \Delta_0)^{1/2}} ( \tilde{\chi}_L - \widetilde{\mathds 1}_L ) \right| \right|_{L^2_{G_0} \rightarrow L^2_{G_0}} \leq C_N h^N , $$
for all $ |t| \leq t_0 $, $ L \geq 0 $ and $ h \in (0,1] $.
This allows to replace $ \tilde{\chi}_L $ by $ \widetilde{\mathds 1}_L $ in (\ref{localisationsharptilde}). By choosing $ N_2 $ large enough and using Proposition \ref{aprioriL1}, we can drop the operator $ (1-\Delta)^{N_2} $ and get the $ L^2_G \rightarrow L^q_{G_0} $ boundedness as well as the gain  $ e^{-N \phi (L)} $. Using the second condition in (\ref{firstandsecondconditions}), we can then replace $ \chi_L $ by $ {\mathds 1}_L $ in the  $ L^2_G \rightarrow L^q_{G_0} $ estimate and we get Proposition \ref{propspeed2eq}.
The proof for $ \nu = 1 $ is similar. \finpreuve

\bigskip

\noindent {\it Proof of Proposition \ref{propspeed2eq}.} By Lemma \ref{reductionEgorov} it suffices to prove {\bf (P)}. We thus fix $ \chi_0 $ and $ \breve{\chi}_0 $ as in the assumption of {\bf (P)}. By  (\ref{pourunitaritereference}), we may replace $ - \Delta_0 $ by $ P $ and $ L^2_{G_0} $ by $ L^2 $. Note that conjugating $ \partial_r $ by $ {\mathcal U} $ changes $ \partial_r$ into $ \partial_r + \phi^{\prime}(r)/2 $ which is harmless by (\ref{order}). To manipulate only bounded operators, we write
$$ \varphi (h^2 P) e^{i\frac{t}{h}(h^2 P)^{\nu}} = \varphi (h^2 P) e^{i\frac{t}{h} \theta_{\nu }(h^2 P)} $$
where $ \theta_{\nu} (\lambda) = \lambda^{\nu} $ near the support of $ \varphi $, $ \theta_{\nu} $ is smooth (since $ \mbox{supp} (\varphi) \Subset (0,\infty) $), real valued and constant for $ \lambda \gg 1 $. By non-negativity of $ P $, the definition of $ \theta_{\nu} $ on $ \Ra^- $ does not matter so we may choose $ \theta_{\nu} $ as the sum of a constant and of a $ C_0^{\infty} $ function. In particular, we can use Proposition \ref{calculfonctionneltheoreme} to describe $ \theta_{\nu}(h^2 P) $. Also, to handle the unboundedness of $ \partial_r^{N_1} $ and $ (1+P)^{N_2} $, we introduce $ \tilde{\varphi} \in C_0^{\infty}(\Ra) $ such that $ \tilde{\varphi} \varphi = \varphi $ and consider
$$ \partial_r^{N_1} \breve{\chi}_L \Pi^c \tilde{\varphi}(h^2 P) e^{i \frac{t}{h} \theta_{\nu}(h^2 P) }  \varphi (h^2 P) \Pi^c \chi_L (1+P)^{N_2}  . $$
By Proposition \ref{calculfonctionneltheoreme}, up to $ O (h^{\infty}) $ terms in operator norm (uniformly in $t$), the study of this operator is reduced to the one of operators of the form
$$ h^{-N_1 - 2 N_2} \Pi^c {\mathbf O}{\mathbf p}_h^{\kappa}(\tilde{a}_L) e^{i \frac{t}{h} \theta_{\nu}(h^2 P) }  {\mathbf O}{\mathbf p}_h^{\kappa}(b_L) $$
with symbols $ \tilde{a}_L (r,\rho,\mu^2) $ and $ b_L (r,\rho,\mu^2) $ such that
\begin{eqnarray}
 \mbox{supp}(\tilde{a}_L(.,.,\mu^2)) \subset \left\{ (r,\rho) \ | \ r \in \mbox{supp}(\breve{\chi}_L) \ \mbox{and} \ (r,\rho) \in \mbox{supp} \big( \tilde{\varphi} (\rho^2 + \mu^2 e^{2 \phi(r)}) \big) \right\} \label{support1}
\end{eqnarray} 
and
\begin{eqnarray}
 \mbox{supp}(b_L(.,.,\mu^2)) \subset \left\{ (r,\rho) \ | \ r \in \mbox{supp}( \chi_L) \ \mbox{and} \ (r,\rho) \in \mbox{supp} \big( \varphi (\rho^2 + \mu^2 e^{2 \phi(r)}) \big) \right\} , \label{support2} 
\end{eqnarray} 
and satisfying the bounds
\begin{eqnarray}
  \big| \partial_r^{\gamma} \partial_{\rho}^{\beta} \tilde{a}_L (r,\rho,\mu^2) \big| + \big| \partial_r^{\gamma} \partial_{\rho}^{\beta} b_L (r,\rho,\mu^2) \big| \leq C_{\gamma \beta} , \label{uniformbounds}
\end{eqnarray}  
uniformly with respect to $ L $ and $ \mu^2 $.   
By separation of variables it is then sufficient to show the exitence of $ t_0 > 0 $ with the property that for all $ N  $ there exists $ C > 0 $ such that
\begin{eqnarray}
h^{-2N_2 - N_1} \left| \left| O \! p_h^{\kappa} \big(\tilde{a}_L (.,.,h^2 \mu_k^2 ) \big) e^{i\frac{t}{h} \theta_{\nu}(h^2 {\mathfrak p}_k)} O \! p_h^{\kappa} \big(b_L (.,.,h^2 \mu_k^2) \big)  \right| \right|_{L^2 \rightarrow L^2} \leq C h^N , \label{estimationreduite}
\end{eqnarray}
for all $ |t| \leq t_0 $, $ h \in (0,1] $, $ L \geq 0 $ and $ k \geq k_0 $. The main point here is to show that the estimates and the time $t_0$ are uniform in $k$ and $ L $. This is a consequence of the Egorov Theorem (see {\it e.g.} \cite{Robe1}) as follows. Let $ \Phi^t_{h,k} $ be the flow of Hamiltonian vector field $ X_{h,k} $ associated to 
$ \theta_{\nu} (\rho^2 + h^2 \mu_k^2 e^{2 \phi(r)}) $, the principal symbol of $ \theta_{\nu}(h^2 {\mathfrak p}_k) $, {\it i.e.}
$$  X_{h,k} = 2 \theta_{\nu}^{\prime} \big(\rho^2 + h^2 \mu_k^2 e^{2 \phi(r)} \big) \left( \rho \frac{\partial}{\partial r} \\ - \phi^{\prime}(r) h^2 \mu_k e^{2 \phi (r)} \frac{\partial}{\partial \rho} \right) . $$
Since $ \mbox{supp}( \theta_{\nu}^{\prime} ) $ is compact, the components of $ X_{h,k} $ are bounded together with all their derivatives, uniformly with respect to $k$ and $h$. In particular, there exists a constant $ C $ independent of $h$ and $k$ such that
\begin{eqnarray}
 |\Phi_{h,k}^t (r,\rho) - (r,\rho)| \leq C |t| , \qquad r > r_0, \ \rho \in \Ra, \label{estimationflot}
\end{eqnarray} 
as long as $ \Phi_{h,k}^t (r,\rho) $ does not hit the boundary $ \{r_0\} \times \Ra_{\rho} $. This is true in particular for all $ t $ small enough (depending  on $ r_1 $ and $ \delta $) independent of $r $ and $ L $ since $ r \geq r_1 - 2 \delta > r_0 $. The Egorov Theorem implies that for all $ M $ we can write
\begin{eqnarray}
 e^{i\frac{t}{h} \theta_{\nu}(h^2 {\mathfrak p}_k)} O \! p_h^{\kappa} \big(b_L (.,.,h^2 \mu_k^2) \big) = O \! p_h^{\kappa} \big(b_{L,t,h,k,M}  \big) e^{i\frac{t}{h} \theta_{\nu}(h^2 {\mathfrak p}_k)} + O_{L^2 \rightarrow L^2} (h^M) , \label{algebredEgorov}
\end{eqnarray} 
with
$$ \mbox{supp}(b_{L,t,h,k,M}) \subset  \Phi^{-t}_{h,k} \left( \mbox{supp} ( b_L (.,.,h^2 \mu_k^2) ) \right)  $$
provided that $t$ is such that the right hand side is at positive distance from  $ \{r_0 \} \times \Ra_{\rho} $.
The bound on the remainder $ O_{L^2 \rightarrow L^2} (h^M) $ is uniform with respect to $t$ in any compact set (on which the flow does not reach the boundary) and with respect to $ L , k ,h $. Similarly, the $ L^{\infty} $ norms of the symbol $ b_{L,t,h,k,M} $ and its derivatives are bounded locally uniformly in $t$  and in $ h,k,L $.
By (\ref{support1}), (\ref{support2}) and (\ref{estimationflot}), there exists $ t_0 > 0 $ independent of $h,k,L$ such that the supports of $ b_{L,t,h,k,M} $ and $ \tilde{a}_L (.,.,h^2 \mu_k^2 ) $ are disjoint for $ |t| \leq t_0 $. Then, by standard pseudo-differential calculus, the composition of the corresponding operators is $ O (h^{\infty}) $ in $ L^2 $ operator norm, uniformly in $k,L $. Thus (\ref{estimationreduite}) follows from (\ref{algebredEgorov}) which completes the proof. \finpreuve

\section{Littlewood-Paley estimates} \label{section4}
\setcounter{equation}{0}
In this section, we provide a convenient Littlewood-Paley decomposition which will allow to localize the Strichartz estimates in frequency. We consider a spectral partition of unity,
\begin{eqnarray}
 1 = \varphi_0 (\lambda) + \sum_{ l \geq 0} \varphi (2^{-l} \lambda), \qquad \lambda \in \Ra , \label{spectralpartition}
\end{eqnarray}
with $ \varphi_0 \in C_0^{\infty} (\Ra) $ and $ \varphi \in C_0^{\infty} (\Ra \setminus 0) $. We also let $ \xi $ be the cutoff introduced in (\ref{introductiondexi}).

 This section is devoted to the proof of the following proposition.                

\begin{prop} \label{propositionLittPale} For all $ q \in [ 2 , \infty ) $, there exists $ C_q >0 $ such that, for all $ \psi \in  L^q_G $,
\begin{eqnarray}
|| \Pi^c \xi \psi ||_{L^q_{G_0}} \leq C_q \left( \sum_{h} || \Pi^c \xi \varphi (-h^2 \Delta) \psi  ||_{L^q_{G_0}}^2 \right)^{1/2} + C_{q} || \psi ||_{L^2_G} ,
\label{BonLP}
\end{eqnarray}
the sum being taken over all $ h $ such that $ h^2 = 2^{-l} $ with integers $ l \geq  0$.
\end{prop}

Note that the last term in the right hand side of (\ref{BonLP}) is well defined since $ L^q_G \subset L^2_G $ for $ {\mathcal S} $ has finite area.

Localized Littlewood-Paley decompositions on non compact manifolds have already been considered in \cite{BoucletLP} but in the context of manifolds with large ends.  We are here considering small ends deserving a different analysis; in particular, we  use the projection $ \Pi^c $ to avoid zero angular modes.

We explain first how to reduce Proposition \ref{propositionLittPale} to Proposition \ref{graceaKhintchine} below. For $ \psi \in C_0^{\infty} ({\mathcal S}) $ and $ \phi \in C_0^{\infty} ({\mathcal S}_0) $, we  can always write
$$ \big( \Pi^c \xi \psi ,  \phi \big)_{L^2_{G_0}} = \big( \Pi^c \xi \varphi_0 (-\Delta) \psi ,  \phi \big)_{L^2_{G_0}} + \sum_{h} \big( \Pi^c \xi  \varphi (-h^2 \Delta) \psi , \phi \big)_{L^2_{G_0}} , $$
by the Spectral Theorem and (\ref{spectralpartition}).
Using (\ref{onS}), we have on one hand
\begin{eqnarray}
 \big| \big| \Pi^c \xi \varphi_0 (-\Delta) \psi  \big| \big|_{L^q_{G_0}} \leq C || \psi ||_{L^2_G} . \label{fivefour}
\end{eqnarray} 
On the other hand, in the sum, let us write
$$ \Pi^c \xi \varphi (-h^2 \Delta) =  \tilde{\varphi} (-h^2 \Delta_0) \Pi^c \xi \varphi (-h^2 \Delta) + \big(1 - \tilde{\varphi} (-h^2 \Delta_0) \big) \Pi^c \xi \varphi (-h^2 \Delta)  $$
with $ \tilde{\varphi} \in C_0^{\infty} (0,+\infty) $ such that $ \tilde{\varphi} \equiv 1 $ near the support of $ \varphi $.  The second term of the right hand side is negligible  according to the following lemma.
\begin{lemm} \label{pseudofonc} There exists $ C$ such that
$$ \left| \left| \big(1 - \tilde{\varphi} (-h^2 \Delta_0) \big) \Pi^c \xi \varphi (-h^2 \Delta) \right| \right|_{L^2_G \rightarrow L^q_{G_0}} \leq C h , \qquad h \in (0,1 ] . $$
\end{lemm}

\noindent {\it Proof.} By Proposition \ref{aprioriL1}, it suffices to show that for all $ N $ there exists $ C $ such that
\begin{eqnarray}
 \left| \left| (1-\Delta_0)^N \big(1 - \tilde{\varphi} (-h^2 \Delta_0) \big) \Pi^c \xi \varphi (-h^2 \Delta) \right| \right|_{L^2_G \rightarrow L^2_{G_0}} \leq C h , \qquad h \in (0,1 ] . \label{ofthenorm} 
\end{eqnarray} 
Indeed, commuting $ (1-\Delta_0)^N $ with $ \big(1 - \tilde{\varphi} (-h^2 \Delta_0) \big) $, we compute $ (1-\Delta_0)^N \Pi^c \xi \varphi (-h^2 \Delta) $ using Proposition \ref{calculfonctionnelreecrit}. We get a sum of operators of the form
$ h^{j-2N} e^{\phi(r)/2} \Pi^c {\mathbf O} {\mathbf p}_h^{\kappa} (a_j) e^{-\phi(r)/2} $ with symbols $ a_j (r,\rho,\mu^2) $ supported in 
\begin{eqnarray}
 \{ (r,\rho) \ | \ r > r_0, \ \rho^2 + \mu^2 e^{2 \phi(r)} \in \mbox{supp}(\varphi) \} \label{premiersupport} ,
\end{eqnarray} 
and a remainder of size $ O_{ L^2_G \rightarrow L^2_{G_0} } (h) $ (actually $ O (h^M) $ for any $M$). For simplicity, we drop the dependence on $j$ in the sequel.  Inserting $ \tilde{\xi} $ such that $ \tilde{\xi}\equiv 1 $ near $ \mbox{supp}(\xi) $, we expand  similarly $ \big(1 - \tilde{\varphi} (-h^2 \Delta_0) \big) \tilde{\xi} \Pi^c $ as a term of order $ O_{L^2_{G_0} \rightarrow L^{2}_{G_0}} (h^{2N+1}) $ and a sum of terms of the form $  e^{\phi(r)/2} {\mathbf O} {\mathbf p}_h^{\kappa} (b)^* e^{-\phi(r)/2} \Pi^c $ with $ b (r,\rho,\mu^2) $ supported in
\begin{eqnarray}
 \{ (r,\rho) \ | \ r > r_0, \ \rho^2 + \mu^2 e^{2 \phi(r)} \in \mbox{supp}(1 - \tilde{\varphi}) \} . \label{deuxiemesupport}
\end{eqnarray} 
Up to terms of order $h$ (and actually $h^M$ for all $ M $) the estimate of the norm in (\ref{ofthenorm}) is reduced to the one of
$$ h^{-2N} \left| \left| e^{\frac{\phi(r)}{2}} \Pi^c {\mathbf O} {\mathbf p}_h^{\kappa} (b)^* {\mathbf O} {\mathbf p}_h^{\kappa} (a) e^{- \frac{\phi(r)}{2}} \right| \right|_{L^2_G \rightarrow L^2_{G_0}} = h^{-2N} \left| \left| \Pi^c {\mathbf O} {\mathbf p}_h^{\kappa} (b)^* {\mathbf O} {\mathbf p}_h^{\kappa} (a)  \right| \right|_{L^2 \rightarrow L^2} . $$
Since the sets (\ref{premiersupport}) and (\ref{deuxiemesupport}) are disjoint, it follows from standard pseudo-differential calculus and separation of variables that the norm above is of size $ h^{\infty} $, which completes the proof. \finpreuve

\bigskip

Let next $ \tilde{\xi} = \tilde{\xi} (r) $ be supported in $ [\tilde{r}_1,\infty)  $, with $ r_0 < \tilde{r}_1 < r_1 $, such that $ \tilde{\xi} \equiv 1 $ near $ \mbox{supp} (\xi) $. Lemma \ref{pseudofonc} implies that we can rewrite
\begin{eqnarray}
 \big(\Pi^c \xi \psi , \phi \big)_{L^2_{G_0}} & = &  \big(\Pi^c \xi \varphi_0 (-\Delta) \psi , \phi \big)_{L^2_{G_0}} +  \sum_{h} \big( \Pi^c \xi  \varphi (-h^2 \Delta) \psi , \Pi^c \tilde{\xi} \tilde{\varphi} (-h^2 \Delta_{0}) \phi \big)_{L^2_{G_0}} \nonumber \\
& & \ + \ O \big( || \psi ||_{L^2_G} || \phi ||_{L^{q^{\prime}}_{G_0}} \big) . \label{mieuxapreciser}
\end{eqnarray}
We then introduce the square functions
$$ S \psi := \left( \sum_h | \Pi^c \xi \varphi (-h^2 \Delta) \psi |^2  \right)^{1/2} , \qquad
 \widetilde{S}_0 \phi := \left( \sum_h | \Pi^c \tilde{\xi} \tilde{\varphi} (-h^2 \Delta_0) \phi |^2  \right)^{1/2} . $$
 Since $ q \geq 2 $, we recall that
\begin{eqnarray}
 ||S \psi ||_{L^q_{G_0}} \leq \left( \sum_{h} || \Pi^c \xi \varphi (-h^2 \Delta) \psi  ||_{L^q_{G_0}}^2 \right)^{1/2} . \label{elementary2}
\end{eqnarray} 
Assume for a while we have shown the following  
\begin{prop} \label{graceaKhintchine} For all $ q^{\prime} \in (1,2 ] $, there exists $ C_{q^{\prime}} $ such that
$$ || \widetilde{S}_0 \phi ||_{L^{q^{\prime}}_{G_0}} \leq C_{q^{\prime}} ||  \phi ||_{L^{q^{\prime}}_{G_0}} , $$
for all $ \phi \in C_0^{\infty}({\mathcal S}_0) $.
\end{prop}
Then we can prove Proposition \ref{propositionLittPale}.
 
 \bigskip
 
\noindent {\it Proof of Proposition \ref{propositionLittPale}.} By (\ref{fivefour}) and (\ref{mieuxapreciser}), there exists $ C > 0 $ such that
\begin{eqnarray}
 \left| \big( \Pi^c \xi \psi ,  \phi \big)_{L^2_{G_0}} - \sum_{h} \big( \Pi^c \xi \varphi (-h^2 \Delta) \psi , \Pi^c \tilde{\xi} \tilde{\varphi}(-h^2 \Delta_0) \phi \big)_{L^2_{G_0}} \right| \leq C || \psi ||_{L^2_G} || \phi ||_{L^{q^{\prime}}_{G_0}} . \label{reductionLP1}
\end{eqnarray}
On the other hand, using that
\begin{eqnarray}
\left| \sum_{h} \big( \Pi^c \xi  \varphi (-h^2 \Delta) \psi , \Pi^c \tilde{\xi} \tilde{\varphi} (-h^2 \Delta_0) \phi \big)_{L^2_{G_0}} \right| \leq ||S \psi ||_{L^q_{G_0}} || \widetilde{S}_0 \phi ||_{L^{q^{\prime}}_{G_0}}, \nonumber
\end{eqnarray}
together with (\ref{elementary2}), Proposition \ref{graceaKhintchine} and (\ref{reductionLP1}), we obtain
$$  \left| \big( \Pi^c \xi \psi ,  \phi \big)_{L^2_{G_0}} \right| \leq C || \phi ||_{L_{G_0}^{q^{\prime}}} \left( || \psi ||_{L^2_G}^2 +  \sum_{h} || \Pi^c \xi \varphi (-h^2 \Delta) \psi  ||_{L^q_{G_0}}^2 \right)^{1/2} . $$
The result follows by taking the supremum over those $ \phi $ such that $ || \phi ||_{L_{G_0}^{q^{\prime}}} = 1 $. \finpreuve

\bigskip

The rest of this section is  devoted to the proof of Proposition \ref{graceaKhintchine}. 
Let $ ( \epsilon_h )_{h^2 = 2^{-l} } $ be the usual Rademacher sequence, realized as a sequence of functions on $ [0,1] $ (see {\it e.g.} \cite{Stein}), and introduce the family of operators
\begin{eqnarray}
 B (t) := \sum_h \epsilon_h (t) \Pi^c \tilde{\xi} \tilde{\varphi} (-h^2 \Delta_0)  . \label{areduireent}
\end{eqnarray}
 Using the Khintchine inequality, Proposition \ref{graceaKhintchine} will follow from the existence of $ C_q > 0 $ such that
\begin{eqnarray}
 \big| \big| B (t) \big| \big|_{L^{q^{\prime}}_{G_0} \rightarrow L^{q^{\prime}}_{G_0}} \leq C_{q} , \qquad t \in [0,1].  \label{pourleL2fort}
\end{eqnarray} 
 For $q = 2 $, (\ref{pourleL2fort}) is a consequence of the spectral theorem and the fact that
\begin{eqnarray}
 \left| \sum_h \epsilon_h (t)  \tilde{\varphi} (h^2 \lambda ) \right| \leq C , \qquad t \in [0,1], \ \lambda \in \Ra , \label{aussipourmoyennenonnulle}
\end{eqnarray}
since $ | \epsilon_h(t)| \leq 1 $ and  at most a fixed finite number of terms in the sum don't vanish. Using the Marcinkiewicz interpolation Theorem (see {\it e.g.} \cite{Stein2}), Proposition \ref{graceaKhintchine} will then follow from (\ref{pourleL2fort}) with $q = 2$ and a weak $ L^1 $ bound on $ B (t) $ which we now prove.

To reduce this problem to an analysis of standard Calder\'on-Zygmund operators acting on sets equipped with the usual Lebesgue measure, it will be convenient to localize the problem in space. 
We thus consider, for $ L \geq 0 $,
$$ \mathds{1}_L (r) = \mathds{1}_{ [ \tilde{r}_1  , \tilde{r}_1  + 1 )} (r-L), \qquad \widetilde{\mathds{1}}_L (r) = \mathds{1}_{ [ \tilde{r}_1 - \delta , \tilde{r}_1  + 1 + \delta )} (r-L) , $$
where $ \delta $ and $ \tilde{r}_1 $ are fixed positive real numbers such that $ \tilde{r}_1 + \delta > r_0 $. In particular $ \mathds{1}_L $ and $ \widetilde{\mathds 1}_L $ are supported in $ (r_0,\infty ) $, $ \sum_L {\mathds 1}_L \equiv 1$ on $ \mbox{supp}(\tilde{\xi}) $ and $ \widetilde{\mathds{1}}_L \equiv 1$ near the support of $ \mathds{1}_L $.

The following definition will be useful.

\begin{defi} A sequence of operators $ (B_L)_{L \geq 0} $  on $ {\mathcal S}_0 $ is said to satisfy uniform weak $(1,1) $ bounds if for some $ C_B > 0 $,
\begin{eqnarray}
 \emph{vol}_{0} \big( \{ |B_L \psi | > \lambda \} \big) \leq C_B \lambda^{-1} || \psi ||_{L^1_{G_0}} , \qquad
\mbox{for all}  \ \psi \in L^1_{G_0}, \  \lambda > 0  \ \mbox{and} \ L \geq 0  . \label{conditionimpliquee}
\end{eqnarray}
\end{defi}

\begin{prop}  \label{sommefaible} Let $ (B_L)_{L \geq 0} $ be a sequence of operators on $ {\mathcal S}_0 $.
 \begin{enumerate} \item{If $ (B_L)_{L \geq 0} $ satisfies uniform weak $ (1,1) $ bounds, then 
$ B:= \sum_{L \geq 0} \mathds{1}_{L} (r) B_{L} \widetilde{\mathds{1}}_L (r) $ is of weak type $ (1,1) $, {\it i.e.}
\begin{eqnarray}
 \emph{vol}_{0} \big( \{ |B \psi | > \lambda \} \big)
  \leq C_{\phi \delta} C_B \lambda^{-1}  ||  \psi ||_{L^1_{G_0}} . \label{firststatement}
\end{eqnarray}} 
\item{ If each $ B_L $ has a range composed of functions supported in $ \{|r-L - \tilde{r}_1| \leq 1 \}  $  and if there exists  $ C > 0 $ such that
\begin{eqnarray}
 \emph{meas}_{dr d {\mathcal A}} \big( \{ | e^{-\phi(r)} B_L e^{\phi (r)} u | > \lambda \} \big) \leq C \lambda^{-1} || u ||_{L^1} , \label{implies}
 \qquad \mbox{for all} \ u \in L^1, \ \lambda > 0 , \ L \geq 0 ,
\end{eqnarray} 
then $ (B_L)_{L \geq 0} $ satisfies uniform weak $ (1,1) $ bounds, with a constant $ C_B = C C_{\phi \delta} $.}
\end{enumerate}
In both cases, $ C_{\phi \delta} $ are constants depending only on $ \phi $ and $ \delta $.
\end{prop}

We recall that  $ L^1 $ stands for $ L^1 ({\mathcal S}_0, dr d {\mathcal A}) $ while $ L^1_{G_0} = L^1 ({\mathcal S}_0, e^{-\phi(r)}dr d {\mathcal A}) $.

\bigskip

\noindent {\it Proof.} 1. The inequality (\ref{firststatement}) is a simple consequence of 
$$ \mbox{vol}_{0} \big( \{ |B \psi | > \lambda \} \big)  \leq  C_B \lambda^{-1} \sum_{L \geq 0} || \widetilde{\mathds{1}}_L \psi ||_{L^1_{G_0}} ,  $$
which follows from (\ref{conditionimpliquee}) and the fact that $ \{ |B \psi | > \lambda \} $ is contained in $ \bigcup_{L \geq 0} \{ |B_L \widetilde{\mathds{1}}_L \psi | > \lambda \} $.

\noindent 2. Let $ \psi = e^{\phi(r)} u $ so that $ || u ||_{L^1} = || \psi ||_{L^1_{G_0}} $. Using  (\ref{order}), there exists $ 0 < c = c (\phi,\delta) < 1 $ such that 
$$ c e^{-\phi(r)} \leq e^{-\phi(L)} \leq  c^{-1 }e^{-\phi(r)} , \qquad |r-L-\tilde{r}_1| \leq 1 . $$ 
This implies that
\begin{eqnarray*}
  \mbox{vol}_{0} \big( \{ |B_L \psi | > \lambda \} \big) & = &\int_{ \{ |r-L-\tilde{r}_1| \leq 1 \} \cap \{ | B_L \psi | > \lambda \}} e^{-\phi(r)}dr d {\mathcal A}  \\ & \leq & c^{-1}  e^{-\phi(L)} \mbox{meas}_{dr d {\mathcal A}} \big( \big\{ | e^{-\phi(r)}B_L e^{\phi(r)} u | > \lambda e^{-\phi(r)} \big\} \big)  
\end{eqnarray*}
and that
$$ \big\{ | e^{-\phi(r)}B_L e^{\phi(r)} u | > \lambda e^{-\phi(r)} \big\} \subset \big\{ | e^{-\phi(r)}B_L e^{\phi(r)} u | > c \lambda e^{-\phi(L)} \big\} $$
Using (\ref{implies}) with $ c e^{-\phi(L)} \lambda $ instead of $ \lambda $, we get (\ref{conditionimpliquee}). 
 \finpreuve

\bigskip

By Proposition \ref{calculfonctionneltheoreme}, we can write for any $ N $ and $ M $ 
\begin{eqnarray}
 \Pi^c \tilde{\xi} \tilde{\varphi} (- h^2 \Delta_0) = \Pi^c \tilde{\xi} e^{\frac{\phi(r)}{2}} \left( \sum_{j=0}^{N-1} h^j {\mathbf O} {\mathbf p}_h^{\kappa}(a_{\tilde{\varphi},j}) \right) e^{- \frac{\phi(r)}{2}}  + h^{N-2M}  B(h)  (1-\Delta_0)^{-M} \Pi^c  \label{inviewof}
\end{eqnarray} 
where $ || B(h) ||_{L^2_{G_0} \rightarrow L^2_{G_0}} \leq C $. Using Proposition \ref{aprioriL1} and the fact that $ L^2_{G_0} $ is contained in $ L^1_{G_0} $, we have for $ M $ large enough,
\begin{eqnarray}
 \big| \big|  B(h)  (1-\Delta_0)^{-M} \Pi^c  \big| \big|_{L^1_{G_0}  \rightarrow L^1_{G_0}} \leq C , \qquad h \in (0,1] , \nonumber
\end{eqnarray}
hence, if we choose $ N  $ such that $ N - 2 M > 0 $, we obtain the uniform $ L^1_{G_0} \rightarrow L^1_{G_0} $ bound 
\begin{eqnarray}
\left| \left| \sum_{h} \epsilon_h (t) h^{N-2M} B(h)  (1-\Delta_0)^{-M} \Pi^c \right| \right|_{L^1_{G_0} \rightarrow L^1_{G_0}} \leq C, \qquad t \in [0,1] . 
\end{eqnarray} 
This bound and (\ref{inviewof}) show that to study  weak type (1,1) estimates for (\ref{areduireent}) it suffices to study operators of the form

$$ B^a (t) = \sum_h \epsilon_h (t) \Pi^c  e^{\frac{\phi(r)}{2}} {\mathbf O}{\mathbf p}_h^{\kappa} (a) e^{- \frac{\phi(r)}{2}} , $$
with an operator valued symbol\footnote{see (\ref{operatorvaluedsymbol})} $a$ of the form
$$ a (r,\rho,h^2 \Delta_{\mathcal A}) = b (r) \rho^j a_0 \big(\rho^2 -  h^2 e^{2 \phi (r)} \Delta_{\mathcal A} \big) , $$
for some $ j \geq 0 $, $ b \in C^{\infty}_b (\Ra) $ with $ \mbox{supp}(b) \subset \mbox{supp}(\tilde{\xi}) $ and $ a_0 \in C_0^{\infty} (\Ra) $, with  $ \mbox{supp} (a_0) \subset \mbox{supp}(\tilde{\varphi}) $.  If the support of $ \kappa $ is small enough (depending on $ \delta  $ chosen in the definition of $ \widetilde{\mathds{1}}_L $), we have  $ \mathds{1}_{L}(r) {\mathbf O}{\mathbf p}_h^{\kappa} (a) =  \mathds{1}_{L} (r) {\mathbf O}{\mathbf p}_h^{\kappa} (a) \widetilde{\mathds{1}}_L (r) $,
so that
\begin{eqnarray}
 B^{a} (t) = \sum_{L \geq 0} \mathds{1}_{L} (r) B^{a}_L (t)  \widetilde{\mathds{1}}_L (r) \label{delocalisationspatiale}
\end{eqnarray} 
with
$$ B^{a}_L (t) = \sum_{h } \epsilon_h (t) \Pi^c    e^{\frac{\phi(r)}{2}} {\mathds 1}_L (r) {\mathbf O}{\mathbf p}_h^{\kappa} (a) e^{-\frac{\phi(r)}{2}} . $$ 
We will prove weak type $(1,1)$ estimates on $ B^a_L (t) $ by using the item 2 of Proposition \ref{sommefaible}. 
We thus consider $ e^{-\phi(r)} B_L^a(t) e^{\phi(r)} $ and the related kernels,
$$ K_{h,L} := \mbox{ Schwartz kernel of } \Pi^c e^{-\frac{\phi(r)}{2}} {\mathds 1}_L (r) {\mathbf O}{\mathbf p}_h^{\kappa} (a) e^{\frac{\phi(r)}{2}}  \mbox{ with respect to } dr d {\mathcal A} . $$
According to the standard theory of Calder\'on-Zygmund operators (see {\it e.g.} \cite{MuscSchl}), the weak $(1,1)$ estimates would follow from  the $  L^2 $ boundedness of $ e^{-\phi(r)} B^a_L (t) e^{\phi(r)} $ (uniformly in $L$ and $t$) and Calder\'on-Zygmund bounds on its Schwartz kernel. The $ L^2 $ boundedness is a consequence of the Calder\'on-Vaillancourt theorem as follows. Note first that the weights  $ e^{\pm \phi (r)} $ are harmless since our pseudo-differential operators are properly supported. We then observe that the full symbol of $ e^{-\phi(r)} B_L^a(t) e^{\phi(r)} $ obtained by summation over $h$ is bounded on $ \Ra^2 $, uniformly in $t$ and $ L $, thanks to the support properties of $a_0$  and the argument leading to (\ref{aussipourmoyennenonnulle}).  The same holds for the derivatives of the symbol, and this yields the $ L^2 $ boundedness.

We now focus on kernel bounds.

 We let $ d_{\mathcal A} $ be the geodesic distance on $ {\mathcal A} $ and, in the proposition below, denote by  $ |\cdot| $  either the usual modulus of a complex number or  the length of a covector with respect to the cylindrical Riemannian metric $ d r^2 + g_{\mathcal A} $.

\begin{prop} \label{noyauoperateur} \begin{enumerate} \item{There exists $ C > 0$ such that $ K_{h,L} \equiv 0 $ if $ h e^{\phi(L)} \geq C $.}
\item{ For all $ N $, there exists $ C $ such that
\begin{eqnarray}
 \big| K_{h,L}(r,\alpha,r^{\prime},\alpha^{\prime}) \big| & \leq & C e^{-\phi(L)} h^{-2} \left( 1 + \frac{|r-r^{\prime}|}{h} + \frac{d_{\mathcal A}(\alpha,\alpha^{\prime})}{e^{\phi(L)} h} \right)^{-N}  \nonumber \\
 \big| \partial_{r^{\prime}} K_{h,L}(r,\alpha,r^{\prime},\alpha^{\prime}) \big| & \leq & C e^{-\phi(L)} h^{-3} \left( 1 + \frac{|r-r^{\prime}|}{h} + \frac{d_{\mathcal A}(\alpha,\alpha^{\prime})}{e^{\phi(L)} h} \right)^{-N}  \nonumber \\  
 \big| d_{\alpha^{\prime}} K_{h,L}(r,\alpha,r^{\prime},\alpha^{\prime}) \big| & \leq & C e^{-2 \phi(L)} h^{-3} \left( 1 + \frac{|r-r^{\prime}|}{h} + \frac{d_{\mathcal A}(\alpha,\alpha^{\prime})}{e^{\phi(L)} h} \right)^{-N}  \label{length}
\end{eqnarray}
for all $ L \in \Na  $, all $ h \in (0,1] $ and $ (r,\alpha), (r^{\prime},\alpha^{\prime}) \in \Ra \times {\mathcal A} $.}
\end{enumerate} 
\end{prop}

To prove Proposition \ref{noyauoperateur}, we record a classical result in the next lemma.

\begin{lemm} \label{propzeta} Let $ \zeta \in C_0^{\infty} (\Ra) $ and $ \varepsilon_0 > 0 $. Then, for all $ N \geq 0 $, the Schwartz kernel of $ \zeta (-\varepsilon^2 \Delta_{\mathcal A}) $ satisfies
$$ \big| \big[ \zeta (-\varepsilon^2 \Delta_{\mathcal A}) \big] (\alpha,\alpha^{\prime}) \big| \leq C_{N \zeta} \varepsilon^{-1} \left( 1 + \frac{d_{\mathcal A} (\alpha,\alpha^{\prime})}{\varepsilon} \right)^{-N} , \qquad 0 < \varepsilon < \varepsilon_0 , $$
where $ d_{\mathcal A} $ is the Riemannian distance on $ {\mathcal A}$. We also have the estimate
$$ \big| d_{\alpha^{\prime}} \big[ \zeta (-\varepsilon^2 \Delta_{\mathcal A}) \big] (\alpha,\alpha^{\prime}) \big|_{g_{\mathcal A}} \leq C_{N \zeta} \varepsilon^{-2} \left( 1 + \frac{d_{\mathcal A} (\alpha,\alpha^{\prime})}{\varepsilon} \right)^{-N} , \qquad 0 < \varepsilon < \varepsilon_0 , $$
where $ d_{\alpha^{\prime}} $ is the differential acting on the second variable. The constant $ C_{N \zeta} $ remains bounded as long as $ \zeta $ belongs to a bounded set of $ C_0^{\infty}(\Ra) $.
\end{lemm}

This result actually holds for any compact Riemannian manifold of dimension $n$, up to the replacement of $ \varepsilon^{-1} $ (resp. $ \varepsilon^{-2} $) by $ \varepsilon^{-n} $ (resp. $ \varepsilon^{-n-1} $). It is a consequence of the standard semi-classical pseudo-differential expression of $ \zeta (- \varepsilon^2 \Delta_{\mathcal A}) $ (see {\it e.g.} \cite{BGT}).

\bigskip

\noindent {\it Proof of Proposition \ref{noyauoperateur}.} Define $ \varepsilon = e^{\phi(L)} h $. Then, $ K_{h,L} (r,\alpha,r^{\prime},\alpha^{\prime}) $ reads
$$ h^{-1} \mathds{1}_L (r) \left( \int e^{\frac{i}{h}(r-r^{\prime}) \rho} b(r)\rho^j \left[ (1-\pi_0) a_0 (\rho^2   - e^{2 (\phi(r) - \phi(L))} \varepsilon^2 \Delta_{\mathcal A} ) \right](\alpha,\alpha^{\prime}) \frac{ d \rho}{2 \pi} \right) \kappa (r-r^{\prime}) e^{\frac{\phi(r^{\prime}) - \phi(r)}{2}}   $$
where the bracket $ [ \cdots ] $ inside the integral corresponds to the Schwartz kernel on $ {\mathcal A} $ (according to the notation of Lemma \ref{propzeta}) and $ \pi_0 $ is the projection on $ \mbox{Ker}(\Delta_{\mathcal A}) $.
We observe that the above integal vanishes if  $ \varepsilon  $ is too large. Indeed, since  $ |r-L| \leq 1 $ thanks to $ \mathds{1}_L (r) $, it follows from (\ref{order}) that $  e^{\phi(r)-\phi(L)} $ is bounded from below hence, if $ \varepsilon $ is too large, we have
$$ (1-\pi_0) a_0 (\rho^2   - e^{2 (\phi(r) - \phi(L))} \varepsilon^2 \Delta_{\mathcal A} ) = 0 , $$
since this is equivalent to $ a_0 \big(\rho^2 +  e^{2 (\phi(r) - \phi(L))} \varepsilon^2 \mu_k^2 \big) \equiv 0 $ for all $  k \geq k_0 $. This proves the item 1 and shows, for the item 2, that may  assume that $ 0 < \varepsilon < \varepsilon_0 $ and so can use Lemma \ref{propzeta}. Using standard integrations by part in $ \rho $ to get a fast decay with respect to $ (1 + |r-r^{\prime}|/h) $, the result follows easily from Lemma \ref{propzeta} and the fact that
$$ \left( 1 + \frac{|r-r^{\prime}|}{h} \right)^{-N} \left( 1 + \frac{d_{\mathcal A}(\alpha,\alpha^{\prime})}{e^{\phi(L)} h} \right)^{-N} \leq  \left( 1 + \frac{|r-r^{\prime}|}{h} + \frac{d_{\mathcal A}(\alpha,\alpha^{\prime})}{e^{\phi(L)} h} \right)^{-N} . $$
We also use that $ \kappa (r-r^{\prime}) e^{\frac{\phi(r^{\prime}) - \phi(r)}{2}} $ and its derivatives are bounded, thanks to (\ref{order}) and the compact support of $  \kappa$. This completes the proof. \finpreuve

\bigskip

\noindent {\bf Proof of Proposition \ref{graceaKhintchine}.} By using a partition of unity on $ {\mathcal A} $, we  consider operators of the form $ \Theta_1 B_L^a (t) \Theta_2 $ with $ \Theta_1 , \Theta_2 \in C^{\infty}({\mathcal A}) $ either supported in the same coordinate patch or in disjoint coordinate patches. Let $ \theta_j $ be the coordinates defined on the support of $ \Theta_j $, $ j = 1,2 $. Define the operator $ D_L $ on $ \Ra^2 $ by
$$ ( D_L v ) (r,\theta) = e^{\phi(L)} v \big( r , e^{\phi(L)} \theta \big) . $$
The Schwartz kernel of $ D_L \theta_{1 * } \big( \Theta_1 \Pi^c  \mathds{1}_L (r) e^{-\frac{\phi(r)}{2}} {\mathbf O}{\mathbf p}_h^{\kappa} (a) e^{\frac{\phi(r)}{2}} \Theta_2 \big) \theta_2^* D^{-1}_L $ is of the form
$$ \tilde{K}_{L,h}(r,\theta,r^{\prime},\theta^{\prime}) =  e^{\phi(L)} K_{L,h} \left( r , \theta_1^{-1}(e^{\phi(L)}\theta) , r^{\prime} , \theta_2^{-1} (e^{\phi(L)}\theta^{\prime}) \right)  \beta_{1,2} \big(e^{\phi(L)}\theta, e^{\phi(L)}\theta^{\prime}\big) $$
for some compactly supported function $ \beta_{1,2} $. If $ \Theta_1 $ and $ \Theta_2 $ have disjoint supports, we may assume that $ \theta_1 , \theta_2 $ have disjoint ranges and that $ \beta_{1,2} $ is supported in $ I_1 \times I_2 $ with $ I_1, I_2 $ disjoint compact subsets of $ \Ra $. In any case, using Proposition \ref{noyauoperateur}, it is not hard to check that
$$ |\partial_{r^{\prime},\theta^{\prime}}^{\gamma} \tilde{K}_{L,h} (r,\theta,r^{\prime},\theta^{\prime})| \leq C_N h^{-2-|\gamma|} \left(1 + \frac{|r-r^{\prime}|+|\theta-\theta^{\prime}|}{h}\right)^{-N}, \qquad  |\gamma| \leq 1 . $$
We use basically that, if $ \theta_1 = \theta_2 $, then $ d_{\mathcal A} (\theta_1^{-1}(e^{\phi(L)}\theta) , \theta_2^{-1} (e^{\phi(L)}\theta^{\prime})) \approx e^{\phi(L)}|\theta - \theta^{\prime}| $ on the support of $ \tilde{K}_{L,h} $, while if $ \Theta_1 $ and $ \Theta_2 $ have disjoint supports, $ d_{\mathcal A} (\theta_1^{-1}(e^{\phi(L)}\theta) , \theta_2^{-1} (e^{\phi(L)}\theta^{\prime})) \gtrsim | e^{\phi(L)}\theta - e^{\phi(L)}\theta^{\prime}| $ since $ e^{\phi(L)} \theta \in I_1  $ and $ e^{\phi(L)} \theta^{\prime} \in I_2 $ on the support of $ \tilde{K}_{h,L} $. Therefore after summation in $h$ and using standard arguments, we see that the Schwartz kernel $ K_{t,L} $ of $ D_L \theta_{1 *} \big( \Theta_1 e^{-\phi(r)} B^a_L(t) e^{\phi(r)}  \Theta_2 \big) \theta_2^* D^{-1}_L  $ satisfies
$$  |\partial_{r^{\prime},\theta^{\prime}}^{\gamma} K_{t,L} (r,\theta,r^{\prime},\theta^{\prime})| \leq C \left(|r-r^{\prime}|+|\theta-\theta^{\prime}|\right)^{-2-\gamma}, \qquad  |\gamma| \leq 1 , $$
with a constant independent of $t$ and $ L $. Thus, by the usual Calder\'on-Zygmund theory, we have the uniform weak (1,1) estimates
\begin{eqnarray}
 \mbox{meas}_{dr d \theta} \left( \big\{ | D_L \theta_{1 *} \big( \Theta_1 e^{-\phi(r)} B^a_L(t) e^{\phi(r)}  \Theta_2 \big) \theta_2^* D^{-1}_L  v | > \lambda \big\} \right) \leq C \lambda^{-1} || v ||_{L^1 (\Ra^2)} \label{itfollowsfrom} 
\end{eqnarray} 
with $ C $ independent of $t$ and $ L $. Using on one hand that $ D_L^{-1} $ is an isometry on $ L^1 (\Ra^2) $ and on the other hand that
$$ \mbox{meas}_{dr d \theta} \big( \{ |D_L w| > \lambda \} \big) = e^{-\phi(L)} \mbox{meas}_{dr d \theta} \big( \{ |w| > e^{-\phi(L)} \lambda \} \big) $$
it follows from (\ref{itfollowsfrom}) that
$$ \mbox{meas}_{dr d {\mathcal A}} \left( \big\{ | \Theta_1 e^{-\phi(r)} B^a_L(t) e^{\phi(r)}  \Theta_2   u | > \lambda \big\} \right) \leq C \lambda^{-1} || u ||_{L^1 } $$
with a (possibly new) constant $ C $ independent of $ t $ and $ L $. By compactness of $ {\mathcal A} $, we only have to consider finitely many $ \Theta_1 , \Theta_2 $ so the same holds for $ e^{-\phi(r)} B_L^a (t) e^{\phi(r)} $ itself. Using Proposition \ref{sommefaible} and (\ref{delocalisationspatiale}), we obtain the expected weak $ (1,1) $ estimates for $ B^a (t) $ hence for (\ref{areduireent}). This leads to (\ref{pourleL2fort}) and thus completes the proof of Proposition \ref{graceaKhintchine}.  
\finpreuve

\section{Strichartz estimates} \label{proofs}
\setcounter{equation}{0}
In this section, we prove Theorems \ref{theoremeonde} and \ref{semiclassicalSchrodinger} as well as Corollary \ref{corollaireSchrodinger}. 

\subsection{Reduction of the problem}

In this paragraph, we explain how to reduce Theorems \ref{theoremeonde} and \ref{semiclassicalSchrodinger} and Corollary \ref{corollaireSchrodinger} to localized versions thereof. We will not only localize the estimates in frequency, as is classical, but also in space to handle the vanishing of the injectivity radius at infinity.
We will use the same kind of spatial localization as in previous sections namely we set
\begin{eqnarray}
 {\mathds 1}_L (r) = {\mathds 1}_{[r_1,r_1+1)} (r- L) , \qquad L \geq 0 , \label{unaile}
\end{eqnarray}
for a given $ r_1 > r_0 $.

\begin{prop}[Microlocal Schr\"odinger-Strichartz estimates] \label{reductionSchrodingerprop}  Let $ \varphi \in C_0^{\infty} (\Ra ) $.  There exists $ t_0 > 0 $ with the following property: for all  Schr\"odinger admissible pair $ (p,q) $ with $ p >2 $, there exists $ C > 0 $ such that, if we set
$$ \Psi_{h,L} (t) := {\mathds 1}_{L} (r) \varphi (-h^2 \Delta_0) e^{it  \Delta_0}  \psi_0 , $$
then
\begin{eqnarray}
|| \Pi^c \Psi_{h,L} ||_{L^p ( [0,h t_0] ; L^q_{G_0} ) } \leq C h^{-\sigma_{\rm S}}  || \psi_0 ||_{L^2_{G_0}}  , \label{microlocalSchrod}
\end{eqnarray}
 for all $ \psi_0 \in L^2_{G_0}
 $, $ h \in (0,1] $ and $ L \geq 0 $.
\end{prop}

We recall that $ \sigma_{\rm S} = \frac{1}{2p} $. There is a similar statement for the wave  equation for which we recall that $ \sigma_{\rm w} = \frac{3}{2p} $.

\begin{prop}[Microlocalized wave-Strichartz estimates] \label{reductionondesprop}  Let $ \varphi \in C_0^{\infty} (\Ra \setminus 0) $. There exists $ t_0 > 0 $ with the following property: for all sharp wave admissible pair $ (p,q) $, there exists $ C > 0 $ such that, if we set,
$$ \Psi_{h,L} (t) = {\mathds 1}_L (r) \varphi (-h^2 \Delta_0) e^{i t\sqrt{|\Delta_0|}} \psi_0  , $$
we have
\begin{eqnarray}
|| \Pi^c \Psi_{h,L} ||_{L^p ( [-t_0,t_0] ; L^q_{G_0} ) } \leq C h^{- \sigma_{\rm w}}  || \psi_0 ||_{L^2_{G_0}} 
, \label{microlocalwave}
\end{eqnarray}
 for all $ \psi_0 \in L^2_{G_0}
 $, $ h \in (0,1] $ and $ L \geq 0 $.
\end{prop}

We postpone the proofs of these two propositions to subsections \ref{soussectionSchrodinger} and \ref{soussectionondes} and first show how to use them to get Theorems \ref{theoremeonde}, \ref{semiclassicalSchrodinger} and Corollary \ref{corollaireSchrodinger}. We start with  the Schr\"odinger equation.

\bigskip

\noindent {\bf Proof of Theorem \ref{semiclassicalSchrodinger}.} For $ \psi \in L^2_G $ we let $ \Psi_h (t) =  e^{it \Delta} \varphi (-h^2 \Delta) \psi $ as in (\ref{notationPsih}). Write $$ {\mathds 1}_{[r_1 , \infty )}(r) = \sum_{L \geq 0} {\mathds 1}_L (r) . $$ Then using that $ p , q \geq 2 $ and the Minkowski inequality, we have for any $ T > 0 $,
\begin{eqnarray}
 \big| \big| \Pi^c {\mathds 1}_{[r_1 , \infty )}(r) \Psi_h \big| \big|_{L^p([0,T];L^q_{G_0})} \leq \left( \sum_{L \geq 0 } \big| \big| \Pi^c {\mathds 1}_{L}(r) \Psi_h \big| \big|_{L^p([0,T];L^q_{G_0})}^2 \right)^{1/2} . \label{deuxiemeorthogonalite}
\end{eqnarray} 
Let next $ \delta > 0 $ be such that $ r_1 - \delta > r_0 $ and set $ \widetilde{\mathds 1}_L (r) = {\mathds 1}_{[r_1 - \delta , r_1 + 1 + \delta)}(r-L) $. By Proposition \ref{propspeed2eq} (with $ \nu = 1 $) we can find $ \tau_0 > 0 $ small enough and $ C > 0 $ (both independent of $ \psi $, $ L $ and $h$) such that
\begin{eqnarray}
 \big| \big| \Pi^c {\mathds 1}_L (r) \Psi_h (t) \big| \big|_{L^q_{G_0}} \leq \big| \big| \Pi^c {\mathds 1}_L (r) e^{i t \Delta_0} \varphi (-h^2 \Delta_0) \widetilde{\mathds 1}_L (r) \psi \big| \big|_{L^q_{G_0}} + C  e^{- \phi (L)} || \psi ||_{L^2_G} , \label{avecresteSchrosinger1}
\end{eqnarray} 
for all $ \psi \in L^2_G $, all $ h \in (0,1] $, all $ L \geq 0 $ and all $ |t| \leq \tau_0 h $ (note that we are considering here $ e^{it \Delta} $ rather than $ e^{ith \Delta} $ in Proposition \ref{propspeed2eq}). Choose an integer $ N_0 > 0 $ large enough so that $ 1/N_0 \leq \min (\tau_0,t_0) $ with $t_0$ as in Proposition \ref{reductionSchrodingerprop}. Using Proposition \ref{reductionSchrodingerprop} and (\ref{avecresteSchrosinger1}), we have
$$  \big| \big| \Pi^c {\mathds 1}_{L}(r) \Psi_h \big| \big|_{L^p([0,h/N_0];L^q_{G_0})} \leq C h^{- \sigma_{\rm S}} || \widetilde{\mathds 1}_L (r) \psi ||_{L^2_{G_0}} + C e^{-\phi(L)} || \psi ||_{L^2_G} . $$
By using (\ref{deuxiemeorthogonalite}), the quasi-orthogonality of  the functions  $\widetilde{\mathds 1}_L (r) \psi $ and the summability of $ e^{-\phi(L)} $ given by Proposition \ref{propositionlemme}, we conclude that
$$ \big| \big| \Pi^c {\mathds 1}_{[r_1 , \infty )}(r) \Psi_h \big| \big|_{L^p([0,h/N_0];L^q_{G_0})} \leq C_0 h^{- \sigma_{\rm S}} || \psi ||_{L^2_G}  . $$
Using the group property and the unitarity of $ e^{it\Delta} $, we conclude that
\begin{eqnarray}
 \big| \big| \Pi^c {\mathds 1}_{[r_1 , \infty )}(r) \Psi_h \big| \big|_{L^p([0,h];L^q_{G_0})} \leq C_0 N_0^{1/p} h^{- \sigma_{\rm S}} || \psi ||_{L^2_G}  , \label{cumulationavecunitaire}
\end{eqnarray} 
using $ N_0 $ times the trick of \cite{BGT} which allows to derive Strichartz estimates on $ [0,h] $ from Strichartz estimates on $ [0,h/N_0] $. This completes the proof of Theorem \ref{semiclassicalSchrodinger}. \finpreuve

\bigskip

\noindent {\bf Proof of Corollary \ref{corollaireSchrodinger}.}  This step is basically as in \cite{BGT} up to minor technicalities due to the spatial localization. We choose the same $ \varphi $ as in Proposition \ref{propositionLittPale}. Using Theorem \ref{semiclassicalSchrodinger} and the same trick as in the end of the previous proof ({\it i.e.} cumulating $ O (h^{-1}) $ estimates on intervals of size $ h $ to get estimates on $[0,1] $), we have
\begin{eqnarray}
  \big| \big| \Pi^c {\mathds 1}_{[r_1 , \infty )}(r) \Psi_h \big| \big|_{L^p([0,1];L^q_{G_0})} \leq C h^{- \sigma_{\rm S} - \frac{1}{p}} || \psi ||_{L^2_G} , \label{etape1finalS}
\end{eqnarray}  
with $ \Psi_h $ given by (\ref{notationPsih}). 
Choosing $ \tilde{\varphi} \in C_0^{\infty}(\Ra \setminus 0) $ such that $ \tilde{\varphi} \varphi = \varphi $ we can replace $ || \psi ||_{L^2_G} $ by $ || \tilde{\varphi} (-h^2 \Delta_G) \psi ||_{L^2_G} $ and, by the spectral theorem, we have
\begin{eqnarray}
  h^{- \sigma_{\rm S} - \frac{1}{p}} || \tilde{\varphi} (-h^2 \Delta_G) \psi ||_{L^2_G} \leq C \left| \left| \tilde{\varphi} (-h^2 \Delta_G) (1-\Delta)^{\frac{\sigma_{\rm S}}{2} + \frac{1}{2p}} \psi \right| \right|_{L^2_G} . \label{etape2finalS}
\end{eqnarray}  
   On the other hand, choosing $ \tilde{r_1} $ such that $ r_0 < \tilde{r}_1 < r_1 $ (and a smooth cutoff $ \xi $ such that $ {\mathds 1}_{[r_1,\infty)} \leq \xi \leq {\mathds 1}_{[ \tilde{r}_1 , \infty )} $),
the Littlewood-Paley decomposition of Proposition \ref{propositionLittPale} yields
\begin{eqnarray}
 \big| \big| \Pi^c {\mathds 1}_{[r_1 , \infty )}(r) \Psi (t) \big| \big|_{L^q_{G_0}} \leq C \left( \sum_{h  } \big| \big| \Pi^c {\mathds 1}_{[ \tilde{r}_1,\infty )}(r) \Psi_h (t) \big| \big|_{L^q_{G_0}}^2 \right)^{1/2} + C || \psi ||_{L^2_G} , \nonumber
\end{eqnarray} 
with $ \Psi $ defined by (\ref{defPsiSchrodinger}).
Using the Minkowski inequality, this implies that
\begin{eqnarray}
 \big| \big| \Pi^c {\mathds 1}_{[r_1 , \infty )}(r) \Psi  \big| \big|_{L^p ( [0,1] ; L^q_{G_0} )} \leq C \left( \sum_{h  } \big| \big| \Pi^c {\mathds 1}_{[ \tilde{r}_1,\infty )}(r) \Psi_h  \big| \big|_{L^p ( [0,1] ; L^q_{G_0} )}^2 \right)^{1/2} + C || \psi ||_{L^2_G} . \nonumber
\end{eqnarray} 
Using (\ref{etape1finalS}) with $ \tilde{r}_1 $ instead of $r_1$ combined with (\ref{etape2finalS}) and the quasi-orthogonality of the operators $ \tilde{\varphi}(-h^2 \Delta) $, we get the result.  \finpreuve

\bigskip

\noindent {\bf Proof of Theorem \ref{theoremeonde}.} This proof is similar the previous one. We thus sketch the main ideas and only record in passing the modifications.  Let $ \psi \in L^2_G $ and set $ \Psi_h (t) = \Pi^c {\mathds 1}_{[r_1 , \infty )}(r) e^{it  \sqrt{\Delta}} \varphi (-h^2 \Delta) \psi $. By  (\ref{deuxiemeorthogonalite}) and Proposition \ref{propspeed2eq} with $ \nu = 1/2 $, there exists an integer $ N_0 > 0 $ and $ C > 0 $ such that
$$ ||  \Pi^c {\mathds 1}_{[r_1 , \infty )}(r) \Psi_h (t) ||_{L^q_G} \leq C \left( \sum_L \big| \big|  \Pi^c {\mathds 1}_{L}(r) e^{it \sqrt{|\Delta_0|}}\varphi (-h^2 \Delta_0) \tilde{\mathds 1}_L (r) \psi \big| \big|_{L^q_G}^2 + e^{-  \phi(L)} || \psi ||_{L^2_G}^2 \right)^{1/2}  $$
for all $ h \in (0,1 ] $, $ |t| \leq 1/N_0 $ and $ \psi \in L^2_G $. Here $ \tilde{\mathds 1}_L $ is as in the proof of Theorem \ref{semiclassicalSchrodinger}. By taking the $ L^2 ([-1/N_0;1/N_0],dt) $ norm in the above estimate combined with the Minkowski inequality, using the summability of $ e^{-\phi(L)} $ and Proposition \ref{reductionondesprop} (since we may assume that $ 1/N_0 \leq t_0 $), we get
\begin{eqnarray*}
 \big| \big|  \Pi^c {\mathds 1}_{[r_1 , \infty )}(r) \Psi_h  \big| \big|_{L^p ([-1/N_0,-1/N_0] ; L^q_G) } & \leq & C h^{- \sigma_{\rm w}} \left( \sum_L \big| \big|  \tilde{\mathds 1}_L (r) \psi \big| \big|_{L^2_G}^2  \right)^{1/2} + || \psi ||_{L^2_G} \\
 & \leq & C h^{- \sigma_{\rm w}} || \psi ||_{L^2_G} .
\end{eqnarray*} 
By  unitarity of $ e^{it \sqrt{|\Delta|}} $ and its group property, we can replace the interval $ [-1/N_0,1/N_0] $ by $ [1,1] $ up to the multiplication of $ C $ by the constant $ 2 N_0^{1/p} $. We may also replace $ || \psi ||_{L^2_G} $ by $ || \tilde{\varphi}(-h^2 \Delta) \psi ||_{L^2_G}  $ with $ \tilde{\varphi} \in C_0^{\infty}(\Ra \setminus 0) $ equal to $1$ on the support of $ \varphi $. Using the spectral theorem for the right hand side, this implies that
\begin{eqnarray}
  \big| \big|  \Pi^c {\mathds 1}_{[r_1 , \infty )}(r) \cos \big( t \sqrt{|\Delta|} \big) \varphi (-h^2 \Delta) \psi  \big| \big|_{L^p_t ([0,1] ; L^q_G) } \leq C \big| \big| \tilde{\varphi}(- h^2 \Delta ) (1-\Delta)^{\frac{\sigma_{\rm w}}{2}} \psi \big| \big|_{L^2_G} . \label{cosspectral} 
\end{eqnarray}  
We may also replace $ \cos $ by $ \sin $ and $ \varphi (-h^2 \Delta) $ by 
$$  \frac{1}{\sqrt{|\Delta|}} \varphi (-h^2 \Delta) = h \varphi_1 (-h^2 \Delta), \qquad \varphi_1 (\lambda) = \varphi (\lambda)/|\lambda|^{1/2} $$
to get
\begin{eqnarray}
  \left| \left|  \Pi^c {\mathds 1}_{[r_1 , \infty )}(r) \frac{\sin \big( t \sqrt{|\Delta|} \big)}{ \sqrt{|\Delta|}} \varphi (-h^2 \Delta) \psi  \right| \right|_{L^p_t ([0,1] ; L^q_G) } \leq C \big| \big| \tilde{\varphi}(- h^2 \Delta ) (1-\Delta)^{\frac{\sigma_{\rm w}-1}{2}} \psi \big| \big|_{L^2_G} . 
  \label{sinspectral}
\end{eqnarray}  
Using (\ref{cosspectral}), (\ref{sinspectral}) and Proposition \ref{propositionLittPale} (as in the proof of Corollary \ref{corollaireSchrodinger}), we get the estimates of Theorem \ref{theoremeonde}. 
\finpreuve

\bigskip

Before proving of Propositions \ref{reductionSchrodingerprop} and \ref{reductionondesprop} in the next two paragraphs, we proceed to some additional reductions and record useful results or remarks which will serve in both cases.

\medskip

\noindent {\bf First reduction.} By Proposition \ref{calculfonctionnelreecrit} (which we use for $ \Delta_0 $) combined with the Sobolev estimates of Proposition \ref{aprioriL1} to handle the remainders, it suffices to consider the terms of the pseudo-differential expansion of $ \varphi (-h^2 \Delta_0) $, namely 
$$  \Psi^{(\nu)}_{h,L}(t) := \Pi^c \mathds{1}_L (r) e^{\frac{\phi(r)}{2}} {\mathbf O}{\mathbf p}_h^{\kappa} (a_{\varphi,j}) e^{-\frac{\phi(r)}{2}} e^{-i \frac{t}{h} (-h^2 \Delta_0)^{\nu}} \psi_0 $$
rather than $ \Pi^c \mathds{1}_L (r) \varphi (-h^2 \Delta_0) e^{i t (-h^2 \Delta_0)^{\nu}/h} \psi_0 $. We omit the dependence on $j$ in the left hand side, since $j$ belongs to a finite set. More importantly notice that when $ \nu = 1 $ ({\it i.e.} for Schr\"odinger), we are considering a semiclassical time scaling (this will be eventually eliminated by (\ref{normeLph})). Let us remark that
 replacing $ \varphi (-h^2 \Delta_0) $ by its pseudo-differential expansion is standard, however   to handle the $ L^q $ norms of the remainders by Sobolev inequalities, we need the projection $ \Pi^c $ to be able to use Proposition \ref{aprioriL1}.

\medskip

\noindent {\bf Second reduction.} To estimate the $ L^q $ norms, we shall use first a Sobolev estimate in the angular variable, namely use the general fact that for $q \geq 2$,
\begin{eqnarray*}
\left| \left| \Pi^c  \Psi^{(\nu)}_{h,L}(t) \right| \right|_{L^q_{G_0}} & = & \left| \left| || \Pi^c \Psi^{(\nu)}_{h,L}(t,r,.) ||_{L^q ({\mathcal A})} \right| \right|_{L^q ((r_0,\infty), e^{-\phi(r)}dr)} \\
& \leq & C_{\mathcal A} \left| \left| || \Pi^c \sqrt{|\Delta_{\mathcal A}|}^{\frac{1}{2}-\frac{1}{q}} \Psi^{(\nu)}_{h,L} (t,r,.) ||_{L^2 ({\mathcal A})} \right| \right|_{L^q ((r_0,\infty),e^{-\phi(r)}dr)} \\
& \leq & C_{\mathcal A} \left( \sum_{k \geq k_0} \big| \big| \mu_k^{ \frac{1}{2} - \frac{1}{q}} \Psi^{(\nu)}_{h,L,k}(t) \big| \big|_{L^q ((r_0,\infty),e^{-\phi(r)}dr)}^2 \right)^{1/2} ,
\end{eqnarray*}
using the Minkowski inequality to get the third line since $ q \geq 2 $, and
where
\begin{eqnarray}
 \Psi^{(\nu)}_{h,L,k} (t,r) := \int_{\mathcal A} \overline{e_k(\alpha)} \Psi^{(\nu)}_{h,L} (t,r,\alpha) d {\mathcal A} ,  \nonumber
\end{eqnarray} 
according to the notation used in (\ref{Pythagore}).
For any $ p \geq 2 $ and any interval $ I $, the above estimate and the Minkowski inequality also yield
\begin{eqnarray}
 || \Psi^{(\nu)}_{h,L} ||_{L^p ( I; L^q_{G_0} )} \leq  C_{\mathcal A} \left( \sum_{k \geq k_0} \big| \big| \mu_k^{\frac{1}{2} - \frac{1}{q}} \Psi^{(\nu)}_{h,L,k} \big| \big|_{L^p (I ; L^q ((r_0,\infty),e^{-\phi(r)}dr) )}^2 \right)^{1/2} . \label{Strichartzl2}
\end{eqnarray}
This reduces the problem to get Strichartz inequalities for $ \Psi^{(\nu)}_{h,L,k} $.

\medskip

\noindent {\bf Third reduction.}  Using the definition of $ \Psi_{h,L,k}^{(\nu)} $ and the unitary equivalences given in (\ref{pourunitaritereference}) and Proposition \ref{separationdevariablesprop}, we have
\begin{eqnarray}
 \Psi^{(\nu)}_{h,L,k} (t,r) = e^{\frac{\phi(r)}{2}}\mathds{1}_L (r)  O \! p_h^{\kappa} (a_{\varphi,j}(.,.,h^2 \mu_k^2 ))  e^{-i \frac{t}{h} (h^2 {\mathfrak p}_k)^{\nu}} u_k , \label{expressiondeveloppeePsihLk}
\end{eqnarray} 
where,   according to (\ref{theunitarymapping}) and (\ref{Pythagore}),
\begin{eqnarray}
 u_k = ({\mathcal U} \psi_0 )_k . \label{unitaritefinale}
\end{eqnarray} 
 We refer to (\ref{operatorvaluedsymbol}) for $ O \! p_h^{\kappa} (a(.,.,h^2 \mu^2_k ) $. The form of $ a_{\varphi,j} $ given in Proposition \ref{calculfonctionneltheoreme} implies that
$$ \mbox{supp} ( {\mathds 1}_L (r) a_{\varphi,j}(.,.,h^2 \mu_k^2 )) \subset \big\{ (r,\rho) \ | \ |r-L| \leq C_0, \ \ \rho^2 + h^2 e^{2\phi(r)} \mu^2_k \in \mbox{supp}(\varphi) \big\} . $$
By (\ref{order}), $ \phi (r)- \phi (L) $ is bounded if $ |r-L| \leq C_0 $, so it follows that the set in the left hand side is empty if $ h^2 e^{2 \phi (L)} \mu_k^2 $ is too large. We may therefore assume that, for some $ C > 0 $, 
\begin{eqnarray}
 h \mu_k e^{\phi (L)} \leq C , \label{regionclassiquementinterdite}
\end{eqnarray}
since otherwise $ \Psi^{(\nu)}_{h,L,k} $ vanishes identically. This is a refinement the observation in the item 1 of Proposition \ref{noyauoperateur}.


\medskip

To prove the (one dimensional) Strichartz estimates, we will use the usual $ T T^* $ criterion \cite{KeelTao}. Since we are considering operators depending on several parameters (namely  $h,L,k$), we record a suitable version of this criterion. For notational simplicity, we let 
$$ {\mathcal H} := L^2 ((r_0,\infty) , dr ) , \qquad L^q (X) = L^q ( (r_0,\infty) , e^{-\phi(r)}dr ) . $$ If $  T (t) $ are time dependent operators from $ {\mathcal H} $ to $ L^2 (X) $, we set
$$ (T \psi ) (t,x) := \big( T (t) \psi \big) (x) .  $$

\begin{prop}[From dispersion to Strichartz] \label{translationStrichartz}  Let $ t_0 $, $ \sigma > 0 $ and $  \beta \geq 0 $ be real numbers. Then, for all real numbers $ p_1 > 2 $ and $q_1 \geq 2 $ such that
\begin{eqnarray}
\sigma \left( \frac{1}{2} - \frac{1}{q_1} \right) p_1 = 1 , \label{HLSouconvolution}
\end{eqnarray}
there exists a constant $ C $ such that for all family of operators $  \big( T_{h,L,k} (t) \big)_{h,L,k} $  satisfying
\begin{eqnarray}
\big| \big| T_{h,L,k} (t_1) T_{h,L,k} (t_2)^*  \big| \big|_{ L^1 (X) \rightarrow L^{\infty} (X)} & \leq & D_{h,L,k} h^{-\beta} \left( \frac{h}{|t_1-t_2| } \right)^{\sigma} \label{dispersionsemiclassique} \\
\big| \big| T_{h,L,k} (t) \big| \big|_{{\mathcal H} \rightarrow L^2 (X)} & \leq &  B_{h,L,k}
\end{eqnarray}
for all $ h \in (0,1] $, $ L \geq 0 $, $ k \geq k_0 $ and $ |t|, |t_1|, |t_2| \leq t_0 $, we have
\begin{eqnarray}
\big| \big| T_{h,L,k} u_k \big| \big|_{L^{p_1} ([0,t_0];L^{q_1} (X))} \leq C B_{h,L,k}^{\frac{2}{q_1}} D^{\frac{1}{2} - \frac{1}{q_1}}_{h,L,k} h^{-\beta \left( \frac{1}{2} - \frac{1}{q_1} \right)} h^{\frac{1}{p_1}}|| u_k ||_{\mathcal H} 
\end{eqnarray}
for all $ h \in (0,1] $, $ L \geq 0 $, $ k \geq k_0 $  and $ u_k \in {\mathcal H} $.
\end{prop}


We omit the proof of this proposition which follows from a standard interpolation argument (see {\it e.g.} \cite[Section 3]{KeelTao}) by tracking the dependence on the constants. We only note that the condition $ p_1 > 2 $ allows to use the Hardy-Littlewood-Sobolev inequality.

It follows from (\ref{expressiondeveloppeePsihLk}) that we have to consider
\begin{eqnarray}
 T_{h,L,k} (t) & := &   e^{\frac{\phi(r)}{2}}\mathds{1}_L (r)  A_k (h)  e^{-i \frac{t}{h} (h^2 {\mathfrak p}_k)^{\nu}} , \label{pourpreciserlenu} 
\end{eqnarray} 
with  $ A_k (h) =  O \! p_h^{\kappa}( a_{\varphi,j}(.,.,h^2 \mu_k^2 ) ) $. On one hand, we have
\begin{eqnarray}
\big| \big| T_{h,L,k} (t) \big| \big|_{{\mathcal H} \rightarrow L^2 (X)} & = & \big| \big| e^{- \frac{\phi(r)}{2}} T_{h,L,k} (t) \big| \big|_{L^2 \rightarrow L^2 } \nonumber \\
& = & \big| \big| \mathds{1}_L (r)  A_k(h) \big| \big|_{L^2 \rightarrow L^2 } \nonumber \\
& \leq & C , \label{estimeeL2H}
\end{eqnarray}
the last estimate being a consequence of the Calder\'on-Vaillancourt theorem since the symbol of $ A_k (h) $ and its derivatives are uniformly bounded with respect to $h,k$.
On the other hand
\begin{eqnarray}
\big| \big| T_{h,L,k} (t_1) T_{h,L,k} (t_2)^* \big| \big|_{L^1 (X) \rightarrow L^{\infty} (X)} & = & \big| \big|  T_{h,L,k} (t_1) T_{h,L,k} (t_2)^* e^{\phi(r)}  \big| \big|_{L^1 \rightarrow L^{\infty} } \nonumber \\
& = & \big| \big| e^{\frac{\phi(r)}{2}} \mathds{1}_L   A_k (h) e^{i\frac{(t_2-t_1)}{h} (h^2{\mathfrak p}_k)^{\nu}}  A_k (h)^* \mathds{1}_L e^{\frac{\phi(r)}{2}} \big| \big|_{L^1 \rightarrow L^{\infty} } \nonumber \\
& \lesssim &  e^{\phi(L)} \big| \big|  \mathds{1}_L   A_k (h) e^{i\frac{(t_2-t_1)}{h} (h^2{\mathfrak p}_k)^{\nu}}  A_k (h)^* \mathds{1}_L  \big| \big|_{L^1 \rightarrow L^{\infty} } \label{pourephiaile}
\end{eqnarray}
using in the last line that $ \phi (r)-\phi(L) $ is bounded on the support of $ \mathds{1}_L $ by (\ref{order}).
We recall that when no domain is specified $ L^q $ stands for $ L^q ( (r_0,\infty) ,d r ) $. In the next two paragraphs,  we derive explicit upper bounds for (\ref{pourephiaile}), first for $ \nu = 1 $ and then for $ \nu = 1/2 $.

\subsection{Proof of Proposition \ref{reductionSchrodingerprop}} \label{soussectionSchrodinger}
This corresponds to the case $ \nu = 1 $. As is well known, we shall derive the dispersion estimate by writing an approximation of $   \mathds{1}_L   A_k (h) e^{i t h{\mathfrak p}_k}  A_k (h)^* \mathds{1}_L  $  of the form $ \mathds{1}_L F (S,b) \mathds{1}_L  $ where $ F (S,b) $ will be our notation for  any Fourier integral operator with phase $S$ and amplitude $ b $,
$$ F (S,b) u (r) = (2 \pi h)^{-1} \int \int e^{ \frac{i}{h} \big( S(t,r,\rho) - r^{\prime} \rho \big)} b (t,r,\rho) d \rho u (r^{\prime})d r^{\prime} . $$
 The construction of such an approximation is standard (see {\it e.g.} \cite[ch. 10]{DiSj1} or \cite{Robe1}) up to the fact that we need to control it with respect to the parameters $ h , k , L $. We recall the main points of the construction. We let the principal symbol of $ h^2 {\mathfrak p}_k $ be 
\begin{eqnarray} 
 H_{h,k} (r,\rho) :=  \rho^2 + h^2 \mu_k^2 e^{2 \phi (r)} . \label{principalnuegal1}
\end{eqnarray}
Occasionally, the phase space variables $ (r,\rho) $ will be replaced by $ (r,\eta) $ or $ (x,\xi) $.
By the property (\ref{order}), for all multi-index $ \gamma $ there exists a constant $ C_{\gamma} > 0 $ such that 
\begin{eqnarray}
 \big|\partial_{r,\rho}^{\gamma}H_{h,k}(r,\rho) \big| \leq C_{\gamma} (1+H_{h,k}(r,\rho)), \qquad (r,\rho) \in \Ra^2, \ h > 0, \ k \geq k_0 . \label{orderaparametre}
\end{eqnarray}
The uniformity of this estimate with respect to $h$ and $ k $ imply that for all open relatively compact intervals $ J_1 \Subset J_2 \Subset \Ra  $, there exists $ \epsilon > 0 $ such that for all $ h , k $,
\begin{eqnarray}
 H_{h,k}^{-1} (J_1) + (-\epsilon,\epsilon)^2 \subset H_{h,k}^{-1} (J_2) . \label{voisinageenergie}
\end{eqnarray} 
There exists also a time $ T > 0 $ independent of $h$ and $k$ such that the Hamiltonian flow $ \Phi^t_{h,k} $ of $ H_{h,k} $ is defined on $ H_{h,k}^{-1} (J_2) $ for $ |t| \leq T $ and satisfies
\begin{eqnarray}
 \big| \big| d_{r,\rho} \Phi_{h,k}^t - I \big| \big| \leq C |t|, \qquad h> 0, \ k \geq k_0 , \ (r,\rho) \in  H_{h,k}^{-1} (J_2) , \ |t| \leq T . \label{borne1phase}
\end{eqnarray} 
Here $ || \cdot || $ is the matrix norm associated to the euclidean norm on $ \Ra^2 $. Denoting $ \Phi_{h,k}^t = \big( x^t_{h,k} , \xi_{h,k}^t \big) $,  (\ref{voisinageenergie}) and (\ref{borne1phase}) allow to choose $ t_0 > 0 $ small enough independent of $h$ and $k$ such that there exists a smooth function $ \eta_{h,k}^t $ defined on $ H_{h,k}^{-1} (J_1) $ for $ |t|\leq t_0 $ such that  
\begin{eqnarray}
 (r,\eta_{h,k}^t (r,\rho)) \in  H_{h,k}^{-1} (J_2)  , \qquad  \xi_{h,k}^t \big( r , \eta_{h,k}^t (r,\rho) \big) = \rho   , \label{borne2phase}
\end{eqnarray} 
for all $ (r,\rho) \in H_{h,k}^{-1}(J_1) $ and $ |t|\leq t_0 $.
Then, by the standard Hamilton-Jacobi theory there exists  a smooth function $ S_{h,k} : (-t_0,t_0) \times H_{h,k}^{-1}(J_1) \rightarrow \Ra $ which solves the Hamilton-Jacobi equation
$$ \partial_t S_{h,k}  =  H_{h,k} (r, \partial_r S_{h,k}) , \qquad S_{h,k}(0,r,\rho) = r \rho , $$
and which, by construction, also satisfies
\begin{eqnarray}
 \partial_{\rho} S_{h,k} (t,r,\rho) = x^t_{h,k} \big( r , \eta^t_{h,k}(r,\rho) \big) . \label{borne3phase}
\end{eqnarray} 
Then,  by solving the relevant transport equations for the amplitude according to the standard procedure \cite{Robe1} (with a uniform control on $ h,k,L$ which comes from (\ref{regionclassiquementinterdite})), we obtain the following result.
\begin{prop} \label{optiquegeometriqueSchrodinger} Let $ \nu = 1 $. Let $ J_0 \Subset J_1 \Subset J_2  $ be relatively compact intervals with $ \emph{supp}(\varphi) \subset J_0 $. Let $ \delta > 0 $ be such that $ r_1 - \delta > r_0  $.\footnote{recall that $r_1$ is used in (\ref{unaile})} Then there exists  $ t_0 > 0 $ such that for all $ N \geq 0 $,  we can find $ b_{h,k,L} \in C^{\infty}  \big( (-t_0,t_0) \times \Ra^2 \big) $ such that
\begin{eqnarray}
 \left| \left|  \mathds{1}_L  \left( A_k (h) e^{i t ( h^2 {\mathfrak p}_k )^{\nu} / h }  A_k (h)^*  - F_h(S_{h,k}, b_{h,k,L})(t) \right) \mathds{1}_L \right| \right|_{L^1 \rightarrow L^{\infty}} \leq C h^N , \label{ecartcontrole}
\end{eqnarray} 
for all $ h \in (0,1 ] $, $ k \geq k_0 $, $ L \geq 0 $ and $ |t| \leq t_0 $. The amplitude satisfies
$$  b_{h,k,L}(t,.,.) \in  C_0^{\infty} \big( H^{-1}_{h,k}(J_1) \cap \{  r_1 - \delta \leq r - L \leq r_1 + 2 \} \big) $$
with bounds on $ || \partial_{\rho} b_{h,k,L}(t,r,.) ||_{L^1} $ uniform with respect to $ h,k,L, t$ and $r$. In addition,
 (\ref{borne1phase}), (\ref{borne2phase}) and (\ref{borne3phase}) hold.
\end{prop}

To get a $ L^1 \rightarrow L^{\infty} $ estimate for $  \mathds{1}_L F_h(S_{h,k}, b_{h,k,L})(t)  \mathds{1}_L $, we will use
 the following classical result (see {\it e.g.} \cite{Stein2}, Section 1.2 of Chapter VIII).
\begin{prop}[Van der Corput estimates] \label{VanderCorput} There exists a universal constant $ C > 0 $ such that for all real numbers $ a < b $, all  $ S \in C^2 ([a,b],\Ra) $ such that $ S^{\prime \prime} > 0 $ and all $ b \in C^1 ([a,b]) $,
$$ \left| \int_a^b e^{i S (\rho)} b(\rho) d \rho \right| \leq C  \frac{||b||_{L^{\infty}} + ||b^{\prime} ||_{L^1} }{\sqrt{\min_{[a,b]} S^{\prime \prime}}} . $$
\end{prop}

We thus have to estimate $ \partial_{\rho}^2 S_{h,k} (t,r,\rho) $ for $t$ small enough, independent of $ h,k,L $, of $r$ in the support of $ \mathds{1}_L $ and $ \rho $ in the set $ \{ \rho \in \Ra \ | \ (r,\rho) \in H^{-1}_{h,k} (J_1) \} $. Note that this set is the union of at most two open intervals (it may be empty or consists of only one interval).  We present here a general method which we can use again for the wave equation in the next subsection. For simplicity, we drop the dependence on $ h,k $ from the notation.  By differentiating (\ref{borne3phase}) and using (\ref{borne1phase})-(\ref{borne2phase}), we may assume that $t$ is small enough so that
\begin{eqnarray}
 \partial_{\rho}^2 S (t,r,\rho) = \partial_{\eta} x^t (r,\eta^t(r,\rho)) \partial_{\rho}  \eta^t (r,\rho) \geq  \frac{1}{2} \big( \partial_{\eta} x^t \big) (r,\eta^t(r,\rho)) . \label{refonde1}
\end{eqnarray} 
On the other hand, the Hamilton equations yield
$$ \partial_{\eta} x^t (r,\eta) = \int_0^t \frac{\partial^2 H}{\partial \xi^2} (x^s,\xi^s) \partial_{\eta} \xi^s ds +  \int_0^t \frac{\partial^2 H}{\partial x \partial \xi} (x^s,\xi^s) \partial_{\eta} x^s ds  $$
where $ (x^s,\xi^s) = (x^s,\xi^s ) (r,\eta) $. Letting $ {\mathcal T} = {\mathcal T}(t,r,\rho) $ be the trajectory, {\it i.e.} the range of $ (x^s , \xi^s) $ when $s$ varies between $ 0$ and $t$, and using (\ref{borne1phase}) for $ \partial_{\eta} \xi^s $, we see that if $ t$ is small enough
\begin{eqnarray}
\partial_{\eta} x^t (r,\eta) \geq \frac{|t|}{2} \inf_{{\mathcal T}} \frac{\partial^2 H}{\partial \xi^2} - C t^2 \sup_{{\mathcal T}} \left| \frac{\partial^2 H}{\partial x \partial \xi} \right| . \label{borneinferieurephaseexplicite}
\end{eqnarray}
Here the constant $ C $ follows from the estimate $ |\partial_{\eta} x^s | \leq C |s| $ which follows  from (\ref{borne1phase}), (\ref{refonde1}) and the fact that we have upper bounds on $ |\partial_{\xi}^2 H| $ and $ |\partial_{x \xi}^2 H| $ on the energy shell  $ H^{-1}(J_2) $ according to (\ref{orderaparametre}).  In the present case, we have $ \partial_{\xi}^2 H = 2 $ and $ \partial_{x,\xi}^2 H \equiv 0 $. Putting back the parameters, we thus conclude that there exists $ t_0 > 0 $ such that
\begin{eqnarray}
 \partial_{\rho}^2 S_{h,k}(t,r,\rho) \geq |t| , \qquad |t| \leq t_0,  \label{lowerboundSchrodinger}
\end{eqnarray} 
for all $ h \in (0,1 ] $, $ k \geq k_0 $, $ L \geq 0 $, $ r \in \mbox{supp} ( \mathds{1}_L ) $ and $ \rho $ such that $ (r,\rho) \in H^{-1}_{h,k} (J_1) $.  By Proposition \ref{VanderCorput}, we conclude that
$$ \big| \big|  \mathds{1}_L F_h(S_{h,k}, b_{h,k,L})(t_1 - t_2)  \mathds{1}_L  \big| \big|_{L^1 \rightarrow L^{\infty}} \leq C h^{-1} \left( \frac{h}{|t_1 - t_2|} \right)^{1/2}, \qquad |t_1|,| t_2| < t_0 .$$
By (\ref{ecartcontrole}), the same holds for $  \mathds{1}_L  A_k (h) e^{i (t_1-t_2)  h {\mathfrak p}_k  }  A_k (h)^* {\mathds 1}_L $,
so by using (\ref{estimeeL2H}) and (\ref{pourephiaile}), we obtain the following

\begin{prop} \label{clarificationSchrodinger}  Let $ (p,q) $ be sharp Schr\"odinger admissible in dimension 2 with $ p,q $ real and consider
$$ q_1 = q , \qquad p_1 = 2 p , \qquad \beta = 1, \qquad \sigma = \frac{1}{2} . $$
Then there exists $ C > 0 $ such (\ref{dispersionsemiclassique}) holds true with $ B_{h,L,k} = C $ and $ D_{h,L,k} = C e^{\phi(L)} $.
\end{prop}

We can now complete the proof of Proposition \ref{reductionSchrodingerprop}.

\bigskip

\noindent {\bf Proof of Proposition \ref{reductionSchrodingerprop}.} We observe that
\begin{eqnarray}
\mu_k^{\frac{1}{2}-\frac{1}{q}}\big| \big| \Psi_{h,L,k}^{(1)} \big| \big|_{L^p ([0,t_0],L^q (e^{-\phi(r)}dr))} & \leq & C_{t_0} \mu_k^{\frac{1}{2}-\frac{1}{q}} \big| \big| \Psi_{h,L,k}^{(1)} \big| \big|_{L^{2p} ([0,t_0],L^q (e^{-\phi(r)}dr))} \\
& \leq & C  \big( \mu_k e^{\phi(L)} \big)^{\frac{1}{2} - \frac{1}{q}} h^{- \frac{1}{2} \left( \frac{1}{2} - \frac{1}{q} \right)} || u_k ||_{L^2} \\
& \leq &  C h^{- \frac{3}{2 p}} ||u_k||_{L^2}
\end{eqnarray}
the second line following from Propositions \ref{translationStrichartz} and \ref{clarificationSchrodinger}, and the 
 third one from (\ref{regionclassiquementinterdite}) and the admissibility condition. We conclude by using (\ref{Strichartzl2}) together with 
\begin{eqnarray}
 \big| \big| f \big| \big|_{L^p ([0,h t_0])} = h^{\frac{1}{p}} \big| \big| f_h \big| \big|_{L^p ([0,t_0])} , \qquad f_h (t) = f (ht)
 \label{normeLph}
\end{eqnarray}
to eliminate the semiclassical time scaling, and with the fact that $ \sum_{k \geq 0} ||u_k ||_{L^2}^2 = || \psi_0 ||_{L^2_{G_0}}^2 $ by (\ref{unitaritefinale}). \finpreuve

\subsection{Proof of Proposition \ref{reductionondesprop}} \label{soussectionondes}

The global strategy is the same as in subsection \ref{soussectionSchrodinger} so we only explain the main changes. We are considering $ \nu = 1/2$, namely $ e^{it \sqrt{{\mathfrak p}_k}} $ and use the classical hamiltonian 
\begin{eqnarray*}
 H_{h,k}^{1/2} (r,\rho) = \sqrt{ \rho^2 + h^2 \mu_k^2 e^{2 \phi (r)}  } 
\end{eqnarray*} 
which is smooth on $ H_{h,k}^{-1}(J) $ for any $ J \Subset (0,\infty) $ (hence in particular near  $ H_{h,k}^{-1} \big( \mbox{supp}( \varphi) \big) $). Choosing open intervals $ J_0 \Subset J_1 \Subset J_2 \Subset (0,+\infty) $ with $ J_0 $ containing the support of $ \varphi $, we obtain the analogue of Proposition \ref{optiquegeometriqueSchrodinger} for $ \nu = 1/2 $ with a phase $ S_{h,k}^{(1/2)} $ solution to
$$ \partial_t S_{h,k}^{(1/2)} = H_{h,k}^{1/2}(r,\partial_r S_{h,k}^{(1/2)}), \qquad S_{h,k}^{(1/2)}(0,r,\rho) = r \rho , $$
 The reduction to this case can be seen either by observing that $ \sqrt{h^2 {\mathfrak p}_k}  $ is well approximated, near the region where one is microlocalized, by a pseudo-differential operator with $ H_{h,k}^{1/2} $ as principal symbol, which is well known. Or, to avoid this approximation step, one may also consider directly $ \cos ( t \sqrt{{\mathfrak p}_k}) $ (which turns out to be sufficient here) and prove an approximation of $ \cos ( t \sqrt{{\mathfrak p}_k}) A_k (h) $ as the sum of two Fourier integral operators  associated respectively to the phases $ S_{h,k}^{(1/2)} (t) $ and $ S_{h,k}^{(1/2)}(-t) $.

Proceeding as in (\ref{refonde1})-(\ref{borneinferieurephaseexplicite}) with $ H = H_{h,k}^{1/2} $ and using the crucial observation that
\begin{eqnarray*}
\partial_{\rho \rho}^2 H_{h,k}^{1/2} & = & \frac{1}{\big( H_{h,k}^{1/2} \big)^3} h^2 \mu_k^2 e^{2 \phi (r)} \\
\partial_{r \rho}^2 H_{h,k}^{1/2} & = & - \frac{\rho \phi^{\prime}(r)}{\big( H_{h,k}^{1/2} \big)^3} h^2 \mu_k^2 e^{2 \phi (r)} ,
\end{eqnarray*}
we see that if $ t_0 $ is small enough,  we have
$$  \partial_{\rho}^2 S_{h,k}^{(1/2)}(t,r,\rho) \geq C |t| h^{2} \mu_k^2 e^{2 \phi(L)} , \qquad |t| < t_0 $$
uniformly in $ h, k  , L $ and $ (r,\rho) \in H_{h,k}^{-1}(J_1)  $ such that $ r \in \mbox{supp} ({\mathds 1}_L) $. Using Proposition \ref{VanderCorput}, we get a dispersion estimate for the Fourier integral operator of order $ h^{-1} (h/|t|h^{2} \mu_ k^2 e^{2\phi(L)})^{1/2} $. Using the analogue of (\ref{ecartcontrole}) for $ \nu = 1/2 $ and the fact that $ h^N \lesssim h^{N-1} \mu_k^{-1} e^{-\phi(L)} $ by (\ref{regionclassiquementinterdite}), we get
$$ \left| \left|  \mathds{1}_L   A_k (h) e^{i (t_1-t_2) {\mathfrak p}_k ^{1/2}  }  A_k (h)^*  \mathds{1}_L \right| \right|_{L^1 \rightarrow L^{\infty}} \leq C e^{-\phi (L)} \mu_k^{-1} h^{-1} \left( \frac{h}{h^2 |t_1-t_2|} \right)^{1/2} . $$ 
Using (\ref{estimeeL2H}) and (\ref{pourephiaile}), the above estimate yields
\begin{prop} \label{clarificationwave} Let $ (p,q) $ be sharp wave-admissible in dimension 2 with $ p,q $ real. Consider
$$ q_1 = q , \qquad p_1 = p , \qquad \beta = 2, \qquad \sigma = \frac{1}{2} . $$
Then there exists $ C > 0 $ such (\ref{dispersionsemiclassique}) holds true with $ B_{h,L,k} = C $ and $ D_{h,L,k} = C \mu_k^{-1} $.
\end{prop}

\noindent {\bf Proof of Proposition \ref{reductionondesprop}.} Note first that since $p,q $ are real and sharp wave admissible we have $ p > 4 $ hence $ p_ 1 $ is greater than 2. By Propositions \ref{translationStrichartz} and \ref{clarificationwave}, we then have directly
$$ \mu_k^{\frac{1}{2}-\frac{1}{q}}\big| \big| \Psi_{h,L,k}^{(1/2)} \big| \big|_{L^p ([0,t_0],L^q (e^{-\phi(r)}dr))} \leq C h^{(\sigma-2) \left( \frac{1}{2} - \frac{1}{q} \right)} || u_k ||_{L^2}  $$ so we conclude, as for Schr\"odinger, by using (\ref{Strichartzl2}) and the fact that $ \sum_{k \geq 0} ||u_k ||_{L^2}^2 = || \psi_0 ||_{L^2_{G_0}}^2 $. Here there is no semiclassical rescaling of time. \finpreuve

\section{Counterexamples} \label{modezero}
\setcounter{equation}{0}

\subsection{Proof of Theorem \ref{theorem1}}

\noindent {\it Proof.} Let $ e_0  $ be an eigenfunction in  $  \mbox{Ker}( \Delta_{\mathcal A} ) $ and  fix a non zero $ u_0 \in C_0^{\infty} (\Ra^+) $. For $ n $ sufficiently large define 
$$ \psi_n (r,\alpha)  = e^{\frac{\phi(r)}{2}} u_0 (r - n) e_0 (\alpha) . $$
Using (\ref{theunitarymapping}) and Proposition \ref{separationdevariablesprop}, it is easy to check that $ \psi_n  $ belongs to the domain of all powers of $ \Delta_0 $ since $ u_0 (r - n) e_0 (\alpha) $ belongs to the domains of all powers of $ P $, which  follows from the fact that $ u_0 (r - n) $ belongs to the domains of all powers of $ \mathfrak{p}_0 $, which in turn is obvious. By the translation invariance of $ \partial_r^2 $ and the boundedness of $ w $ on $ \Ra $ (see (\ref{referenceaw}) for $w$), we see that for all $ k \in \Na $, 
$ \mathfrak{p}_0^{k}(  u_0 (\cdot - n) ) $ is bounded in $ L^2 ((r_0,\infty),dr) $ as $ n  $ grows. Thus, for all $ \sigma > 0 $, 
\begin{eqnarray}
  || \psi_n ||_{H_{G_0}^{\sigma}} \leq C_{\sigma} , \qquad n \gg 1 . \label{bornedanslesSobolev}
\end{eqnarray}
{\it Proof of 1.} We observe that, for $q > 2$,
\begin{eqnarray}
 || \psi_n ||_{L^q_{G_0}} = c_q \big| \big| e^{\left(\frac{1}{2} - \frac{1}{q} \right) \phi(\cdot)} u_0 (\cdot-n) \big| \big|_{L^q (\Ra)} \gtrsim e^{ \left( \frac{1}{2} - \frac{1}{q} \right)\phi(n)} \rightarrow + \infty, \qquad n \rightarrow \infty , \label{lowerboundstationnaire}
\end{eqnarray}
the lower bound being a consequence of \ref{order} and the support property of $ u_0 (r-n) $. Using that $ || \psi_n ||_{H^\sigma_{G_0}}^{-1} $ is bounded from below by (\ref{bornedanslesSobolev}), we get the result.

\medskip

\noindent {\it Proof of 2.} Using the same computation as above
\begin{eqnarray*}
 || \Psi_n (t) - \psi_n ||_{L^q_{G_0}} & = & \left| \left|e^{\left(\frac{1}{2}-\frac{1}{q} \right)\phi(\cdot)} (\cos t \sqrt{P_0} - I ) u_0 (\cdot - n) \otimes e_0 \right| \right|_{L^q ( (r_0,\infty) \times {\mathcal A},dr d{\mathcal A}) }  \\
 & \leq & C \left| \left|e^{\left(\frac{1}{2}-\frac{1}{q}\right) \phi (\cdot)} (\cos t \sqrt{\mathfrak{p}_0} - I ) u_0 (\cdot - n)  \right| \right|_{L^q ((r_0,+\infty),dr)}.  
\end{eqnarray*} 
By finite speed of propagation,  $ (\cos t \sqrt{\mathfrak{p}_0} - I ) u_0 (\cdot - n) $ is supported in $ [ n - C , n + C ] $ for $ |t| \leq 1 $ and some constant $ C $ independent of $n$. Therefore, on the support of this function, (\ref{order}) allows to use $ e^{\phi(r)} \lesssim e^{\phi(n)} $ and thus, using rough Sobolev embeddings, we get
\begin{eqnarray*}
 || \Psi_n (t) - \psi_n ||_{L^q_{G_0}} & \leq & C e^{\left(\frac{1}{2}-\frac{1}{q} \right) \phi (n)} || (\cos t \sqrt{\mathfrak{p}_0} - I ) u_0 (\cdot - L)  ||_{H^1_0 (r_0,\infty)} , \\
 & \leq &  C e^{\left(\frac{1}{2}-\frac{1}{q} \right) \phi (n)} || (\cos t \sqrt{\mathfrak{p}_0} - I ) (\mathfrak{p}_0 + 1)^{1/2} u_0 (\cdot - L)  ||_{L^2 ((r_0,\infty),dr)}
\end{eqnarray*} 
by (\ref{injectionSobolevexplicite}). Finally, by writing $ \cos t \sqrt{\mathfrak{p}_0} - I = - \int_0^t \sin s \sqrt{\mathfrak{p}_0} ds \sqrt{\mathfrak{p}_0} $, we  get
\begin{eqnarray*}
 || \Psi_n (t) - \psi_n ||_{L^q_{G_0}} & \leq & C |t| e^{\left(\frac{1}{2}-\frac{1}{q} \right) \phi (n)} \big| \big| \sqrt{\mathfrak{p}_0} (\mathfrak{p}_0 + 1)^{1/2} u_0 (\cdot - n)  \big| \big|_{L^2 ((r_0,\infty),dr) } \\
 & \lesssim & |t|  e^{\left(\frac{1}{2}-\frac{1}{q} \right) \phi (n)} ,
\end{eqnarray*}
since  $ \big| \big|  (\mathfrak{p}_0 + 1) u_0 (\cdot - n)  | |_{L^2  } $ is bounded in $ n $.
Therefore, using (\ref{lowerboundstationnaire}), we see that for some small enough $ t_0 > 0 $ and $c> 0$ we have
$$ || \Psi_n (t) ||_{L^q_{G_0}} \geq c e^{ \left( \frac{1}{2} - \frac{1}{q} \right) \phi (n)}, \qquad n \gg 1, \ \ |t|<t_0  .  $$
 Using this lower bound and the fact that that
$$ || \Psi_n ||_{L^p ( [0,1] , L^q_{G_0} ) } \geq t_0^{1/p} \inf_{(-t_0,t_0)} || \Psi_n (t) ||_{L^q_{G_0}} , $$
the result follows as in the previous item.

\medskip





 
\noindent {\it Proof of 3.} Here we consider $ \phi (r) = r $ and $r_0 = 0$. In this case $ w (r) = 1/4 $ is a constant so that we are reduced to a free Schr\"odinger equation on the half line. We consider
\begin{eqnarray}
u_n (r) & = & \exp \left( - \frac{(r-n)^2}{2} \right) - \exp \left( - \frac{(r+n)^2}{2}\right), \nonumber \\
  \psi_n (r,\alpha) & = &  e^{\frac{r}{2}} u_n (r) e_0 (\alpha) , \nonumber
\end{eqnarray}
Obviously, $ u_n $ is an odd Schwartz function on $ \Ra $ hence so is $ \partial_r^{2k} u_n $. This implies that the restriction of $ u_n $ to $ \Ra^+ $ satisfies the Dirichlet condition at $ r= 0 $ as well as $ {\mathfrak p}_0^k u_n $ for all $ k $. This implies that $ u_n $ belongs to the domains of all powers of $ {\mathfrak p}_0 $, with uniform bounds in $n$ so that we still have an upper bound of the form (\ref{bornedanslesSobolev}). On the other hand, a direct computation shows that
\begin{eqnarray*}
 \Psi_n (t,r,\alpha) & = & e^{\frac{r}{2}} U_n (t,r) e_0 (\alpha) , 
\end{eqnarray*}
with
\begin{eqnarray*} 
   U_n (t,r) & = & \frac{e^{i \frac{t}{4}}}{\sqrt{1-2it}} \left[ \exp \left( - \frac{1}{1 - 2 it} \frac{(r-n)^2}{2} \right) - \exp \left( - \frac{1}{1 - 2 it} \frac{(-r-n)^2}{2} \right) \right] .
\end{eqnarray*}
%
Now, we observe independently that given two real numbers $ \varepsilon \geq 0 $ and $ \tau > 0 $, one has
$$ \int_0^{ + \infty} \left| e^{\varepsilon r} \exp \left( - \tau \frac{(\pm r -n)^2}{2} \right)  \right|^q d r = e^{\pm q \varepsilon n + \frac{q \varepsilon^2}{2 \tau}} \int_{\mp n - \frac{\varepsilon}{\tau} }^{+ \infty} e^{- q \frac{\tau}{2}  x^2 } dx . $$
In particular, for fixed $ \varepsilon > 0 $ and $q > 0 $, there exists $ C > 0 $ such that, for all $ n \gg 1 $ and $ \tau \in (1/2,3/2) $ we have
$$  \int_0^{+ \infty} \left| e^{\varepsilon r} \exp \left( - \tau \frac{( r -n)^2}{2} \right)  \right|^q d r \geq C e^{ q \varepsilon n }   $$
and
$$ \int_0^{+ \infty} \left| e^{\varepsilon r} \exp \left( - \tau \frac{(- r -n)^2}{2} \right)  \right|^q d r \leq C e^{- q \varepsilon n } . $$
Using these estimates with $ \varepsilon = \frac{1}{2} - \frac{1}{q} $, $ q>2 $, and $ \tau = \frac{1}{1+4t^2} $, it is not hard to see that for some $ t_0 > 0 $ and $ c > 0 $ small enough we have
$$ || \Psi_n (t) ||_{L^q_{G_0}} \geq c e^{\varepsilon n}, \qquad n \gg 1, \ \ t \in (-t_0,t_0) .  $$
With this estimate at hand, we complete the proof as in the end of the proof of the item 2. \finpreuve

\bigskip

\noindent {\bf Remark.} The counterexamples of this section are not specific to surfaces. They would work exactly the same in higher dimension, {\it i.e.} with $ {\mathcal A} $ of dimension $ n- 1 \geq 2 $, up to obvious  natural modifications such as the replacement of $ \frac{1}{2} - \frac{1}{q} $ by $ \frac{n-1}{2} - \frac{n-1}{q} $.


\subsection{Proof of Theorem \ref{theorem5}}
For notational simplicity we assume that $ k_0 = 1 $ and thus consider  the first non zero eigenvalue $ \mu_1^2 $ of $ - \Delta_{\mathcal A} $ and an associated eigenfunction $ e_1 (\alpha) $. We consider 
 $$ \psi_0^h (r,\alpha) = e^{\frac{r}{2}} u_0^h (r) e_{1} (\alpha) , $$
 where, for a given $ \chi \in C_0^{\infty}(\Ra) $ which is equal to $ 1 $ near $ 0 $, we have set
$$ u_0^h (r) = (\pi h)^{-1/4} \chi (r + \log h) \exp \left( \frac{-(r+ \log h)^2}{2 h} \right) . $$
Ovviously $ \psi_0^h $ belongs to the range of $ \Pi^c $ and so does
\begin{eqnarray}
  e^{it \Delta_0} \psi_0^h = e^{\frac{r}{2}}  e^{- it {\mathfrak p}_1} u_0^h \otimes e_1 , \label{formulereecrite}
\end{eqnarray}  
where we recall that $ {\mathfrak p}_1 = - \partial_r^2 + \mu_1^2 e^{2r} + \frac{1}{4} $. Up to the cutoff $  \chi (r + \log h) $, which ensures that $u_0^h $  satisfies the Dirichlet condition,
we can  interpret $ u_0^h $ as a wave packet (or coherent state) centered at $ (- \log h , 0) \in T^* (r_0,\infty) $.  By rescaling the time as $ t = h s $, we will compute
$$  e^{- it {\mathfrak p}_1} u_0^h = e^{- \frac{i}{h}s  h^2{\mathfrak p}_1} u_0^h , $$
by seeing $ h^2{\mathfrak p}_1 $ as a semiclassical Schr\"odinger operator with ($h$ dependent) principal symbol
$$ H_h (r,\rho) = \rho^2 + \mu_1^2 h^2 e^{2r} . $$
The propagation of coherent states by the semiclassical Schr\"odinger equation is a well known topic. We recall here the main points of the analysis.
We let $ \Phi_h^s = \big(x_s^h , \xi_s^h \big)  $ be the Hamiltonian flow of $ H_h $. The classical action associated to a trajectory starting at $ (x,\xi) $ at $s=0$ is 
$$ S_s^h = S_s^h (x,\xi) = \int_0^s \dot{x}^h_{\tau} \xi_{\tau}^h - H_h (x_{\tau}^h,\xi_{\tau}^h) d \tau . $$ 
We also consider
$$ a_s^h = \frac{\partial x_s^h}{\partial x}, \qquad b_s^h = \frac{\partial x_s^h}{\partial \xi}, \qquad c_s^h = \frac{\partial \xi_s^h}{\partial x}, \qquad d_s^h = \frac{\partial \xi_s^h}{\partial \xi} . $$ 
and let
$$ \Gamma_s^h = \frac{c_s^h + i d_s^h}{a_s^h + i b_s^h} . $$
Then, according to a well known procedure (see {\it e.g.}  \cite{Robe2,CoRo} for a detailed presentation), we can write for any fixed $ N $,
\begin{eqnarray}
 e^{ - i s h {\mathfrak p}_1} u_0^h  =   U_N^h (s,r)    + h^N R_N^h (s,r) , \label{approximationannoncee}
\end{eqnarray} 
with
$$ U_N^h (s,r) = (\pi h)^{-1/4} A_N^h (s,r) \exp \frac{i}{h} \left(  S_s^h + \xi_h^s (r-x_h^s)  +  \frac{\Gamma_s^h}{2}(r-x_h^s)^2 \right) $$ 
where the trajectory $ (x_h^s,\xi_s^h) $, the action $S_s^h $ and $  \Gamma_s^h$ are associated to the starting point  $ (- \log h , 0) $ at $s=0$, and where the amplitude is of the form
$$ A_N^h (s,r) =  \big(a_s^h + i b_s^h \big)^{-1/2} \chi \big(r - x_s^h \big) \left(1 + h^{1/2} \sum_{k=0}^{k(N)} h^{k/2} Q^h_k (s,y) \right)_{| y = \frac{r-x_s^h}{h^{1/2}}}  $$
for some $ k (N) $ large enough and some suitable functions $ Q_k^h(s,y) $ which are polynomials in $y$ and vanish at $s=0$. We summarize this analysis as a proposition.
\begin{prop} \label{etatcoherentparticulier} If $ h_0$ is small enough and  $ N \geq 0$ is arbitrarily fixed, the approximation (\ref{approximationannoncee}) holds with a remainder such that, for all $ j \leq N $,
\begin{eqnarray}
  || {\mathfrak p}_1^j R_N^h (s) ||_{L^2} \leq C_j h^{-2j} , \qquad |s| \leq 1 , \ h \in (0,h_0] . \label{decroissancedomaine}
\end{eqnarray}
The expansion $ U_N^h $ is such that
$$ Q_k^h (s,y) = \sum_{l = 0}^{3(k+1)} c_{kl}^h (s) y^l, $$
with coefficients such that
$$ |c_{kl}^h(s)| \leq C_{kl}   \qquad \mbox{for} \ |s| \leq 1, \ h \in (0,h_0] . $$
Furthermore, there exists a  constant $ C >1 $ such that, for $ |s| \leq 1 $ and $ h \in (0,h_0 ] $,
\begin{eqnarray}
 C^{-1} \leq  \emph{Re} (i \Gamma_s^h) \leq C  , \qquad |x_h^s + \log h| \leq C . \label{bornesurGamma}
\end{eqnarray} 
\end{prop}



In the next proposition, we check that $ \psi_0^h $ is approximately spectrally localized.

\begin{prop} \label{controleerreur1} Let $ \varphi \in C_0^{\infty}(\Ra) $ be such that $ \varphi \equiv 1 $ near $ \mu_1^2 $. Then, for any $ N \geq 0$,
\begin{eqnarray}
 \big| \big| \varphi (-h^2 \Delta_0) \psi_0^h - \psi_0^h \big| \big|_{H^{2}_{G_0}} \leq C_{N} h^N . \label{propositionacitergaussienne}
 \end{eqnarray}
Moreover, for any $ \sigma \in [0,2] $,
\begin{eqnarray}
 || \psi_0^h ||_{H_{G_0}^{\sigma}} \lesssim h^{-\sigma} . \label{bornesuperieureSobolev}
\end{eqnarray}
\end{prop}

\noindent {\it Proof.} We recall first  that the asymptotic expansion of the action of a pseudo-differential operator on a wave packet is well known (see \cite{Robe2}), and basically given by linear combinations of the symbol and its derivatives evaluated at the center  (here $ (- \log (h),0) $) times polynomials  and derivatives of the wave packet. Thus by Proposition \ref{calculfonctionnelreecrit} and the fact  that $ \varphi (\rho^2 + h^2 \mu_1^2 e^{2r}) $ is equal to $1$ near  $ (- \log (h),0) $, all terms of the expansion  of $ \varphi (-h^2 \Delta_0) \psi_0^h - \psi_0^h $ vanish which yields easily (\ref{propositionacitergaussienne}). We omit the calculations but point out that the dependence of the center $ (- \log (h),0) $ on $h$ does not cause any problem.
To get (\ref{bornesuperieureSobolev}) we observe that $ | | \psi_0^h | |_{L^2_{G_0}} \lesssim 1 $ 
and that $ \Delta_0  \psi_0^h = e^{\frac{r}{2}} {\mathfrak p}_1 u_0^h \otimes e_1 $ with   $ || {\mathfrak p}_1 u_0^h ||_{L^2(\Ra^+)} \lesssim h^{-2} $. Then a simple interpolation argument yields the result. \finpreuve

\bigskip

\noindent {\bf Proof of Theorem \ref{theorem5}.} From now on, we consider a  Schr\"odinger admissible pair $ (p,q) $ in dimension $2$.
Thanks to Proposition \ref{controleerreur1} and Corollary \ref{corollaireSchrodinger}, we see that
$$ \big| \big| \Pi^c {\mathds 1}_{[r_1,\infty)} e^{it\Delta_0} \psi_0^h - \Pi^c {\mathds 1}_{[r_1,\infty)} e^{it \Delta_0} \varphi (-h^2 \Delta_0) \psi_0^ h \big| \big|_{L^p ([0,1],L^q_{G_0})} \leq C_N h^N , $$
which is {\it a fortiori} true if we restrict the time interval to $ [0,h] $.
Using this error estimate together with the upper bound (\ref{bornesuperieureSobolev}), we will get the expected counterexample if we show that
\begin{eqnarray}
 \big| \big|  {\mathds 1}_{[r_1,\infty )} e^{it\Delta_0} \psi_0^h  \big| \big|_{L^p( [0,h]_t,L^q_{G_0} )} \gtrsim  h^{-\frac{1}{2p}} ,
 \nonumber
\end{eqnarray}
since assuming $ \sigma < \sigma_{\rm S} $ means $ \sigma < \frac{1}{2p} $. Equivalently, in semiclassical time scaling, the above lower bound reads
\begin{eqnarray} 
 \big| \big|  {\mathds 1}_{[r_1,\infty )} e^{is h \Delta_0} \psi_0^h  \big| \big|_{L^p( [0,1]_s,L^q_{G_0} )} \gtrsim  h^{-\frac{3}{2p}} . \label{borneinferieureducontrexemplefinal}
\end{eqnarray} 
Let us prove this lower bound. By (\ref{formulereecrite}), (\ref{approximationannoncee}) , we have
$$ \big( e^{is h \Delta_0} \psi_0^h \big) (r,\alpha) = e^{r/2} \left( U_N^h (s,r) + h^N R_N^h (s,r) \right) e_1 (\alpha) $$
By Propositions \ref{aprioriL1} and  \ref{etatcoherentparticulier}, we have for some $ N_0 > 0 $,
\begin{eqnarray}
 \big| \big| h^N e^{\cdot /2} R_N^h (s, \cdot)  \otimes e_1  \big| \big|_{L^q_{G_0}} \lesssim h^{N-N_0} , \label{derniereerreurestimee}
\end{eqnarray} 
uniformly with respect to $ |s| \leq 1 $. In particular, if we choose $ N $ large enough, this is a bounded quantity. On the other hand, using that $ U_N^h (s,.) $ is supported on a set where $ r - x_h^s  $ is bounded and that $ x_h^s = - \log (h) + O (1) $ by the second estimate of (\ref{bornesurGamma}), we have
\begin{eqnarray*}
 \big| \big| {\mathds 1}_{[r_1,\infty)} e^{\cdot /2}  U_N^h (s,\cdot) \otimes e_1 \big| \big|_{L^q_{G_0}}  \gtrsim  e^{-  \left(\frac{1}{2} - \frac{1}{q} \right) \log h} || U_N^h (s,.) ||_{L^q (\Ra)} 
\end{eqnarray*}
Using the form of $ U_N^h (s,h) $ together with the first estimate of (\ref{bornesurGamma}), we obtain
\begin{eqnarray*} 
\big| \big| {\mathds 1}_{[r_1,\infty)} e^{\cdot /2}  U_N^h (s,\cdot) \otimes e_1 \big| \big|_{L^q_{G_0}}   \gtrsim  h^{- \left( \frac{1}{2} - \frac{1}{q} \right)} h^{-\frac{1}{4} + \frac{1}{2q}} = h^{-\frac{3}{2p}}  
\end{eqnarray*} 
This estimate and (\ref{derniereerreurestimee}) imply (\ref{borneinferieureducontrexemplefinal}) which completes the proof. \finpreuve


\begin{thebibliography}{99}


\bibitem{Anton} {\sc R. Anton}, {\it Strichartz inequalities for Lipschitz metrics on manifolds and nonlinear Schr\"odinger equation on domains},
Bull. SMF 136, fascicule 1, 27-65 (2008) 

\bibitem{BaCh} {\sc H. Bahouri, J.Y. Chemin}, {\it Equations d'ondes quasilin\'eaires et estimations de Strichartz},
Amer. J.  Math.  121, no. 6, 1337-1377 (1999)

\bibitem{BanicaDuyckaerts} {\sc V. Banica, T. Duyckaerts}, {\it Weighted Strichartz estimates for radial Schr\"odinger equation on noncompact manifolds},
  Dyn. P.D.E. 4, no. 4, 335-359 (2007)

\bibitem{BFHM} {\sc M. Blair, G.A. Ford, S. Herr, J.L. Marzuola}, {\it Strichartz estimates for the Schr\"odinger equation on polygonal domains} (with G. A. Ford, S. Herr, and J. L. Marzuola), Journal of Geometric Analysis, 22 (2), 339-351 (2012)

\bibitem{BSS1}  {\sc M. Blair, H. Smith, C. Sogge}, {\it Strichartz estimates for Schr\"odinger operators in compact manifolds with boundary}, Proc. AMS, 136 (1),  247-256 (2008)


\bibitem{BSS2} {\sc M. Blair, H. Smith, C. Sogge}, {\it Strichartz estimates and the nonlinear Schr\"odinger equation on manifolds with boundary}, Mathematische Annalen, 354 (4), 1397-1430 (2012) 

\bibitem{BoucletLP} {\sc J.-M. Bouclet}, {\it Littlewood-Paley decompositions on manifolds with ends}, Bull. SMF 138, fascicule 1, 1-37 (2010)

\bibitem{Boucletref} {\sc \name}, {\it Strichartz estimates on asymptotically hyperbolic manifolds},
Analysis \& PDE 4-1, 1-84 (2011)

\bibitem{Bourgain} {\sc J. Bourgain}, {\it Fourier transform restriction phenomena for certain lattice subsets and application
to nonlinear evolution equations I. Schr\"odinger equations}, GAFA 3, 107-156 (1993)

\bibitem{BGT} {\sc N. Burq, P. G\'erard, N. Tzvetkov}, {\it Strichartz inequalities and the nonlinear Schr\"odinger equation on compact manifolds}, Amer. J. Math. 126, no. 3, 569-605 (2004)

\bibitem{BGH} {\sc N. Burq, C. Guillarmou, A. Hassell}, {\it Strichartz estimates without loss on manifolds with hyperbolic trapped geodesics}, GAFA  20, no. 3, 627-656 (2010)

\bibitem{CoRo} {\sc M. Combescure, D. Robert}, {\it Coherent states and applications in mathematical physics}
Theoretical and Mathematical Physics. Springer, Dordrecht (2012)

\bibitem{DiSj1} {\sc M. Dimassi, J. Sj\"ostrand},
{\it Spectral asymptotics in the semiclassical limit},  London
Mathematical Society Lecture Note Series, 268. Cambridge
University Press (1999)


\bibitem{HaZh} {\sc A. Hassell, J. Zhang}, {\it Global-in-time Strichartz estimates on non-trapping asymptotically conic manifolds}, arXiv:1310.0909
   

\bibitem{Hebey} {\sc E. Hebey}, {\it Nonlinear analysis on manifolds: Sobolev spaces and inequalities},
Courant Lecture Notes in Mathematics,  AMS, Providence, RI (1999)

\bibitem{OIvan} {\sc O. Ivanovici}, {\it On the Schr\"odinger equation outside strictly convex domains}, Analysis \& PDE vol.3, no.3, 261-293 (2010)

\bibitem{Ivanovici} {\sc \name}, {\it Counter-examples to the Strichartz estimates for the wave equation in general domains with boundary}, JEMS, vol. 14, issue 5,   1357-1388 (2012)

\bibitem{ILP} {\sc O. Ivanovici, G. Lebeau, F. Planchon}, {\it Dispersion for the wave equation inside strictly convex domains I: the Friedlander model case},
 Annals of Mathematics, to appear.

\bibitem{Kapi}  {\sc L.V. Kapitanski}, {\it Some generalizations of the Strichartz-Brenner inequality}, Algebra i Analiz 1:3, 127-159 (1989) 

\bibitem{KeelTao} {\sc M. Keel, T. Tao}, {\it Endpoint Strichartz estimates}, Amer. J. Math. 120, no. 5, 955-980 (1998)

\bibitem{MMTa} {\sc J. Marzuola, J. Metcalfe, D. Tataru}, {\it Strichartz estimates and local smoothing estimates for asymptotically flat Schr\"odinger equations},
J. Funct. Anal. 255, no. 6, 1497-1553 (2008) 

\bibitem{MetcTata} {\sc J. Metcalfe, D. Tataru}, {\it Global parametrices and dispersive estimates for variable coefficient wave equations}, Math. Ann. 353, no. 4, 1183-1237 (2012)

\bibitem{MuscSchl} {\sc C. Muscalu, W. Schlag}, {\it Classical and Multilinear Harmonic Analysis. Vol. 1},
 Cambridge Studies in Advanced Mathematics (2013)
 
\bibitem{Robe1} {\sc D. Robert},
{\it Autour de l'approximation semi-classique}, Progress in
mathematics, {\bf 68}, Birkha\"user (1987)

\bibitem{Robe2} {\sc \name}, {\it Propagation of coherent states in quantum mechanics and applications}, Partial differential equations and applications, 181-252,
S\'emin. Congr., 15, Soc. Math. France (2007)

\bibitem{SmTa} {\sc H.F. Smith, D. Tataru}, {\it Sharp counterexamples for Strichartz estimates for low regularity metrics}, Math. Res. Lett. 9, no. 2-3, 199-204 (2002)

\bibitem{Stein} {\sc E. Stein}, {\it Singular integrals and differentiability of functions}, Princeton Univ. Press (1970)

\bibitem{Stein2} {\sc \name}, {\it Harmonic analysis}, Princeton Univ. Press (1993)



\bibitem{TaII} {\sc D. Tataru}, 
{\it Strichartz estimates for second order hyperbolic operators with nonsmooth coefficients. II},
Amer. J. Math. 123, no. 3, 385-423 (2001)

\bibitem{TaIII} {\sc \name}, {\it Strichartz estimates for second order hyperbolic operators with nonsmooth coefficients. III}, J. Amer. Math. Soc. 15, no. 2, 419-442 (2002)


\bibitem{Zhang} {\sc J. Zhang}, {\it Strichartz estimates and nonlinear wave equation on nontrapping asymptotically conic manifolds},  arXiv:1310.4564
  


\end{thebibliography}
\end{document}